\input amstex
\documentstyle{amsppt}
\NoBlackBoxes
\nopagenumbers
\magnification=1200
\define \a{\alpha}
\define \be{\beta}
\define \Dl{\Delta}
\define \dl{\delta}

\define \g{\gamma}
\define \G{\Gamma}

\define \k{\kappa}
\define \lm{\lambda}
\define \Lm{\Lambda}
\define \Om{\Omega}
\define \om{\omega}
\define \r{\rho}
\define \s{\sigma}
\define \Si{\Sigma}
\define \th{\theta}

\define \ve{\varepsilon}
\define \vp{\varphi}

\define \z{\zeta}

\define \cd{\cdot}

\define \fa{\forall}

\define \iy{\infty}

\define \p{\partial}
\define \sm{\setminus}

\define \Llr{\Longleftrightarrow}

\define \sbt{\subset}

\define \ep{\endproclaim}

\define \bk{\bigskip}

\define \hf{\hat{f}}
\define \hg{\hat{g}}

\define \BT{\Bbb T}
\define \BC{\Bbb C}

\define \BN{\Bbb N}

\define \BR{\Bbb R}
\define \BZ{\Bbb Z}

\define \Cp{\Cal P}

\define \Lip{\operatorname{Lip}}
\define \Lc{Lebesgue constants} 
\define \Ft{Fourier transform}
\define \Fs{Fourier series}
\define \Fc{Fourier coefficient}
\define \Il{\int\limits}
\define \Sl{\sum\limits}
\define \Ll{\lim\limits}
\define \spl{\sum\limits}
\define \inl{\inf\limits}

\define \mes{\operatorname{mes}}
\define \supp{\operatorname{supp}}
\define \sign{\operatorname{sign}}
\define \esp{\operatornamewithlimits{ess\,sup}}

\document
\font\sssb=cmssbx10 scaled\magstep2
\centerline{\bf Donetsk State University}
\vskip 3cm
\centerline{\sssb R. M. TRIGUB}
\vskip 2cm
\centerline{\sssb SOME TOPICS IN FOURIER ANALYSIS}
\centerline{\sssb AND APPROXIMATION THEORY}
\vskip 3cm
\centerline{\bf PREPRINT}
\vskip 4cm
\centerline{\bf Donetsk\qquad $\bullet$\qquad 1995}
\vfill\eject
\centerline{\bf CONTENTS}
\vskip 1cm
Preface

1. The Fourier transform. Absolute convergence

2. Multipliers and comparison of various methods 

\qquad of summability of Fourier series as a whole

3. Regularity of summubility methods of Fourier series

4. Lebesgue constants and approximation of function classes

5. Two-sided estimates of approximation. Moduli of smoothness

\qquad and K-functionals

6. Hardy spaces $H_p$

7. Positive definite functions and splines

8. Multiplicative Walsh system

9. Pointwise approximation of functions by polynomials.

\qquad Approximation by polynomials with integral coefficients

10. Approximation by trigonometric polynomials with 

\qquad a given number of harmonics

11. The Pompeiu problem

12. Entire functions of exponential type and the Fourier integral

13. Unsolved problems

\quad

References
\vfill\eject
\head Preface.\endhead
This manuscript presents shortly the results obtained by participants
of the scientific seminar which is held more than twenty years
under leadership of the author at Donetsk University. In the list of
references main publications are given. These results are published
in serious scientific journals and reported at various conferences,
including international ones at Moscow, ICM66; Kaluga, 1975; Kiev, 1983;
Haifa, 1994; Z\"urich, ICM94; Moscow, 1995. As for the area of
investigation this is the Fourier analysis and the theory of approximation
of functions. Used are methods of classical analysis including
special functions, Banach spaces, etc., of harmonic analysis in
finitedimensional Euclidean space, of Diophantine analysis, of random
choice, etc.

The results due to the author and active participants of the seminar,
namely E. S. Belinskii, O. I. Kuznetsova, E. R. Liflyand,
Yu. L. Nosenko, V. A. Glukhov, V. P. Zastavny, Val. V. Volchkov,
V. O. Leontyev, and others, are given. Besides the participants of the
seminar and other mathematicians from Donetsk, many mathematicians
from other places were speakers at the seminar, in particular,
A. A. Privalov, Z. A. Chanturia, Yu. A. Brudnyi, N. Ya. Krugljak,
V. N. Temlyakov, B. D. Kotlyar, A. N. Podkorytov, M. A. Skopina,
A. A. Ligun, A. S. Romanyuk, V. A. Martirosyan.

Besides the papers of the participants of the seminar, only monographs and
survey papers are given in the list of references. Relative
results of other mathematicians are described in the text without
details, but the year of publication is given.

The outline of the work is seen from the contents. Some unsolved problems
are given as well.

Main notation is given in the very beginning of \S1. Bibliographical remarks
are given in the end of every section.

The most part of this preprint will appear soon in extended form
in the author's book "Fourier Analysis and Approximation Theory",
World Federation Publishers, Inc., P.O.Box 48654, Atlanta, GA 30362-0654, USA.

The author thanks E. R. Liflyand (Bar-Ilan University, Israel) for the
translation from Russian to English.
\vfill\eject

\head 1. The Fourier transform. Absolute convergence.\endhead
We denote absolute positive constants by the letter $C.$ Let
$\g(\a,\be)$ be a positive constant depending only on $\a$ and $\be.$
Let $\BR^m$ be the $m$-dimensional real Euclidean space, $\BT^m=
[-\pi,\pi)^m$ be the torus, and $\BZ^m$ be the lattice of points
in $\BR^m$ with integer coordinates.
For $x=(x_1,...,x_m),$ $y=(y_1,...,y_m)\in\BR^m$ denote by $(x,y)
=\sum\limits_{j=1}^m x_j  y_j$ the inner product of $x$ and
$y,$ by $|x|=(x,x)^{1/2}$ the Euclidean norm of $x.$ We write
$x\le y$ when $x_j\le y_j$ for all $j\in[1,m].$ Denote also
$\BR^m_+=\{x\in\BR^m: x\ge 0\},$ and $||\cd||_p$ means the $L_p$-norm.
We omit superscript $m$ for $m=1.$ The abbreviation {\it a.e.} means
{\it almost everywhere} with respect to the Lebesgue measure.

For every function $f\in L_p(\BR^m),$ with $p\in[1,2],$ the Fourier transform
$$\hat f(x)=(2\pi)^{-{m\over 2}}\int\limits_{\BR^m}f(u)e^{-i(x,u)}\,du$$
is defined as usual (see e.g., [M22]). Let
$$B(\BR^m)=\{f: f(x)=\int\limits_{\BR^m}e^{-i(x,u)}\,d\mu(u),\qquad
||f||_B=\operatorname{var}\mu<\infty\}$$
be the space of the Fourier transforms of all the finite complex valued
Borel measures on $\BR^m.$ If also a measure $\mu$ is absolutely continuous
with respect to Lebesgue measure, that is $d\mu(u)=g(u)du,$ where
$g\in L_1(\BR^m),$ then we shall write $f\in A(\BR^m)$ and $||f||_A=||g||_1.$

For $f: \BT^m\to\Bbb C$ set for $p>0$
$$A_p(\BT^m)=\{f: f(x)=\sum\limits_{k\in\BZ^m}c_ke^{i(k,x)},\qquad
||f||_{A_p}=\biggl(\sum\limits_k|c_k|^p\biggr)^{{1\over p}}<\infty\}.$$
Notice that $||\cd||_{A_p}$ is the quasinorm for $p\in(0,1).$ Set also
for $p\in(0,1)$
$$A_p(\BR^m)=\{f: f=\hat g,\qquad ||f||_{A_p}=||g||_p<\infty\}$$
in assumption, e.g., $g\in L_p\cap L_1(\BR^m).$ In particular, for a function 
with compact support $f\in A_p$ is equivalent to $\hat f\in L_p(\BR^m).$
It is well known that $B(\BR^m),$ $A(\BR^m)$ and $A(\BT^m)$ are Banach
algebras, and their local structure is the same (see e.g. [M10]).
A difference in behavior of functions from $B(\BR^m)$ and $A(\BR^m)$
occurs only near infinity.
\flushpar
\proclaim{1.1} If $f\in B(\BR^m)$ and is of bounded Vitali variation
outside some neighborhood of the origin,
and $\Ll_{|x|\to\infty}f(x)=0,$ then $f\in A(\BR^m)$ and $||f||_B=
(2\pi)^{-m}||\hf||_1.$ \ep   \bk
The following two statements are connected with discretization
and we shall give them only for $m=1$ (for $p=1$ the case {\bf 1.2}b) 
is Wiener's result).
\flushpar
\proclaim{1.2} a) If $f\in L_1(\BR)$ and supported in $[-\pi,\pi],$ then
for every $p>0$ $$\hf\in L_p(\BR)\Llr f\quad\text{and}\quad f_1\in A_p(\BT)$$
where $f_1(x)=e^{ix/2}f(x).$

b) If either $p>1,$ or $p\in(0,1]$ and diameter of the support of $f$ is less 
than $2\pi,$ then for $p\in(0,1]$ $$\hf\in L^p(\BR)\Llr f \in A_p(\BT).$$
\ep\bk
\proclaim{1.3} Let $n\in\BZ$ and for some $r\in\BZ_+$ the functions
$f$ and $f^{(r)}$ be of bounded variation on $[n,+\infty),$ and
$\lim f^{(\nu)}(x)=0$ as $x\to+\infty$ and $\nu\in[0,r].$ Then for
$0<|x|\le\pi$ $$\Sl_{k=n}^\infty f(k)e^{ikx}=\Il_n^\infty f(u)e^{iux}\,du
+{1\over 2}f(n)e^{inx}+e^{inx}\Sl_{p=o}^{r-1}{(-i)^{p+1}\over p!}
h^{(p)}(x)f^{(p)}(n)+{\th\over\pi^r}V_n^\infty f^{(r)},$$
where $V_n^\infty f^{(r)}$ is the total variation of $f^{(r)}$ on
$[n,\infty),$ the function $h(x)={1\over x}-{1\over 2}\cot{x\over 2},$
and $|\th|\le 3.$\ep
Tending to the limit as $x\to 0$ we get the classical Euler-Maclaurin
formula as a corollary.

Let us go on to necessary and sufficient conditions of belonging to
the algebras defined above.
\flushpar
\proclaim{1.4} Let $E$ be a closed set in $\BR^m$ and $A(E)$ be a linear
set of functions $E\to\BC,$ which is containing restrictions of all
boundedly supported functions from $C^{\iy}(\BR^m)$ as well as
those multiplied by elements of $A(E).$ Then $A(E)$ possesses a local
property, namely if for each $x\in E,$ including the infinity point when
$E$ is non-compact, there exist a neighborhood $V_x$ of $x$ and a
function $g_x\in A(E)$ such that $f=g_x$ on $V_x,$ then $f\in A(E).$\ep
\flushpar
\proclaim{1.5} Let $f\in A(\BR^m)$ and $f_0$ be its radial part, that is
$f_0(t)$ is an integral average of $f$ over the sphere $|x|=t.$
Then $f_0\in C^{m_1}(0,\iy),$ where $m_1=\biggl[{m-1\over 2}\biggr],$
when $t\to\iy$ we have $\lim t^p f_0^{(p)}(t)=0$ for $0\le p\le m_1.$
Besides that, for all $t>0$ the following integral converges:
$$\Il_0^t u^{m_1-{m+1\over 2}}\biggl[f_0^{(m_1)}(t+u)-f_0^{(m_1)}(t-u)
\biggr]\,du.$$ \ep 

The following (approximate) criterion is similar for
$p=1$ to that for the absolute convergence of orthogonal series
due to S. B. Stechkin (1955; see e.g., [M10]).

Taking $f\in L_2(\BR^m),$ let us assume that for each $\s>0$
$$a_{\s}(f)_2=\inf\{||f-g||_2\quad\text{for}\quad g\quad\text
{such that}\quad\mes\supp\hg\le\s\},$$ where $\mes$ denotes the
Lebesgue measure, be best approximation to $f$ in $L_2(\BR^m)$
by functions with the spectrum in a set of measure less or equal
to $\s.$
\flushpar
\proclaim{1.6} If $f\in L_p(\BR^m),$ then for $p\in(0,2)$ we have
$\hf\in L_p(\BR^m)$ if and only if
$$\Il_0^{\iy}\s^{-p/2}\,(a_{\s}(f)_2)^p\,d\s<\iy.$$\ep
\flushpar
\proclaim{1.7} Let $f\in C(\BR^m),$ $\lim f(x)=0$ as $|x|\to\iy,$ and
for $r\in\BR$ $$\Om_r(h_1,...,h_m)=\Om_r(f;h)=||\dot\Dl_h^r f(\cdot)||_2,$$
where $$\dot\Dl_h^r f(x)=\biggl(\prod\limits_{j=1}^m\dot\Dl_{h_j}^r\biggr)f(x),
\qquad\dot\Dl_{h_j}f(x)=f(x+h_je_j^0)-f(x-h_je_j^0)$$ is the symmetric
mixed difference, $e_j^0$ is the orth of the axis $Ox_j.$

a) If for at least one $r$
$$\Sl_{s_1=-\iy}^\iy...\Sl_{s_m=-\iy}^\iy 2^{{1\over 2}\sum s_j}
\Om_r(f;\pi/2^{s_1},...,\pi/2^{s_m})<\iy$$
and for some $p\in(0,2)$ and some $r$
$$\Sl_{s_1=-\iy}^\iy...\Sl_{s_m=-\iy}^\iy 2^{(1-{p\over 2})\sum s_j}
\Om_r^p(f;\pi/2^{s_1},...,\pi/2^{s_m})<\iy,$$
then $f=\hg,$ where $g\in L_1\cap L_p(\BR^m).$

b) If $f=\hg,$ where $g\in L\cap L_p(\BR^m)$ for some $p\in(0,2),$
and in addition $|u_j|\ge|v_j|,$ with $\sign u_j=\sign v_j$ and
$1\le j\le m,$ implies $|g(u)|\le|g(v)|,$ then the first series in a)
converges for each $r$ while the second one for $r>(1/p-1/2)m.$\ep

The result {\bf 1.7} immediately yields a generalization to the multiple
case of one theorem of A. Beurling,1949 (see the beginning of \S13).

Bernstein's precise result for a function to be of smoothness greater
than $m/2$ in $C$ is sufficient for belonging to $A(\BR^m)$ in a small
neighborhood of the endpoint, see [M22].
\flushpar
\proclaim{1.8} Let $p\in(0,1]$ and $q=[{m\over p}-{m+1\over 2}]$ (the integral
part). If a boundedly supported function $f\in C^q(\BR^m),$ and its
partial derivative $D_{q,j}(f)={\p^q f\over\p x_j^q}$ considered as
a function of $x_j$ has a finite number of points of inflection with
respect to the other variables (more precisely, points separating
inytervals of convexity) and $\om(D_{q,j};t)_\iy$ is the modulus of
continuity in the space $C(\BR^m),$ for $1\le j\le m,$ then the condition
$$\max\limits_j\Il_0^1 t^{qp+(m+1)(p/2-1)}\om^p(D_{q,j};t)\,dt<\iy$$
implies $f\in A_p(\BR^m)$ (or $\hf\in L_p(\BR^m)$).\ep

In the neighborhood of $\iy,$ and for even $f$ also near the origin,
these conditions can be weakened.
\flushpar
\proclaim{1.9} If for every $x\in\BR^m$ and for all $1\le j\le m$
$$f(x)=\Il_{|x_j|\le|u_j|<\iy}g(u)\,du\qquad\text{and}\qquad
\Il_{\BR^m}\esp_{|u_j|>|v_j|}|g(u)|\,dv<\iy,$$ then $f\in A(\BR^m).$\ep

Given conditions of convexity type, sharper results are obtained.
\flushpar
\proclaim{1.10} Let $f$ be a locally absolutely continuous function on
$[0,+\iy),$ and $\lim f(x)=0$ as $x\to+\iy$ and
$$||f||_{V^*}=\Il_0^\iy\esp_{u\ge x}|f'(u)|\,dx<\iy.$$

a)Let for every $x<0$ we have $f(x)=0.$ Then for every $y\in\BR\setminus
\{0\}$ $$\hf(y)=-{i\over y\sqrt{2\pi}}f({\pi\over 2|y|})+F(y),
\qquad\text{where}\qquad ||F||_1\le C||f||_{V^*}.$$

b) Let $f(x)=f_0(|x|),$ with $f_0\in C^{m_1}[0,\pi]$ for $m_1=[{m-1\over 2}]$
and the derivative $f_0^{(m_1)}$ is convex on $[0,\pi]$ and $|f_0^{(m_1+1)}(0)
<\iy.$ In order that $\hf\in L(\BR^m)$ (or $f\in A(\BR^m)$) it is necessary
and sufficient that $$\biggl|\Il_0^\pi u^{-{m+1\over 2}}f_0(\pi-u)\,du
\biggr|<\iy.$$\ep

Notice that it was G. E. Shilov,1942, who first investigated the
asymptotics of \Fc s of a convex function.

In the theory of partial differential equations, one has to study the
absolute convergence of the spectral resolution of an elliptic operator.
In the case of the Laplace operator on $\BR^m$ this means the absolute
convergence of the Fourier integral with respect to spheres. V. P. Maslov
(Dokl. AN SSSR, 1970) considered such a convergence for functions of two
variables with a weak discontinuity on a given curve and clarified a role
of the evolute of this curve for the divergence problem. The first result
on convergence (see the case $k=0$ in the next theorem) has been obtained
by E. M. Nikishin and G. I. Osmanov (Sib. Math. J.,1976).
\flushpar
\proclaim{1.11} Let $l$ be $k+2$ times continuously differentiable curve in 
$R^2.$ Let $E$ be the evolute of $l$ and $\k$ be the curvature of $l.$
Denote by $E_k$ the set of points of $E$ corresponding to points of $l$
so that $\k^{(p)}(s)=0$ for $p=1,...,k-1,$ and $\k^{(k)}(s)\ne 0,$
where $s$ is natural parametrization of $l.$

Let $f$ be boundedly supported continuous function
in $R^2,$ continuously differentiable outside $l.$ Let also for
$j=1,2$ and some $\varepsilon>0$ and $\sigma<{k+4\over 2(k+2)}$ the
following inequalities hold:
$$|\Delta_h {\partial f\over \partial x_j}(x)|\le{|h|^{\varepsilon}\over
\min(r^{\sigma}(x,l),r^{\sigma}(x+h,l))},$$
where $x,h\in R^2$ and $r(x,l)=\inf\limits_{y\in l}|x-y|.$ Then
$$J_0(x_0)=(2\pi)^{-1/2}\int\limits_0^\infty|\int\limits_{|t|=z}
\hat f(t)e^{itx_0}\,dt|\,dz<\infty$$ for every $x_0\in E^k.$\ep
\bigskip

For $R^m, m>2,$ this theorem is generalized
for functions weakly discontinuous on arbitrary compact $l$-dimensional
surface $L_l, 1\le l\le n-1.$ Sets $E_k$ are defined with respect to
different curves on $L_l$ (for $l>1$). The theorem itself can be formulated
in the same manner, interesting is the inequality for
the order of weak discontinuity of $[m/2]$-th derivative:
$$\sigma<{m-l\over m-1}-{k\over 2(m-1)(k+2)}.$$
One gets the two-dimensional case taking $m=2,$ $l=1.$

Let us consider also a question as to the structure of the set of real
zeros of the \Ft\ of the indicator function of a convex planar body $K.$
\flushpar
\proclaim{1.12} Let $h(\vp)$ be the support function of $K.$ Let $d(\vp)=
h(\vp)+h(\vp+\pi)$ be the width of $K$ in the given direction and $\dl(\vp)=
{1\over 2}(h(\vp+\pi)-h(\vp)).$ If $K$ is not a polygon, then the set
$N(K)$ of real roots of the equation $\Il_K e^{i(u,x)}\,du=0$ is
contained in $\cup_{p=1}^\iy M_p,$ where $M_p$ is a closed
continuous curve symmetric with respect to the origin and defined in the polar
coordinates by the equation $r=r_p(\vp).$ For each $p\in\BN$ we
have $2p\pi<d(\vp)r_p(\vp)<2(p+1)\pi.$ If $\dl\in C^k(\BR)$ (analytic),
then also $r_p\in C^k(\BR)$ (analytic) for every $p\in\BN.$ If the
set $N(K)\cap M_p$ is infinite, then $M_p\subset N(K)$ and $r_p$ is
an analytic function. If besides $K$ possesses a center of symmetry, then
$N(K)=\cup_{p=1}^\iy M_p.$ The curve $M_1$ may be not convex.\ep

In connection with one Beurling's theorem (see the beginning of \S13)
the following function space arised:
$$A_p^*(\BT)=\{f:f(x)=\Sl_{k\in\BZ}c_ke^{ikx},\quad
||f||_{A_p^*}^p=\Sl_{n=0}^\iy\sup\limits_{|k|\ge n}|c_k|^p<\iy\}.$$
The space $A_p^*(\BR)$ is defined similarly.

For $p\in(0,1],$ when passing from the inequality
$$||fg||_{A_p(\BT)}\le||f||_{A_p(\BT)}||g||_{A_p(\BT)}$$
to that for $A_p^*(\BT)$ a factor appears depending only on $p.$

Notice that an analog of the Wiener-Ditkin theorem holds for $A_p^*(\BT)$
as well as analogs of closed results, and the sharp sufficient
condition of belonging to $A_p^*(\BT)$ is the smoothness in $L$ (or $C$)
is greater than $1/p.$

The integral operator of order $\a$ $$f(x)\to\Il_0^1(1-t)^{\a-1}f(xt)\,dt$$
takes $A(\BT)$ continuously into $A(\BT)$ for all $\a>0$ and into $A^*(\BT)$
for all $\a\ge 1.$

The description of the dual space of $A^*(\BT)$ is contained in the
following proposition.
\flushpar
\proclaim{1.13} For each two sequences
$$\align \sup\limits_{\Sl_{n=0}^\iy\sup\limits_{k\ge n}|\a_k|\le 1}\biggl|
\Sl_{k=0}^\iy\a_k\be_k\biggr|&=\sup\limits_{n\ge 0}{1\over n+1}\Sl_{k=0}^n
|\be_k|,\\ \quad\\ \sup\limits_{{1\over n+1}\Sl_{k=0}^n |\be_k|\le 1}\biggl|
\Sl_{k=0}^\iy\a_k\be_k\biggr|&=\Sl_{n=0}^\iy\sup\limits_{k\ge
n}|\a_k|. \endalign$$\ep

Notice in conclusion that there are several definitions of $A^*$ in
the multidimensional case and some of them are applied to the study
of summability of the multiple \Fs\ at Lebesgue points (see 3.4 and
3.5 later on).
\bigskip
\flushpar
{\bf Bibliographical remarks.}\flushpar
For {\bf 1.1}, {\bf 1.7} -- {\bf 1.9}, see [T17]; for {\bf 1.2}
see [T8] (for a generalization to the multiple case, see [L1]); for {\bf
1.3} see [T25] (the case $r=0$ is contained, in essence, in [Be2]);
for {\bf 1.5} -- {\bf 1.6}, see [T13]); for {\bf 1.10}a), see [T10] and
[T24] (much more general result can be found in [L9]); for {\bf 1.10}b),
see [T13] (for a certain refinement, see [L2]); for {\bf 1.11}, see
[L3] in the two-dimensional case and [L5] for $\BR^m$ when $m\ge 2;$
for {\bf 1.12}, see [Z1-4]. Properties of $A^*(\BT)$ are given in [T9];
for a generalization to the multiple case, see [L10]); for {\bf 1.13},
see [BT1]. See also [T18].
\bigskip

\head 2. Multipliers and comparison of various methods of
summability of Fourier series as a whole. \endhead
Let $\{e_k\}$ be an orthogonal system of functions and $\{c_k\}=
\{c_k(f)\}$ be the Fourier coefficients of $f$ with respect to
this system.

A number sequence $\{\lm_k\}$ is called a multiplier in $L_p,$
written $\{\lm_k\}\in M_p,$ if the operator
$$\sum c_k e_k\to\sum\lm_k c_k e_k$$ takes $L_p$ continuously
into $L_p.$ We write $\{\lm_k\}\in M$ if the multiplier operator
takes $C$ into $C.$ Multipliers form a special class of integral
operators. Observe, that the norm of such integral operator can
be expressed explicitly via the kernel only in the cases $p=1$
and $p=\iy.$

Let us start with the trigonometric system $e_k=e^{i(k,x)},$ where
$k\in\BZ^m$ and $x\in\BT^m.$ In this case the multiplier
operator taking $L_p(\BT^m)$ into $L_p(\BT^m)$ commutes with shifts
and is expressed as the convolution (see e.g., [M22]). J. Marcinkiewicz
indicated convenient sufficient conditions for multipliers in $L_p(\BT^m)$
with $p\in(1,+\iy)$ and gave certain applications of multipliers
(see [M21]). For general properties of multipliers of trigonometric
Fourier series, see [M7], Ch.16.

In particular, in order a sequence $\{\lm_k\}_{k\in\BZ^m}$ to be
a multiplier in $C$ ($L_1$ or $L_\iy$), it is necessary and sufficient
that there exists on $\BT^m$ a finite complex valued Borel
measure $\mu$ such that for each $k\in\BZ^m$ we have $\lm_k=\Il_{\BT^m}
e_{-k}\,d\mu.$ In this case the operator is the convolution of a function
with the measure $\mu$ and
$$||\{\lm_k\}||_M=||\{\lm_k\}||_{M_1}=||\{\lm_k\}||_{M_\iy}=
\operatorname{var}\mu$$ (see also {\bf 2.4} below).

Observe that there exist plenty many of such multiplier operators.
\flushpar
\proclaim{2.1} For every linear bounded operator $A$ taking $L_p(\BT^m)$
into $L_p(\BT^m)$ for some $p\in[1,+\iy],$ or taking $C(\BT^m)$ into
$C(\BT^m),$ there exists a multiplier operator $A_0$ such that
$$||A_0||\le||A||,\qquad ||f-A_0f||_p\le\sup\limits_\theta
||f^\theta-A(f^\theta)||_p,$$ where $f^\theta(x)=f(x+\theta).$
And if also $A\ge0,$ then $A_0\ge0.$ \ep

As it is well known, every bounded operator is weakly compact
in the space $L_\iy(\BT^m)=(L_1(\BT^m))^*.$ The picture is different
for $C(\BT^m).$
\flushpar
\proclaim{2.2} If the multiplier $\{\lm_k\}_{k\in\BZ^m}$ is weakly
compact in $C(\BT^m),$ then it is expressed as a convolution with
some kernel from $L_1(\BT^m)$ and hence is compact. \ep

The set of those $k\in\BZ^m$ for which the Fourier coefficients
$c_k=c_k(f)$ do not vanish is called the spectrum of a function
$f\in L_1(\BT^m).$

Let $S$ be a non-empty set of $\BZ^m$ and $L_p(\BT^m,S),$ where
$p\in[1,+\iy],$ be a subspace of functions in $L_p(\BT^m)$ with
a spectrum in $S.$ Similarly $C(\BT^m,S)$ is defined. A number
sequence $\{\lm_k\}_{k\in S}$ is called a multiplier in $L_p(\BT^m,S),$
written $\{\lm_k\}_{k\in S}\in M_p(S),$ if for every $f\in L_p(\BT^m,S)$
the series $\sum\lm_k c_k(f)e_k$ is the Fourier series of some function
$\Lm f\in L_p(\BT^m,S)$ and $||\{\lm_k\}||_{M_p(S)}=||\Lm||_{L_p\to L_p}.$
For the space $C(\BT^m,S),$ we shall write $M(S)$ in contrast to
$M_\iy(S)$ when $p=\iy.$
\flushpar
\proclaim{2.3} For every $S\sbt\BZ^m$ and every $p\in(1,+\iy),$
and $q$ such that $1/p+1/q=1$ we have
$$||\{\lm_k\}||_{M_p(S)}=||\{\lm_k\}||_{M_q(S)}.$$\ep
\flushpar
\proclaim{2.4} Every multiplier in $C(\BT^m,S)$ can be extended
to a norm-preserving multiplier in $C(\BT^m)$ and
$$\align ||\{\lm_k\}||_{M(S)}&=\min\{\operatorname{var}\mu:\quad
\text{for $\mu$ such that for all $k\in S$ we have}\quad\Il_{\BT^m}
e_{-k}\,d\mu=\lm_k\}\\
&=\min\limits_{\{\lm_k\}_{k\not\in S}}\sup\limits_n(2\pi)^{-m}
||\s_n(\Lm)||_1=\min\limits_{\{\lm_k\}_{k\not\in S}}\lim\limits_{n\to\iy}
(2\pi)^{-m}||\s_n(\Lm)||_1,\endalign$$ where
$$\s_n(\Lm)=\Sl_{|k_j|\le n}\lm_k\prod\limits_{j=1}^m
\biggl(1-{|k_j|\over n+1}\biggr) e_k.$$ \ep

Notice two corollaries; the first one answers a question posed in
[M7], Vol.2.
\flushpar
\proclaim{2.5} {\rm a)} Given a sequence $\{\lm_k\}_{k\in S},$ in order that
for each function $f\in C(\BT^m)$ there exists at least one function
$\Lm f\in C(\BT^m)$ with the Fourier coefficients $c_k(\Lm f)=\lm_kc_k(f)$
for all $k\in S,$ it is necessary and sufficient that for some finite
Borel measure $\mu$ we have $\lm_k=\Il_{\BT^m}
e_{-k}\,d\mu$ for all $k\in S.$ 

{\rm b)} For every $S\sbt\BZ^m$ and every $p\in[1,+\iy)$ we have
$$||\{\lm_k\}||_{M_p(S)}\le||\{\lm_k\}||_{M(S)}=||\{\lm_k\}||_{M_\iy(S)}.$$
\ep

Will a multiplier be compact after extension from the spectrum
if it was such? Let us give an example in case $m=1.$
\flushpar
\proclaim{2.6} Let $S$ be a finite subset of $\BZ.$ In order that a
multiplier $\{\lm_k\}_{k\in S}$ can be extended to a compact
multiplier on $C(\BT)$ with the same norm, it is necessary and
sufficient that there exist $p\in S,$ $\a\in\BR,$ and a function
$f\in L_1(\BT)$ such that $e^{i\a}e^{-ipx}f(x)\ge0$ a.e. By this,
$||\{\lm_k\}||_{M(S)}=|\lm_p|.$ \ep
\flushpar
\proclaim{2.7} Let $p$ be even and $p\ne2.$

{\rm a)} There exist $S\sbt\BZ$ such that there exist an isometric multiplier
in $L_p(\BT,S)$ which is not extendable to a norm-preserving linear
operator in $L_p(\BT).$

{\rm b)} Let $\{\lm_k\}_{k\in\BZ}\in M_p$ and let exist $k_0\in\BZ$ such
that $|\lm_{k_0}|=|\lm_{k_0+1}|=||\{\lm_k\}||_{M_p}.$ Then the
multiplier $\{\lm_k\}$ is a shift, up to a factor. \ep

The question on extention of a multiplier in $C(\BT,S)$ with
preservation of norm and compactness is equivalent to the following one:
whether $L_1(\BT,\BZ\setminus S)$ is an existence space, written e. s.,
in $L_1(\BT)?$ For the definition of e.s., see e.g., [M13].
\flushpar
\proclaim{2.8} {\rm a)} If a set $S\sbt\BZ$ is such that every measure with
spectrum in $S$ is absolutely continuous with respect to Lebesgue
measure, or $S$ is the Sidon set (for definition, see {\rm [M10])}, then
$L(\BT,S)$ is e.s. in $L(\BT).$

{\rm b)} If $S=-S$ and $\BZ\setminus S$ is lacunary in the Hadamard sense,
then $L(\BT,S)$ is not e.s. in $L(\BT).$ \ep

Obviously, one can always assume $\lm_k=\vp(k)$ for $\vp$
sufficiently smooth. Consider more general case $\lm_k=\vp(\ve k)$
where $\ve$ is a positive parameter.
\flushpar
\proclaim{2.9} {\rm a)} If $\vp\in B(\BR^m),$ that is the \Ft\ of a finite
Borel measure on $\BR^m,$ then for each $\ve>0$ we have $||\{\vp(\ve k)
\}||_M\le||\vp||_B.$

{\rm b)} If $\vp$ is continuous a.e. on $\BR^m$ and for some sequence
$\ve_n\to 0$ $$H=\sup\limits_n||\{\vp(\ve_n k)\}||_M<\iy,$$
then the function $\vp$ can be corrected on the set of its points
of discontinuity so that it will be in $B(\BR^m)$ and $||\vp||_B\le H.$\ep

Hence we established connection between the algebras $M(\BR^m)$ and
$B(\BR^m):$ if $\vp\in C(\BR^m),$ then
$$||\vp||_M=\sup\limits_{\ve>0}||\{\vp(\ve k)\}||_M=||\vp||_B.$$
\flushpar
\proclaim{2.10} Let $\BZ_0=\BZ^m\setminus\{0\}.$

{\rm a)} If a function $\vp$ is continuous a.e. on $\BR^m\sm\{0\},$
and $$\varliminf\limits_{\ve\to+0}||\{\vp(\ve k)\}||_{M(\BZ_0)}=H<\iy,$$
the it can be corrected at the points of discontinuity and defined at
the origin by continuity so that $||\vp||_M=H.$

{\rm b)} If $\vp\in C(\BR^m),$ then
$$\Ll_{\ve\to+0}||\{\vp(\ve k)\}||_{M(\BZ_0)}=
\sup\limits_{\ve>0}||\{\vp(\ve k)\}||_{M(\BZ_0)}=||\vp||_M.$$ \ep

Let us give one sufficient condition of boundedness of the norms of
of a sequence of multipliers.
\flushpar
\proclaim{2.11} Let $|\lm_0|+\ln(n+1)(|\lm_n|+|\lm_{-n}|)\le1.$
If also for some $\dl>0$ one of the following two conditions,
either {\rm a)} or {\rm b)}, hold

{\rm a)} $|\lm_k-\lm_{k+s}|\le\sqrt{s/n}|\ln(s/3n)|^{-1-\dl},\qquad
1\le s\le n,\quad -n\le k\le n-s,$

{\rm b)} $\Sl_{k=-n}^n|\lm_k-\lm_{k+1}|\le 1,\quad|\lm_k-\lm_{k+s}|
\le|\ln(s/3n)|^{-2-\dl},\qquad 1\le s\le n,\quad -n\le k\le n-s,$
\flushpar
then $||\{\lm_k\}_{-n}^n||_M$ are bounded by a constant depending
only on $\dl.$ For $\dl=0$ no one of these statements valids. \ep

Let us go on to comparison of the two multiplier operators
$$\Lm f\sim \sum\lm_k c_k e_k\qquad\text{and}\qquad
\tilde\Lm f\sim \sum\tilde\lm_k c_k e_k$$ and formulate a comparison
principle.
\flushpar
\proclaim{2.12} Let $S_0=\{k\in\BZ^m:\quad\lm_k=0\}.$

{\rm a)} If $\tilde\lm_k=0$ for every $k\in S_0,$ which is also
necessary, and $K=\inf\limits_{0/0}||\{\tilde\lm_k/\lm_k\}||_M<\iy,$ where
the least lower bound is taken with respect to choice of values of
fractions of type $0/0,$ then for each function $f$ such that $\Lm f\in
C(\BT^m)$ we have $||\tilde\Lm f||_\iy\le K||\Lm f||_\iy.$ And vice
versa, if this inequality holds for each function $f$ such that
$\Lm f\in C(\BT^m)$ and $\{1/\lm_k\}\in M(\BZ^m\sm S_0),$ then
$$\inf\limits_{0/0}||\{\tilde\lm_k/\lm_k\}||_M=
\min\limits_{0/0}||\{\tilde\lm_k/\lm_k\}||_M\le K.$$

{\rm b)} If $\tilde\lm_k=0$ for every $k\in S_0$ and $\{1/\lm_k\}\in
M(\BZ^m\sm S_0),$ while $\{\tilde\lm_k/\lm_k\}\in M(\BZ^m\sm S_0)$ and
compact, then $$\sup\limits_{f\in C:\,||\Lm f||_C\le 1}||\tilde\Lm f||_C=
\sup\limits_{f\in L_\iy:\,||\Lm f||_{\iy}\le 1}||\tilde\Lm f||_C=
\min\limits_{0/0}||\{\tilde\lm_k/\lm_k\}||_M.$$

{\rm c)} If $\Lm$ is a compact operator in $C(\BT^m)$ and the equality
$\lm_k=1$ implies $\tilde\lm_k=1,$ then we have $||f-\tilde\Lm f||_C\le
K||f-\Lm f||_C$ for each $f\in C(\BT^m)$ if and only if
$$\inf\limits_{0/0}||\{(1-\tilde\lm_k)/(1-\lm_k)\}||_M\le K.$$ \ep
\flushpar
\proclaim{Example} Continuity of $$D_{2r}(f)=\Sl_{j=1}^m{\p^{2r}f\over
\p x_j^{2r}}$$ for $r\ge 2$ and $m\ge 2$ does not imply continuity of
$\Dl^r f,$ where $\Dl$ is the Laplace operator. \ep
To prove this, we consider the functions to be periodic and compare
the \Fs\ of $D_{2r}$ and $\Dl^r.$ Then {\bf 2.12}a), {\bf 2.10}a) and
{\bf 2.9}b) yield that the function $\vp(x)=(\Sl_j x_j^{2r})
(\Sl_j x_j^2)^{-r}$ must have a limit as $x\to 0.$ But this function
is homogenious of zero order and hence constant.

Let us generalize now the comparison principle.

Let $E$ be a complex Banach space and $L(E)$ the Banach algebra of
linear continuous operators acting in $E.$ Let further $\{\Cp_k\}_0^\iy$
be a complete sequence of mutually orthogonal projectors in $E.$ This
means that

$\a$) $\Cp_k\in L(E)$ for every $k\in\BZ_+;$

$\be$) $\Cp_k f=0$ for each $k\in\BZ_+$ implies $f=0;$

$\g$) $\Cp_k\Cp_s=\dl_{ks}\Cp_s$ for every $k,s\in\BZ_+,$ where $\dl_{ks}$
is the Cronecker delta;

$\dl$) $||\Cp_k||\ne0$ for every $k\in\BZ_+.$
\flushpar
Let us associate with $f\in E$ the series $\Sl_{k=0}^\iy\Cp_kf.$ 
By $\be$) the correspondence between the elements and the series
indicated is one-to-one correspondence.
Such series may be called "the \Fs".
A number sequence $\{\lm_k\}$ is called the multiplier in $E,$ with
respect to a given sequence of orthoprojectors, if for every $f\in E$
there exists $g\in E$ satisfying the condition $\Cp_kg=\lm_k\Cp_kf$ 
for every $k\in\BZ_+,$ that is $g=\Lm f\sim\sum\lm_k\Cp_kf$ and 
$\Lm\in L(E),$ and $||\{\lm_k\}||_M=||\Lm||.$
\flushpar
\proclaim{2.13} {\rm a)} If $\{\lm_k\}_0^\iy$ is a multiplier and
$\lm_k=1$ implies $\tilde\lm_k=1,$ the latter is also necessary, and
$$K=\inf\limits_{0/0}||\{(1-\tilde\lm_k)/(1-\lm_k)\}||_M<\iy,$$
then for every $f\in E$ we have $||f-\tilde\Lm f||\le K||f-\Lm f||.$

{\rm b)} If the inequality indicated holds for each $f\in E$ and
$\Lm$ is a compact operator, then
$$\inf\limits_{0/0}||\{(1-\tilde\lm_k)/(1-\lm_k)\}||_M\le K+
\Sl_{\nu:\,\lm_\nu=1}||\Cp_\nu||,$$ where in case $\lm_k\ne1$ for each
$k\in\BZ_+$ the sum is absent on the right-hand side as well as the
least lower bound on the left-hand side. \ep

With the availability of multiplier theorems it will be possible to
apply the comparison principle to basis expansions, to orthogonal
and biorthogonal expansions, etc.

Let us compare now classical summability methods as regards to the
rate of convergence to elements of Banach space.
\flushpar
\proclaim{2.14} If the $(C,1)$ method is regular, that is
$$\s_n(f)=\Sl_{k=0}^n(1-{k\over n+1})\Cp_k f\to f\qquad\text{as}\quad
n\to\iy$$ in the norm of a given space $E$ for every $f\in E,$ then
it is equivalent to the Abel-Poisson method, that is for every $r\in(0,1)$
and for each $f\in E$ we have
$$||f-\s_n(f)||\asymp||f-\Sl_{k=0}^\iy r^k\Cp_k f||,$$ where $n=[1/(1-r)]$
is the integral part of $1/(1-r)$ and in the two-sided inequality replaced
by $\asymp$ positive constants do not depend on $f$ and $r.$ \ep

To replace $\s_n$ by $\s_n^\a,$ that is $(C,\a)$ method, for small $\a$
more restrictive assumptions on $\sum\Cp_kf$ are needed.
\flushpar
\proclaim{2.15} {\rm a)} Let for some $\a>0$ the following condition
$(S,\a)$ hold: if for $\lm_{n+1}=0$
$$\Sl_{k=0}^n|\lm_k-\lm_{k+1}|\le1\quad\text{and}\quad |\lm_k-\lm_{k+s}|
\le(s/n)^\a$$ for $1\le s\le n$ and $0\le k\le n+1-s,$ then
$$\sup\limits_n\sup\limits_{||f||\le1}||\Sl_{k=0}^n\lm_k\Cp_kf||<\iy.$$
Then $$||f-\s_n^\a(f)||\asymp||f-\s_n(f)||.$$

{\rm b)} Let the condition $(S,\a)$ hold and let $\vp\in\Lip\a$ and
be a function of bounded variation on $[0,1].$ Suppose further that
$0\le\vp(1)\le\vp(x)<\vp(0)$ on $(0,1).$ Then either for each $\dl>1,$ 
or for each $\dl>0$ provided $\vp$ is piecewise monotone on $[0,1],$
for which the condition $(S,\dl\a)$ holds, we have $$||f-\Sl_{k=0}^n
\vp^\dl(k/n)\Cp_kf||\asymp||f-\Sl_{k=0}^n\vp(k/n)\Cp_kf||.$$ \ep
For instance, the condition $(S,\a),$ for each $\a>0,$ is satisfied by
the trigonometric system (see {\bf 2.11}b)) and by the Walsh system
(see {\bf 8.5}).

For $\vp(x)=(1-x^\a)_+^\dl$ we get the Riesz means. Hardy established that
the Riesz methods, as applied to number series, become stronger for
larger $\dl$ and fixed $\a$ and are equivalent one to another, that is
the summability by one of them yields the summability by another, for
different $\dl$ and fixed $\a.$ The comparison of summability methods
defined as above is of other type. For the same Riesz methods as applied
to the \Fs, another picture occurred to be true: the growth of $\dl$
does not change anything while the growth of $\a$ improves the convergence
(more precisely, a speed of convergence). See also \S5.

To obtain more general sufficient conditions than $(S,\a),$ the strong
summability can be applied. This notion was introduced by Hardy and
Littlewood for trigonometric series, the so-called condition $(H,p), p>0:$
$$\sup\limits_n\sup\limits_{||f||\le1}||{1\over n+1}\Sl_{s=0}^n
|\Sl_{k=0}^s\Cp_kf|^p||<\iy.$$
The case $E=C$ may be understood, for example. The case $p=1$ can be
found in essence in {\bf 1.13}.
\flushpar
\proclaim{2.16} Let $x=\{x_k\}_1^\iy$ and $y=\{y_k\}_1^\iy.$ Define
$$||x||_{h_p}=\sup\limits_n({1\over n}\Sl_{k=1}^n|x_k|^p)^{1/p}$$ and
$$||y||_{b_p}=\Sl_{n=1}^\iy({1\over n}\Sl_{k=n}^\iy|y_k|^p)^{1/p}.$$
For $p\in(1,+\iy)$ and $1/p+1/q=1$ the following inequalities hold.
$$\align |\Sl_{k=1}^\iy x_k y_k|&\le\g_1(p)||x||_{b_p}||y||_{h_q},\\
\sup\limits_{||x||_{b_p}\le1}|\Sl_{k=1}^\iy x_k y_k|&\ge\g_2(p)||y||_{h_q},\\
\sup\limits_{||y||_{h_q}\le1}|\Sl_{k=1}^\iy x_k y_k|&\ge\g_3(p)||x||_{b_p}.
\endalign$$   \ep

Let us go on to orthogonal series in the space $L_{2,h}(a,b),$ where
$(a,b)$ is the finite or infinite interval in $\BR.$ Let $\{\vp_k\}_0^
\iy$ be orthonormal system of complex valued functions on $(a,b)$ with
weight $h.$ This means that
$$\Il_a^b\vp_k(t)\bar\vp_s(t)h(t)\,dt=\dl_{ks}.$$
Let us assume that $\vp_k\in L_\iy(a,b)$ for all $k\in\BZ_+$ while
$h\in L(a,b)$ and $h(t)>0$ a.e. on $(a,b).$

To each measurable function $f$ satisfying $\Il_a^b|f(t)|h(t)\,dt<\iy$
its \Fs\ with respect to the system $\{\vp_k\}_0^\iy$ corresponds:
$$f\sim\Sl_{k=0}^\iy c_k\vp_k,\qquad c_k=c_k(f)=\Il_a^b f(t)\bar
\vp_k(t)h(t)\,dt.$$

In the following statement, where also a question on multipliers
with regard to the position of a point $x\in(a,b)$ is considered,
the $(C,1)$ method is used which can be replaced by any Toeplitz
regular method.
\flushpar
\proclaim{2.17} {\rm a)} If
$$\sup\limits_n\Il_a^b|\Sl_{k=0}^n\lm_k(1-{k\over n+1})\bar\vp_k(t)
\vp_k(x)|\,h(t)\,dt=K(x)$$ is in $L_\iy(a,b),$ then for every
$f\in L_\iy(a,b)$ the series $\sum\lm_k c_k(f)\vp_k$ is the \Fs\ of
some function $\Lm f$ to which it converges in average, that is in
the space $L_2$ with the weight $h,$ and a.e. on $(a,b)$
$$|(\Lm f)(x)|\le K(x)||f||_\iy.$$

{\rm b)} If in addition to that $\vp_0\equiv C$ and for all $n\in\BZ_+,$
for all $x\in(a,b),$ and for almost all $t\in(a,b)$ we have
$$\Sl_{k=0}^n(1-{k\over n+1})\vp_k(t)\bar\vp_k(x)\ge 0,$$  then
$$||\{\lm_k\}||_{M_\iy}=\esp_{x\in(a,b)}\sup\limits_n\Il_a^b
|\Sl_{k=0}^n\lm_k(1-{k\over n+1})\bar\vp_k(t)\vp_k(x)|\,h(t)\,dt.$$

{\rm c)} If in addition $\vp_k\in C(a,b)$ for all $k\in\BZ_+$ and
for every $f\in C(a,b)$
$$\Ll_{n\to\iy}||f-\Sl_{k=0}^n(1-{k\over n+1})c_k(f)\vp_k||_C=0,$$
then $\{\lm_k\}_0^\iy$ is a multiplier taking $C$ into $C$ and
$$||\{\lm_k\}||_M=||\{\lm_k\}||_{M_\iy}.$$ \ep

It is known that when $[a,b]=[-1,1]$ for the system of Jacobi polynomials,
that is when $h(t)=(1+t)^\a(1-t)^\be$ (see [M20]), nonnegative are
$(C,\a+\be+2)$ means (see G. Gasper, SIAM J.Math.An.,1977). For the
system of Hermite polynomials, that is for $(a,b)=(-\iy,+\iy)$ and
$h(t)=e^{-t^2},$ nonnegative is the Poisson kernel.

Let us return to the trigonometric system for which a greater
progress can be achieved.
\flushpar
\proclaim{2.18} {\rm a)} If $\{1/\lm_k\}_{k\in\BZ^m\sm S_0}\in
M(\BZ^m\sm S_0)$ (see {\bf 2.12}), then the inequality $||\tilde\Lm f||_\iy
\le K||\Lm f||_\iy$ being true with the finite right-hand side
for each $f\in C(\BT^m)$ implies that for every $p\in[1,+\iy]$ the
inequality $||\tilde\Lm f||_p \le K||\Lm f||_p$ holds with the finite
right-hand side for each $f\in L_p(\BT^m).$

{\rm b)} If $\vp\in B(\BR^m),$ $\psi\in A(\BR^m),$ and $\psi(x)\ne1$ for
all $x\ne0,$ then the least constant in the inequality
$$||f-\Sl_k\vp(k/n)c_k e_k||_\iy\le K||f-\Sl_k\psi(k/n)c_k e_k||_\iy$$
being true for all $f\in C(\BT^m)$ for all $n$ is
$$K=||(1-\vp)/(1-\psi)||_B.$$

{\rm c)} Let the inequality from {\rm b)} hold with some constant
independent of $f$ and $n.$ If in addition $\vp\in C(\BR^m),$ $\psi\in
B(\BR^m),$ and $\psi(x)\ne1$ for all $x\ne0$ and there exist limits
of $\psi$ and $\vp$ as $|x|\to\iy$ and the first of them equals to zero,
and besides, outside of some neighborhood of the origin $\vp$ and $\psi$
are of bounded variation in the sense of Vitali and $|1-\psi(x)|\le
C|1-\vp(x)|$ for all $x\in\BR^m,$ then the inverse inequality to that
indicated in {\rm b)} holds with some other positive constant independent
of $f$ and $n.$  \ep

Let us give an additional result in case of $L$ metric and $m=1.$
\flushpar
\proclaim{2.19}  {\rm a)} If $\lim\lm_k=0$ as $|k|\to\iy$ and
$\Sl_{-\iy}^\iy|\lm_k-\lm_{k+1}|<\iy,$ then
$$||\{\lm_k\}||_{M_1(\BZ_0)}=\inf\limits_{\lm_0}{1\over 2\pi}\Il_{-\pi}^\pi
|\Sl_{-\iy}^\iy\lm_ke^{ikt}|\,dt-\theta\Sl_k|\lm_k-\lm_{k+1}|,$$
where $\theta\in[0,C].$

{\rm b)} Let for $f\in L[a,b]$ we have $||f-\lm_0||_1=\min\limits_{\lm\in\BC}
\Il_a^b|f(x)-\lm|\,dx$ and there exist a subset $e\sbt[a,b]$ on which
$f$ is bounded such that $\g\mes e=b-a$ for some $\g\in[1,2).$ Then
$$\sup\limits_{x\in e}|f(x)-\lm_0|\le{3\over 2\sqrt{2-\g}}\om^e(f)$$
where $\om^e(f)$ is the oscillation of $f$ on $e.$  \ep
\flushpar
{\bf Bibliographic remarks.}
\flushpar
The comparison principle for various trigonometric \Fs\ and first results
like {\bf 2.14} and {\bf 2.15} appeared in [T3]. The comparison principle
in other form was published also by H. Shapiro (Bull.AMS, 1968). It turned
out that the question on comparison of series with respect to their
approximation properties was posed by J. Favard as early as in 1963.
A generalization to series in Banach spaces was formulated by P. Butzer,
R. Nessel, and W. Trebels (see [M29]). Statement 2.13 from [T12] is a
refinement of this result.
\flushpar
For {\bf 2.2--2.6, 2.10, 2.12,} see [T29]; for {\bf 2.7,} see [Le2];
for {\bf 2.8,} see [T22]; for {\bf 2.9, 2.18,} see [T17];
for {\bf 2.11,} see [T29]; for {\bf 2.14, 2.15, 2.17,} see [T18,27];
{\bf 2.16} is due to E. Belinskii (see [BT1]).
\flushpar
It is known that in order that a sequence $\{\lm\}$ is a multiplier
from $L_\iy$ into $C$ with respect to the trigonometric system, it is
necessary and sufficient that the series $\sum\lm_k e_k$ is the \Fs\
(see, e.g., [M7]). Thus various sufficient conditions for multipliers
from $L_\iy$ into $C$ may be found in papers of S. A. Telyakovskii,
G. A. Fomin, Ya. S. Bugrov, Yu. L. Nosenko ([N7,8]), O. I. Kuznetsova
([Ku4,6]), E. R. Liflyand ([L9]), P. V. Zaderey. See also [S14,4].

\head 3. Regularity of summability methods of Fourier series.\endhead
Let $f\in L(\BT^m)$ and $\sum c_ke_k$ be its trigonometric \Fs.
Let us consider for $n\in\BR^m_+$ a matrix $\{\lm_{n,k}\}_{k\in\BZ^m}$ 
and introduce the following polynomial means
$$\Lm_n^0(f)=\Sl_{-n\le k\le n}\lm_{n,k}c_ke_k.$$
For $\lm_{n,k}=1$ or $0$ we obtain various partial sums (rectangular,
circular, etc. for $m=2$).

Let in addition $E\sbt L(\BT^m).$ Summability method is called {\it
regular} (in one sense or another) if for every $f\in E$ we have
$\Lm_n^0(f)\to f$ as $n\to\iy$ at all the points of a subset of $\BT^m$
or in norm if one is defined in $E.$ Such a regularity is not
connected with the Toeplitz regularity and thus is called $F$-regularity.

An investigation of summability methods of simple series, $m=1,$ has a long
history (see, e.g., [S14]).

In 1936, S. Bochner introduced the following spherical means of Riesz type
$$S_n^\dl(f)=\Sl_{|k|\le n}(1-{|k|^2\over n^2})^\dl c_ke_k,$$
with $n\in\BN,$ and clarified a role of the critical order $\dl={m-1\over2}$
in the multiple case (see, e.g., [M22]). In 1972, Ch. Fefferman proved
that $S_n^0,$ that is spherical partial sums, converge on $L_p(\BT^m)$
for $m\ge 2$ when $p=2$ only (see, e.g., [S6]). Up to the present
the regularity in $L-p(\BT^m)$ of the Bochner-Riesz means $S_n^\dl$ is
not investigated in full yet for all $\dl>0$ and $m\ge3$ (see, e.g., [M5]).

If $E$ is Banach space, the question is reduced to the boundedness of
the sequence of operator norms $||\Lm_n^0||,$ the so-called Lebesgue
constants. If the trigonometric system is closed in $E,$ then the
condition $\lm_{n,k}\to1,$ as $n\to\iy,$ for all $k\in\BZ^m$ 
is the one that should be added
to the boundedness of $||\Lm_n^0||.$ The \Lc\ $||\Lm_n^0||$ are the
norms of multipliers $\{\lm_{n,k}\}$ which in the case of spaces $C$ or
$L$ and $\lm_{n,k}=\vp(k/n)$ are closely connected with belonging of
$\vp$ to $B$ or $A$ (see {\bf 2.9} and {\bf 1.2}).
\proclaim{3.1} {\rm a)} If $\vp\in C(\BT)$ and $\vp(\pi)=\vp(0)=0,$ then
$$||\vp||_M=\sup\limits_{\ve>0}||\vp(\ve k)||_M\le C\Sl_{k=-\iy}^\iy
|c_k(\vp)|\ln(|k|+1).$$ And if in addition $\vp$ is continuous real even
functions with alternating, with respect to sign, \Fc s in cosines
starting with the first one, then the opposite inequality valids too.

{\rm b)} If $\vp\in C(\BT)$ and $\vp(\pi)=\vp(-\pi)=\vp(0)=0,$ then
denoting $\vp_0(x)=\vp(x)\sign x$ we have
$$||\vp||_M\le C\Sl_{k=-\iy}^\iy |c_k(\vp_0)|\ln(|k|+1),$$ under
assumption that the series on the right-hand side converges.

{\rm c)} If for $m_1=[{m-1\over 2}]$ we have $\vp_0\in C^{m_1}[0,+\iy)$
and $\vp_0$ being supported in $[0,\pi]$ is such  that for every
$t\in[0,\pi]$ $$\vp_0(t)=\Sl_{\nu=0}^\iy\a_\nu\cos\nu t\quad\text{and}
\quad\Sl_{\nu=0}^\iy\nu^{m-1\over 2}|\a_\nu|\ln(\nu+1)<\iy,$$
then $\vp(x)=\vp_0(|x|)\in M(\BR^m).$ \ep

Several authors, R. P. Boas, J.-P. Kahane, S. I. Izumi and T. Tsuchikura,
I. Wik, independently and almost
simultaneously considered the following question: known is the
behavior of absolute values of the \Fc s of $\vp$, how is it possible
to derive that $\vp_0(x)=\vp(x)\sign x\in A(\BT).$ The answer obtained (see 
[M10, Ch.6, \S2, Cor.4),5)]) immediately follows from the results given above.
\proclaim{3.2} Let $\vp$ be continuous a.e. on $\BR^m$ and $\{\vp(\ve k)
\}\in M_1$ for every $\ve>0.$ In order that $\sum\vp(\ve k)c_k(f)e_k$
converges to $f,$ as $\ve\to+0,$ on $L(\BT^m)$ or $C(\BT^m)$ it is
necessary and sufficient that after correction at the points of
discontinuity $\vp\in B(\BR^m)$ and $\vp(0)=1.$ \ep
\proclaim{Example} $\vp(x)=(1-\Sl_{j=1}^m|x_j|^\a)_+^\dl$ for $\a=1$
and $\dl>0$ or for positive $\a\ne 1$ and $\dl>{m-1\over 2}.$ \ep

Let us go on to summability a.e. Sufficient conditions can be derived
from Marcinkiewicz's theorem on the strong summability a.e. and
{\bf 2.16} (see [BT1]). Let us give a necessary condition.
\proclaim{3.3} Let $p\in[1,2)$ and $\vp$ be continuous a.e. on $\BR^m$
and boundedly supported. If on some subset of $\BT^m$ of positive
measure $\sum\vp(\ve k)c_k(f)e_k$ converges to $f,$ as $\ve\to+0,$
for all $f\in L_p(\BT^m),$ then for each $q>p$
$$\Il_{\BR^m}|\hat\vp(u)|^q\,du<\iy.$$ \ep

This yields immediately K. I. Babenko's result, 1971: for every
$p\in[1,{2m\over m+1})$ for $m\ge 2$ there exists a function
$f\in L_p(\BT^m)$ such that its spherical Bochner-Riesz means
$S_n^\dl(f)$ do not converge to $f$ on a set of positive measure
if $\dl\in[0,{m\over p}-{m+1\over 2}).$

As for the summability at Lebesgue points, criteria do exist.
Let us start with the compact case.
\proclaim{3.4} In order that $\sum\vp(k/n)c_k(f)e^{i(k,x)}$
converges to $f$ as $n\to\iy$ for each $f\in L(\BT^m)$ at all
its Lebesgue points, that is those for which
$$\Ll_{r\to+0}r^{-m}\Il_{|x-u|\le r}|f(x)-f(u)|\,du=0,$$
it is necessary and sufficient that $\vp\in A^*(\BR^m)$ (see
after {\bf 1.12}) and $\vp(0)=1.$ \ep

There exists a similar criterion for "strong" Lebesgue points,
that is those for which the ball in the definition is replaced by
arbitrary parallelepiped with faces parallel to coordinate planes
(see [Be1]), and {\bf 3.4} in terms of belonging to $A^*(\BT^m)$
(see [T9] and [L1]).

The general case is investigated in the following theorem.
\proclaim{3.5} Let $\vp\in B(\BR^m).$ In order that for each
$f\in L(\BT^m)$ its \Fs\ is summable at all Lebesgue points
by a method generated by the function $\vp$ it is necessary and
sufficient that $\vp\in A^*(\BR^m)$ and $\vp(0)=1.$ \ep

Observe that if $\vp$ satisfies convexity type conditions, then
$\vp\in A(\BT)$ and $\vp\in A^*(\BT)$ simultaneously. The point is
that in the asymptotics {\bf 1.10}a) for this case, but not for the
general one, we have in addition that also $\sup\limits_{|y|\ge x}
|F(y)|\in L[0,+\iy).$

Let $\{\nu_k\}_0^\iy$ be strictly increasing integer-valued sequence.
When does the following statement
$$\Ll_{n\to\iy}{1\over n+1}\Sl_{k=0}^n|f(x)-S_{\nu_k}(f;x)|=0$$
valids? R. Salem, 1955, proved that in order that this equality holds
for every $f\in C(\BT)$ it suffices that the sequence $\{\nu_k\}$
is of power growth.
\proclaim{3.6} {\rm a)} If the sequence $\{\nu_k\}$ is convex, that
is $\nu_{k+2}-2\nu_{k+1}+\nu_k\ge0$ for all $k\ge0,$ then in order
that the afore-mentioned limiting equality valids for every $f\in
C(\BT)$ everywhere or uniformly on $\BT,$ it is necessary and
sufficient that $\ln\nu_n=O(\sqrt{n}).$

{\rm b)} For the multiple case, when cubic partial sums of the \Fs\
of $f$ are taken as $S_{\nu_k},$ the answer is the following under
the same conditions: \quad $\ln\nu_n=O(n^{1/2m}).$

{\rm c)} For $m=1$ and $\nu_k=[2^{k^\a}],$ it is impossible to
take $\a>1/2$ even in the case when the sign of absolute value
is removed outside the sign of sum, that is approximation by
arithmetic means with gaps. \ep

Let us investigate as an example the question of regularity of the following
Bernstein-Rogosinski type means in rather general form on $\BT^2:$
$$R_n(f;x)=R_n(f;W;\g,\mu,x)=\Il_{\BR^2}S_n(f;W;x-\g u/n)\,d\mu(u),$$
where for $n\in\BN$ $$S_n(f;W;x)=\Sl_{k\in nW}c_k(f)e^{i(k,x)}$$
is the partial sum of the \Fs\ of a function $f\in L(\BT^2)$ generated by
a bounded set $W\sbt\BR^2,$ $\g>0,$ and $\mu$ is a finite complexvalued
Borel measure on $\BR^2$ and $\Il_{\BR^2}d\mu=1.$

Let us start with the case of discrete measure, that is a measure
concentrated at finite number of points.
\proclaim{3.7} {\rm a)} Let $W$ be a parallelogram symmetric with
respect to the origin. If a measure $\mu$ is concentrated at four
points, which is the least number, then the regularity in $C$ valids
if and only if these points are of equal measure and are the vertices
of a parallelogram $W_1$ with the sides perpendicular to those of $W$
and ratio of the product of lengths of the two perpendicular sides
of $W$ and $W_1$ equals to the ratio of odd numbers.

{\rm b)} Let $W$ be a polygon with the origin inside such that no one of
the sides does not lie on the line passing through the origin. If $m$
is a maximal number of its sides no one of which are not parallel
one to another, then it is possible to choose $2^m$ points with the
samr measure so that the means $R_n$ are regular in $C(\BT^2).$

{\rm c)} If $W$ is a circle, then for every discrete measure and for
every $\g$ the regularity in $C(\BT^2)$ is impossible. \ep
\proclaim{3.8} {\rm a)} Let measure $\mu$ be with compact support and $W$
be a bounded connected set containing a neighborhood of the origin
and satisfying another two conditions: the boundary $\p W$ is of zero
two-dimensional Lebesgue measure and in every neighborhood of each
boundary point there are both inner points from $W$ and $\BR^2\sm W.$
Then for regularity of $R_n$ in $C(\BT^2)$ it is necessary and
sufficient that $$\Il_{\BR^2}e^{-i\g(x,u)}\,d\mu(u)=0$$ for all $x\in\p W.$

{\rm b)} If $K$ is a subgraph of a nonnegative function continuous
on an interval and non-constant, and measure $\mu$ is uniformly
distributed on the area of $K,$ the the method $R_n$ is not regular
in $C(\BT^2)$ for any choice of $\g>0$ and $W.$

{\rm c)} Let $K$ be strictly convex bounded body in $\BR^2$ symmetric
with respect to the origin and the boundary $\p K\in C^\iy.$ Let
measure $\mu$ be uniformly distributed on its area and $\g\p W=Mp$
where $M_p$ is the curve from {\bf 1.12}. If in addition $p$ is big
enough, then for every $f\in L(\BT^2)$ at each of its Lebesgue points
$\Ll_{n\to\iy}R_n(f;W;x)=f(x).$ \ep

Note that in order to get a regular method in $C(\BT^m)$ for $m>2$
measure $\mu$ averaging of $S_n$ should be done several times.
\flushpar
{\bf Bibliographical remarks.}
\flushpar
For results on multiple \Fs\ obtained during the last 20-25 years,
see [S6, M22, M5].
\flushpar
For {\bf 3.1}a) and b); see [T3], for {\bf 3.1}c), see [T18];
for {\bf 3.2}, see [T17]; for {\bf 3.3}, see [Be6];
for {\bf 3.4}, see [Be1,4]; for {\bf 3.5}, see [BT2];
for {\bf 3.6}a), see [ZgT], and independently L. Carleson (conference
in Budapest, 1979); for {\bf 3.6}b), see [Ku3], and also [Ku6,8,10]
where some other problems of strong summability of multiple \Fs\ are
investigated; for {\bf 3.6}c), see [Be9].
\flushpar
The means $R_n$ in the partial case of $W$ being a circle and measure
$\mu$ is uniformly distributed on $\p W$ were investigated by S.
Minakshisundaram and K. Chandrasekharan in 1947. For the general form of
$R_n,$ examples and necessary condition from {\bf 3.8}a), see [T17];
for {\bf 3.7}a),b), see [N4] (see also [N6,9]); for {\bf 3.7}c) and
for {\bf 3.8}a)-c), see [Za1,2,4].

\head 4. Lebesgue constants and approximation of function classes.\endhead
A. Lebesgue found the asymptotics of operator norms of partial sums $S_n$ as
follows: $$\sup\limits_{||f||_C\le1}||S_n(f)||_C=(2\pi)^{-1}
\Il_{\BT}|\Sl_{k=-n}^n e^{ikt}|\,dt=4\pi^{-2}\ln n+O(1).$$
It is easy now to write any number of members of the series
$4\pi^{-2}\ln n+C_0+C_1n^{-1}+....$

E. Landau found the asymptotics
$||S_n||$ for power series, or for trigonometric series with spectrum
in $\BZ_+:$ the factor $4\pi^{-2}$ is replaced by $\pi^{-1}.$

If the \Lc\ $||\Lm_n||$ increase infinitely as $n$ increases, then there is
no regularity of the means $\Lm_n(f).$ Then for the convergence of $\Lm_n(f)$
to $f,$ some smoothness should be added to $f$ so that it is connected
with the rate of growth of $\Lm_n.$

A. N. Kolmogorov, 1935, has found the asymptotics of approximation of the
class $W^r,$ with $r\in\BN,$ by partial sums $S_n.$ Here
$$W^r=\{f:\quad f^{(r-1)}\quad\text{is absolutely continuous and}\quad
||f^{(r)}||_\iy\le1\}$$ and
$$\sup\limits_{f\in W^r(\BT)}||f-S_n(f)||_C=4\pi^{-2}n^{-r}\ln n+O(n^{-r}).$$
This Kolmogorov's problem has a long history (for refinement of this
formula, generalizations, etc., see, e.g., [S15]).
\proclaim{4.1} For every $r>0,$ it is not obligatory for $r$ to be integer,
four numbers $A,$ $B,$ $C,$ and $D$ can be pointed out so that
$$\align\sup\limits_{f\in W^r(\BT)}||f-S_n(f)||_C&= 4\pi^{-2}n^{-r}
\ln n+An^{-r}+Bn^{-r-1}\ln n\\ &+Cn^{-r-1}
+Dn^{-r-2}\ln n+O(n^{-r-2}).\endalign$$  \ep

Let us go on to the case of the space $C(\BT^m),$ or $L(\BT^m),$ for
$m\ge 2.$ The growth of the \Lc\ of the spherical Bochner-Riesz means
(see the beginning of \S3) was obtained by V. A. Ilyin, 1968, for $\dl=0,$
and by K. I. Babenko, 1971, for $\dl\in[0,{m-1\over 2}),$ as follows:
$$\sup\limits_{||f||_C\le1}||S_n^\dl(f)||_C\asymp n^{{m-1\over 2}-\dl}$$
with $\dl\in[0,{m-1\over 2}).$ For $\dl={m-1\over 2}$ the asymptotics
is known: $\g(m)\ln n+O(1),$ E. Stein, 1961.

Let for $n\in\BN$ $$S_n(f;W)=\Sl_{k\in nW}c_k(f)e_k.$$
The growth of the \Lc\ is minimal when $W$ is polygon.
\proclaim{4.2} If $W$ is the $m$-dimensional polygon in $\BR^m,$ then
$$\sup\limits_{||f||_C\le1}||S_n(f;W)||_C=
\sup\limits_{||f||_1\le1}||S_n(f;W)||_1\asymp\ln^m n.$$ \ep
\proclaim{4.3} If for $m=2$ for the rombic partial sums
$$S_n^\diamondsuit(f)=\Sl_{|k_1|/n_1+|k_2|/n_2\le 1}c_k(f)e_k$$
we have $n_2/n_1\in\BN,$ then
$$||S_n^\diamondsuit||_{C\to C}=32\pi^{-4}\ln n_1\ln n_2-
16\pi^{-4}\ln^2n_1+O(\ln n_2).$$ \ep

Note that for $m\ge 3$ a similar result like in Theorem {\bf 4.3} apparently
is not obtained, unlike the following Theorem {\bf 4.4}.
\proclaim{4.4} Let for $m=2,$ $n\in\BN,$ and $\a\ge1$
$$H_n^\a(f)\sim\Sl_{|k_1|^\a|k_2|\le n}c_k(f)e_k.$$

{\rm a)} We have $||H_n^\a||_{C\to C}\asymp n^{1/(2+2\a)}$ with positive
constants depending only on $\a.$

{\rm b)} Let the spectrum of $f$ is outside coordinate axes and
$$D_{r,\a}(f)\sim\Sl_k(ik_1)^{r\a}(ik_2)^r c_k(f)e_k.$$   Then
$$\spl_{||D_{r,\a}(f)||_\iy\le1}||f-H_n^\a(f)||_C\asymp n^{-r+1/(2+2\a)}.$$

{\rm c)} If the matrix $A$ with the entries $\{a_ij\}_{i,j=1,2}$ cannot become
a matrix with integer entries after multiplication its lines by some
numbers, then the operator
$$H_n^A(f)\sim\Sl_{|(a_{11}k_1+a_{21}k_2)(a_{12}k_1+a_{22}k_2)|\le n}
c_k(f)e_k$$ is not bounded in $C(\BT^2)$ for $n$ sufficiently large. \ep

Let us consider now step-wise hyperbolic sums which give the approximation
of the best order for the class of functions with bounded mixed derivative
in $L_p$ with $p\in(1,+\iy).$
\proclaim{4.5} Let for $f\in L(\BT)$
$$H_n(f)=\Sl_{s_j\ge0, \sum s_j\le n}\Sl_{2^{s_j}\le|k_j|<2^{s_j+1}}
c_k(f)e_k.$$ Then $$||H_n||_{L_1\to L_1}\asymp n^{(3m-1)/2}.$$  \ep
Let us go on to general estimates of the \Lc. Let us start with
the upper bound for  the $L_p$ norm of a polynomial for $p\in(0,2).$
\proclaim{4.6}Let $\{e_j^0\}_{j=1}^n$ be the standard basis in $\BR^m,$ $M_0=
(1,...,m),$ and $q=\sum q_j e_j^0$ where the $q_j$ are natural numbers
($j\in M_0$); analogously $h=\sum h_j e_j^0$ where the
$h_j$ are also natural numbers. Set
$$\Delta_{h_j}\lambda_k=\lambda_k-\lambda_{k+h_j e_j}$$
(the difference operator with stepsize $h_j$ in the direction $e_j$) and
$$\Delta_h^q\lambda_k=\biggl(\prod\limits_{j\in M_0}\Delta^{q_j}_{h_j}
\biggr)\lambda_k$$ ("mixed" difference in the direction of all axes).
For every $p\in(0,2)$ and  $q\in\BN^m$ 
$$\align &\Il_{\BT^m}|\Sl_{-n_j\le k_j\le n_j}\lambda_k e^{i(k,x)}|^p
\,dx\\ \le \g(p,q,m)&\prod\limits_j(n_j+1)^{(p-2)/2}\Sl_{0\le s_j\le
[\log_2(n_j+1)]} 2^{(1-p/2)\Sl_j s_j}\biggl(\Sl_k|\Delta_h^q
\lambda_k|^2\biggr)^{p/2}, \endalign$$
where $\lambda_k$ is taken to equal $0$ for $k_j\ne[-n_j,n_j],$ 
at least for one $j,$ in the sum
$\sum\limits_k,$ while $h=h(s,n)$ is defined by the following conditions
$${n_j+1\over 3\cdot 2^{s_j}}\le h_j\le {5(n_j+1)\over 6\cdot 2^{s_j}},
\qquad\qquad {n_j+1\over 3\cdot 2^{s_j}}\le h_j\le {n_j+1\over 2^{s_j}}$$
according as $s_j<[\log_2(n_j+1)]$ or $s_j=[\log_2(n_j+1)].$\endproclaim

In the case $\lm_k=\vp(k/n)$ there appears the \Ft\ $\hat\vp.$
E. Belinskii was apparently the first
began a {\it systematic\/} study of connections between summability and 
integrability of the Fourier transform of a function generating a
method of summability, in the multi-dimensional case. 
\proclaim{4.7}  Let $\vp$ be a bounded
measurable function with a compact support. Then for the norms of a
sequence of linear operators 
$$L_n^\vp:\quad f\to L_n^\vp(f;\cdot)=\Sl_k\vp(k/n)\hf(k)e^{i(k,\cdot)}$$
 we have
$$\align \Vert L_n^{\vp} \Vert_{L^1 (\BT^m) \to L^1 (\BT^m)} & \le
(2\pi)^{-m} \Il_{n\BT^m} \prod \limits_{j=1}^m {x_j \over
2n\sin {(x_j/ 2n)}} |\hat\vp (x)| \, dx \\
& + \Sl_{j=1}^{r-1} (\pi/2)^{(j+1)m} \Il_{n\BT^m} |\hat\vp (x)||x/N|^j \, dx \\
+ \pi^{mr+m/2} 2^{-mr+m/2} \int \limits_{\BT^m/2\pi} & \cdots
\Il_{\BT^m/2\pi} \left( \sum_k \left| \Delta_{k/n}^r
\left( \vp ; {u_1/n},...,{u_r/n} \right) \right|^2
\right)^{1/2} du_1 ... du_r ,  \endalign$$
and
$$\g(p) \Vert L_n^{\vp} \Vert_{L^p (\BT^m) \to L^p (\BT^m)}  
\ge\left\{(2\pi)^{-m} 
\Il_{\ve n\BT^m} \left| \prod \limits_{j=1}^m {x_j \over
2n \sin (x_j/2n)} \hat\vp(x) \right|^p \, dx \right\}^{1/p} $$
$$- \Sl_{j=1}^{r-1} (\pi/2)^{(j+1)m} \left\{ \Il_{\ve n\BT^m}
|\hat\vp (x)|^p |x/n|^{jp} \, dx \right\}^{1/p}$$
$$- \dfrac{\pi^{mr+m{2-p \over 2p}}}{2^{mr-m{2-p \over 2p}}}
\frac{\ve^{m{2-p \over 2p}}}{n^{m{p-1 \over p}}} 
\Il_{\BT^m/2\pi} \cdots
\Il_{\BT^m/2\pi} \left( \sum_k \left| \Delta_{k/n}^r
\left( \fa ; {u_1/n},...,{u_r/n} \right) \right|^2
\right)^{1/2} du_1 ... du_r .$$  \endproclaim

\quad

Here $\ve\in(0,1]$ is an arbitrary real number, $r$ is integer, and 
$1 \le p \le 2.$ The $r$-th difference $\Delta_z^r (\vp; h_1 ,..., h_r)$ 
is defined recursively by the formulas
$$\Delta_z^1 (\vp; h_1) = \vp(z+h_1) - \vp (z) ;$$
$$\Delta_z^r (\vp; h_1 ,..., h_r) = \Delta_{z+h_r}^{r-1}
(\vp; h_1 ,..., h_{r-1}) - \Delta_z^{r-1} (\vp; h_1 ,..., h_{r-1}) ,$$
with $h_j , z \in \BR^r$. When $p > 2$, in view of duality 
the lower bound still
valids with $p' = p/(p-1)$ instead of $p$.
\medskip
\proclaim{4.8}  Let the boundary of the
region $B$ contain a simple (non-intersecting) piece of a surface
of smoothness $[{(m+2)/2}]$ in which there is at least one point with
non-vanishing principal curvatures. Then there exists a positive
constant $C$ depending only on $B$ such that
$$\Il_{\BT^m} \left| \Sl_{k \in nB \cap \BZ^m} e^{i(k,x)}\right| \, dx
\ge C n^{(m-1)/2} $$     for large $n$.         \endproclaim \medskip
\proclaim{4.9} For each $r>0$ and for every $\a>2r$
$$\sup\limits_{||\Dl^rf||_\iy\le1}||f-\Sl_{|k|\le n}(1-|k|^\a/n^\a)^
{{m-1\over 2}}c_ke_k||_C=\g(m,r,\a){\ln n\over n^{2r}}+O(n^{-2r}).$$\ep
Let us return to the one-dimensional case, that is $m=1.$ Let
$$\Lm_n^0(f)=\Sl_{k=-n}^n\lm_{n,k}c_k(f)e_k$$ and
$$||\Lm_n^0||=\sup\limits_{||f||_C\le1}||\Lm_n^0(f)||_C=
{1\over 2\pi}||\Sl_{k=-n}^n\lm_{n,k}e_k||_1.$$
If to replace the \Fc\ $c_k$ (integrals) by $c_k^{(n)}(f)$ by
rectangle formula for the uniform partition $x_p={2p\pi\over 2n+1},$
with $|p|\le n,$ then we get
$$\tilde\Lm_n^0(f;x)=\Sl_{k=-n}^n\lm_{n,k}c_k^{(n)}(f)e^{ikx}$$ where
$$c_k^{(n)}(f)={1\over 2n+1}\Sl_{p=-n}^n f(x_p)e^{-ikx_p}.$$
Coefficients $c_k^{(n)}(f)$ are called the Fourier-Lagrange coefficients.
When $\lm_{n,k}=1$ for all $k\in[-n,n]$ we have that $\tilde\Lm_n^0(f)$
is the interpolation polynomial defined by the values $f(x_p)$ for
$p\in[-n,n].$

The convergence of $\Lm_n^0(f)$ to $f$ at a point $x$ for every $f\in
C(\BT)$ is reduced to the boundedness in $n$ of the norms of functionals,
the Lebesgue functions.
\proclaim{4.10} We have
$$\sup\limits_{||f||_C\le1}|\tilde\Lm_n^0(f;x)|=(\pi/2)|sin(n+1/2)x|\cdot
||\Lm_n^0||+\theta\sup\limits_{||f||_C\le1}|\tilde\Lm_n^0(f;0)|
$$ where $|\theta|\le C.$  \ep

In the following theorem a trigonometric polynomial $T_n$ of order not
greater than $n$ is replaced by piece-wise sinusoidal function (see a)),
which allows to calculate the asymptotics if integral norms of $T_n$
(see b)). For this, let us introduce a sequence of functions $\vp_n$
corresponding to $\{\lm_{n,k}\}_{k=-n}^n$ and satisfying the only
condition $\vp_n(x_k)=\lm_{n,k}$ for $k\in[-n,n].$ For instance,
$\vp_n$ may be a polynomial or piece-wise linear function.
\proclaim{4.11} {\rm a)} For every $x\in[-\pi,\pi]$ and
$p=[{1\over 2}+{2n+1\over 2\pi}x]$ the following inequality holds
$$\align |T_n(x)&-{(-1)^p\over 2n+1}T'_n(x_p)\sin(n+1/2)x|\\ 
\quad\\ &\le C\Sl_{k=-n}^n|T_n(x_k)|
{(2n+1)^2\over(2|p-k|+1)^2(4n+1-2|p-k|)^2}.\endalign$$

{\rm b)} The following asymptotic equality holds
$$\Il_{\BT}|\Sl_{k=-n}^n\lm_{n,k}e^{ikx}|\,dx={4\over\pi}\Sl_{k=-n}^n
|c_k^{(n)}(x\vp_n)|+\theta\Sl_{k=-n}^n|c_k^{(n)}(\vp_n)|$$
with $|\theta|\le C.$ \ep

As an application, one can derive an asymptotics of the \Lc\ for
the general sums of Bernstein-Rogosinski type $\Sl_{k=-n}^n\mu_k
S_n(f;x+x_k)$ with the remainder term $\theta\sum|\mu_k|.$

Let us return to the problem of approximation of the class of
functions with bounded derivative.

Let $\{\psi(k)\}_1^\iy$ be a convex downwards and decrease to zero
sequence and $\be\in\BR.$ Let us denote by $W^{\psi,\be}$ the class of
continuous functions on $\BT$ for which the trigonometric series
$$\Sl_{k\ne0}e^{i\be{\pi\over 2}\sign k}{1\over\psi(|k|)}c_k(f)e_k
\sim f^{\psi,\be}$$ is the \Fs\ of the function $f^{\psi,\be},$
derivative in some sense, and $||f^{\psi,\be}||_\iy\le1.$ If
$\psi(k)=k^{-r}$ and $\be=r$ or $\psi(k)=k^{-r}$ and $\be=r+1$ with
$r>0,$ then this is the class of functions with the $r$-th derivative
in the Weyl sense $f^{(r)}$ or conjugate to it $\tilde f^{(r)}$
bounded by $1$ in both cases. The class $W^{\psi,\be}$ was introduced
and investigated by A. I. Stepanets, 1986, mainly with respect to
approximation by partial sums (see [M24, Ch.3]).

A final general result (see a)) looks as follows. Let $E_n(f)$ be
best approximation of $f$ by polynomials of order not greater than
$n$ in $C(\BT).$
\proclaim{4.12} {\rm a)} Let $n\in\BZ_+,$ and $\{\psi(k)\}_{n+1}^\iy$
and $\be$ be as above. Assume also in the case $\sin\be{\pi\over 2}\ne0$
that $$\Sl_{k=1}^\iy{1\over k}\psi(k)<\iy,$$ which is necessary.
Then treating $0/0=0$ we have
$$\align \sup\limits_{f^{\psi,\be}\in C}{||f-S_n(f)||_C\over E_n(f^{\psi,
\be})}&=\max\limits_{f\in W^{\psi,\be}}{||f-S_n(f)||_C\over E_n(f^{\psi,
\be})_\iy}=\max\limits_{f\in W^{\psi,\be}}||f-S_n(f)||_C\\ \quad\\
={4\over\pi^2}\Sl_{k=1}^{n+1}{\psi(k+n)\over k}&+{2\over\pi}|\sin\be
{\pi\over 2}|\Sl_{k=n+2}^\iy{\psi(k+n)\over k}+\theta\psi(n+1)\endalign$$
with $|\theta|\le C.$

{\rm b)} Under the same assumptions on $\psi$ and $\be$ the following
two-sided estimate holds with absolute constants:
$$\sup\limits_{f^{\psi,\be}\in C}{E_n(f)\over E_n(f^{\psi,\be})}=
\sup\limits_{f\in W^{\psi,\be}}E_n(f)\asymp\psi(n+1)+|\sin\be{\pi\over 2}|
\Sl_{k=n+1}^\iy{\psi(k+n)\over k}.$$             \ep

Let make additional three remarks.
\flushpar
In the results given in {\bf 4.11} norms in $C$ or $L_\iy$ can be replaced
by norm in $L.$ Then in b) $E_n(f^{\psi,\be})$ in the denominator can be
replaced by $\om_p(f^{\psi,\be};{\pi\over n+1}),$ the modulus of smoothness
of order $p.$ But in this case constants may be dependent on $p\in\BN.$
At last, the convexity condition for $\{\psi(k)\}$ can be substituted be
the weaker condition from {\bf 1.10}a).
\bigskip
\flushpar
{\bf Bibliographical remarks.}
\flushpar
For problems on approximation of classes of functions, see, e.g., [S15]
and [M23]. Only recently steps to make a survey on \Lc\ in the multiple
case were taken by E. Liflyand (see [L11]).  Besides the results mentioned
above, those due to V. A. Yudin, A. N. Podkorytov, M. A. Skopina, and
others, are considered in that survey.
\flushpar
For {\bf 4.1}, see [Le1,4]; for {\bf 4.2}, see [Be5];
for {\bf 4.3}, see [Ku1,2] (for $n_1=n_2$ by I. K. Daugavet, 1970);
for {\bf 4.4}, see [Be7,8] (see also [BL3]);
for {\bf 4.5}, see [Be14]; for {\bf 4.6}, see [T16];
for {\bf 4.7}, see [Be5]; for {\bf 4.8}, see [L7];
for {\bf 4.9}, see [BL2] in the case of essentially more general radial
methods; for {\bf 4.10}, see [T31] (in the case $m=2$ a similar result
for rhombic sums from {\bf 4.3} was obtained in [Ku12]);
for {\bf 4.11}, see [T10]; for {\bf 4.12}, see [T24,28] (the last inequality
in {\bf 4.12}a) was independently proved by S. A. Telyakovskii, 1987.

\head 5. Two-sided estimates of approximation. Moduli of
smoothness and $K$-functionals. \endhead
Let us start with the following result.
\proclaim{5.1} For every $r\in\BN$ there exists a continuous
function $\vp$ with compact support such that for each $p\in[1,+\iy]$
and for each $f\in L_p(\BT)$ we have the following two-sided
estimate with positive constaants depending only on $r.$
$$\align  ||f-\sum\vp(k/n)c_ke_k||_p &\asymp\om_r(f;1/n)_p\\ &
=\sup\limits_{0<\dl\le1/n}||\Sl_{\nu=0}^r
\binom r\nu (-1)^\nu f(\cdot+\nu\dl)||_p.\endalign$$  \ep

The upper bound is the well-known Jackson type theorem
($r=1 -$ D. Jackson; $r=2 -$ N. I. Akhiezer, 1947; $r\ge3 -$
S. B. Stechkin, 1951; see [M26] or [M6]). Here the same lower bound is
added. Such results allow to obtain:

1) constructive (approximate) characteristics of basic classes
of functions, which are not described by the rate of decrease of
their best approximations;

2) saturation classes of sequences of linear operators of convolution type;
\flushpar
and, as it turned out,

3) formulas of new type for $K$-functionals of a couple of spaces
(see below).

The first method of proving such inequalities is based on extremal
properties of polynomials, like Bernstein's inequalities, and use
of the Jackson theorem itself (see, e.g., [M6], p.241). On this way
necessary and sufficient (at once) conditions are found for validity
of upper bounds, lower bounds, two-sided estimates, etc. Found is also
the sharp order of approximation expressed by moduli of smoothness for
each of the classical methods of summability of \Fs, namely those of
Riemann, Weierstrass, Lebesgue, Jackson - de la Valle\'e Poussin,
Riesz, Rogosinski, Bernstein, Favard, $(C,\a),$ and Abel-Poisson.
The latter two methods are  equivalent for $\a>0$ (see {\bf 2.14} and
{\bf 2.15}a)). For these, see [T2,18]. 
Exact order of approximation by the $(C,1)$-means is
found by V. V. Zhuk, 1967. See also [M31] where this method of proof
is applied to nonlinear methods as well.

The second method is based on the comparison principle from \S2. On this way
both the Jackson theorem and the needed properties of polynomials are
obtained. Simplicity and exactness is its benefit. It was B. S. Mityagin,
1962, who first applied multipliers to problems of approximation
theory in $L_p$ for $p\in(1,+\iy).$ Observe, that the inequality in the
$C$-metric yields the inequality in $L_p$ for all $p\in[1,+\iy]$ with
the same constant (see {\bf 2.18}a)) but not vice versa. For instance,
for $p=1$ and $p=\iy$ and odd $r$ the result {\bf 5.2}b) is false.
\proclaim{5.2} {\rm a)} Let $r\in\BN$ and $\vp\in C(\BR).$ In order
that $f\in W^r$ (for the definition of $W^r,$ see the beginning of \S4)
and for every $n\in\BN$ $$||f-\sum\vp(k/n)c_ke_k||_\iy\le K n^{-r}
||f^{(r)}||_\iy,$$ it is necessary and sufficient that $g(x)=(1-\vp(x))/
x^r\in B(\BR).$ The minimal value of $K$ is $||g||_B.$

b) For every $r\in\BN,$ for every $f\in L_p(\BT)$ with $p\in(1,+\iy),$
and for every $n\in\BN$ we have the following two-sided inequality
with constants depending on $r$ and $p$
$$||f-\Sl_k(1-|k|^r/n^r)_+c_ke_k||_p\asymp\om_r(f;\pi/n)_p.$$   \ep

Let us linearize the modulus of smoothness in $C(\BT)$ as well as in $L_p
(\BT)$ as follows     $$\tilde\om(f;h)=
||{1\over h}\Il_0^h\Sl_{\nu=0}^r(-1)^\nu\binom{r}\nu f(\cdot+\nu\dl)\,d\dl||,$$
that is taking of the least upper bound in $\dl\in(0,h]$ is replaced
by the integral averaging.
\proclaim{5.3} a) For every $r\in\BN$ and for every $f\in C(\BT)$
$$\g(r)\om_r(f;h)\le\tilde\om_r(f;h)\le\om_r(f;h).$$

b) The inequality $\om_r(f;h)\le h^r$ holds if and only if the inequality
$\tilde\om_r(f;h)\le{1\over r+1}h^r$ holds. \ep

All this was the one-dimensional case. Let us go on to the multiple case.
For the Jackson type theorem, see [M16]. Let $r\in\BN,$ $E\subset\BR^m,$
and $h>0.$ Define $$\om_r(f;E;h)=\sup\limits_{u\in E}||\Sl_{\nu=0}^r
\binom r\nu(-1)^\nu f(\cdot+\nu hu)||,$$ where the norm is taken in
$C(\BT^m)$ or $L_p(\BT^m).$ For the monotonicity in $h,$ it should be
assumed that $E$ is star-like with respect to the origin. Let it also
be compact. Then the biggest, in some sense, modulus is in the case
of $E$ is the ball, the full modulus $\om_r^0;$
while the smallest one when $E$ is the closed interval
starting at the origin, the modulus of smoothness in the given direction.
Denote by $\om_r^+$ the modulus of smoothness defined by the set $E$
consisting of $m$ unit intervals along the axes of the standard basis.
\proclaim{5.4} For every $r\in\BN,$ for every $p\in(1,+\iy),$ for
every $f\in L_p(\BT^m),$ and for every $n\in\BN$ the two-sided estimate
$$||f-\Sl_k\prod\limits_{j=1}^m(1-|k_j|^r/n^r)c_ke_k||_p\asymp\om_r^+(f;
\pi/n)_p$$ holds with the constants depending on $r$ and $p.$  \ep
\proclaim{5.5} For no one $r\ge1$ and $m\ge2$ the following inequalities
$$\g_1(r,m,\vp)\om_r^+(f;\pi/n)\le||f-\sum\vp(k/n)c_ke_k||\le
\g(r,m,\vp)\om_r^0(f;\pi/n)$$ cannot valid in $C(\BT^m),$ at least if
$\vp$ is continuous at its support, the unit cube.  \ep

Thus the usual moduli in $C(\BT^m)$ and $L_p(\BT^m)$ are not suitable
for $m\ge2.$ Let us introduce linearized moduli.

Restricting to the case of even $r,$ let us set for $q\in\BN$
$$\tilde\om_{2q}(f;\mu;h)=||\Il_{\BR^m}\Sl_{\nu=0}^{2q}\binom{2r}\nu
(-1)^\nu f(\cdot+(\nu-r)hu)\,d\mu(u)||,$$ that is taking of the least
upper bound in $h$ is substituted by the integral averaging with
respect to the finite Borel measure. If $d\mu=\chi_E\,du,$ where $\chi_E$
is the indicator of $E,$ then we write $\tilde\om_{2q}(f;E;h).$ If
$E$ is the unit ball or the sum of $m$ closed unit intervals starting
at the origin in the directions of coordinate axes, then we write
$\tilde\om_{2q}^0(f;h)$ or $\tilde\om_{2q}^+(f;h)$ respectively.
\proclaim{5.6} {\rm a)} For every $q\in\BN,$ for every $\dl>{m-1\over 2},$
and for every $f\in C(\BT^m)$
$$||f-\Sl_k(1-|k|^{2q}/n^{2q})_+^\dl c_ke_k||\asymp\tilde\om_{2q}^0(f;\pi/n).$$

{\rm b)} For every $q\in\BN,$ for every $\dl>{m-1\over 2},$
and for every $f\in C(\BT^m)$ $$||f-\Sl_k(1-n^{-2q}\Sl_{j=1}^m k_j^{2q})_+^\dl
c_ke_k||\asymp\tilde\om_{2q}^0(f;\pi/n),$$ where these two-sided inequalities
hold with constants depending on $r,$ $\dl,$ and $m.$  \ep

Observe that the moduli $\tilde\om_{2q}^0$ and $\tilde\om_{2q}^+$ are
almost monotone in $h.$
\proclaim{5.7} For every $q\in\BN,$ for every $p\in(1,+\iy),$ for every
$f\in L_p(\BT^m),$ and for every $h>0$ the following two-sided estimate
$$\tilde\om_{2q}^0(f;h)_p\asymp\om_{2q}^0(f;h)_p$$ holds with constants
depending on $r,$ $m,$ and $p.$  \ep

When $m\ge2$ the cases $q=1$ and $q\ge 2$ differ for $\tilde\om_{2q}(f;E;h).$
\proclaim{5.8} {\rm a)} If $E$ is symmetric in the following sense: for
any rearrangement of each two coordinates or change of the sign of
any of the coordinates the point is still in $E,$ then for every
$p\in[1,+\iy]$ $$\tilde\om_2(f;E;h)_p\asymp\tilde\om_2^0(f;h)_p.$$

{\rm b)} For $m\ge2$ and $q\ge2$ the moduli $\tilde\om_{2q}^0$ and
$\tilde\om_{2q}^\square,$ where the latter modulus corresponds to the
unit cube, are uncomparable on $C(\BT^m)$ as $h\to 0.$  \ep

For odd $r$ and $p=\iy$ or $p=1,$ unlike the case $p\in(1,+\iy),$
in {\bf 5.2}b) the two-sided estimates from {\bf 5.1} are impossible
if $\vp$ is continuous and even function with compact support (see [T6]).

Let $\tau_n$ be the Marcinkiewicz means, that is the $n$-th arithmetic means
of the cubic partial sums.
\proclaim{5.9} For every $f\in C(\BT^2)$
$$||f-\tau_n(f)||\asymp||\Il_1^\iy(\Dl^2_{(t/n)(e_1^0+e_2^0)}+
\Dl^2_{(t/n)(e_1^0-e_2^0)})f(\cdot){dt\over t^2}||,$$ where
$\{e_1^0, e_2^0\}$ is the standard basis in $\BR^2$ and
$$\Dl_h^2 f(x)=f(x+2h)-2f(x)+f(x-2h).$$  \ep

In order to find interpolation spaces by the real interpolation method 
J. Petre, 1963, has introduced $K$-functionals (see, e.g., [M4]). 
Z. Ciesielski, 1983, put
the problem of finding of the $K$-functional via linear means of the \Fs.
\proclaim{5.10} Under the above notation for every $\ve>0$

{\rm a)} For $p\in[1,+\iy]$ we have
$$K(\ve^{2r},f,L_p,\Dl^r)\asymp\tilde\om_{2r}^0(f;\ve)_p.$$

{\rm b)} For $p\in[1,+\iy]$ we have
$$K(\ve^{2r},f,L_p,D_{2r})\asymp\tilde\om_{2r}^+(f;\ve)_p,$$
where $\Dl$ is the Laplace operator and $D_{2r}=\Sl_{j=1}^m{\p^{2r}
\over\p x_j^{2r}}.$  \ep

For $p=\iy$ and $p=1$ these formulas for the $K$-functionals of the
indicated couples of spaces were not known earlier.

Observe, that for $p=\iy$ as well as for $p=1$ it is impossible in
{\bf 5.10}b) to put $\om_{2r}^+$ in view of {\bf 5.5}. The modification
needed for this was introduced by Z. Ditzian, TAMS 1991, for $r=1$ in
the space $C.$ Let us generalize this modulus for any $r$ in the two
cases indicated. Let us set
$$\Dl_{r,\dl}^+f(x)=\Sl_{j=1}^m\Sl_{\nu=0}^{2r}\binom{2r}\nu(-1)^\nu
f(x+(\nu-r)\dl e_j^0).$$
\proclaim{5.11} For every $r\in\BN$ and for every $p\in[1,+\iy]$

{\rm a)} $\tilde\om_{2r}^+(f;h)_p\asymp\sup\limits_{0<\dl\le h}
||\Dl_{r,\dl}^+ f(\cdot)||_p,$

{\rm b)} $\tilde\om_{2r}^0(f;h)_p\asymp\sup\limits_{0<\dl\le h}
||(\Dl_{1,\dl}^+)^r f(\cdot)||_p.$   \ep

Problems dealing with hyperbolic differential operators are usually
more difficult than those for elliptic operators; see [S16].

A formula for the $K$-functional is found till now only as the
approximation by stepwise-hyperbolic Riesz means and only for $p\in(1,+\iy).$

Let $D_r={\p^{r_1+...+r_m}\over\p x_1^{r_1}...\p x_m^{r_m}}$ be the
mixed derivative with $0<r_1=...=r_\nu<r_{\nu+1}\le...\le r_m.$
Let us assume that $\Il_{-\pi}^\pi f(x)\,dx_j=0$ for each $j\in[1,m]$
and set for $s\in\BZ_+^m$
$$\dl_s(f)=\Sl_{k\in\rho(s)}c_k(f)e_k,$$
where $\rho(s)=\{k\in\BZ^m:2^{s_j-1}\le|k_j|<2^{s_j}; 1\le j\le m\}.$
Denote by $I$ the identity operator.
\proclaim{5.12} For $[1/\ve]=2^n$ and $p\in(1,+\iy)$
$$K(\ve^{r_1};f;L_p,D_r)\asymp||f-\Sl_{(s,r)\le r_1n}\biggl(
I-{1\over 2^{r_1n}}D_r\biggr)\dl_j(f)||_p.$$  \ep

the same methods are applicable also in the case of differential
operators of fractional order. Let us restrict ourselves to the case $m=1.$
\proclaim{5.13} Let us set for fractional $r>0$ and for $f\in C(\BT)$
$$\tilde\om_r(f;h)=||\Il_1^\iy(\dot\Dl_{hu}^{2p}f(\cdot)+\g
\dot\Dl_{hu}^{2p+1}f(\cdot)){du\over u^{1+r}}||,$$ where $p\in\BN$
satisfies $p>r/2,$ and $\dot\Dl_\dl f(x)=f(x+\dl)-f(x-\dl)$ and
$$\g={1\over 2}\tan(r\pi/2)\Il_0^\iy{\sin^{2p}t\over t^{1+r}}\,dt
\biggl(\Il_0^\iy{\sin^{2p+1}t\over t^{1+r}}\,dt\biggr)^{-1}.$$

{\rm a)} Setting $\vp_r(x)=(1-|x|^r)_+-i\tan(r\pi/2)|x|^r(1-|x|)_+\sign x,$
we have for every $f\in C(\BT)$ and for every $n$
$$||f-\Sl_k\vp_r(k/n)c_ke_k||\asymp\tilde\om_r(f;1/n).$$

{\rm b)} $\tilde\om_r(f;h)=O(h^r)$ as $h\to 0$ if and only if
$f^{(r)}\in L_\iy(\BT).$

{\rm c)} $K(\ve^r;f;C,W^r)\asymp\tilde\om_r(f;\ve).$  \ep

Let us give one more comparatively general result on saturation class.
\proclaim{5.14} If a positive measure $\mu$ satisfies $\int|u|^{2r}\,d\mu
<\iy$ and the support of this measure does not lie in any hyperplane
passing through the origin, then for every $\dl>{m-1\over 2}$
$$||f-\Sl_k\biggl[1-\int{(k,u)^{2r}\over n^{2r}}\,d\mu(u)\biggr]_+^\dl
c_ke_k||_\iy=O(n^{-2r})$$ if and only if  $$\int\biggl(\Sl_{j=1}^m
u_j{\p\over\p x_j}\biggr)^{2r}f(x)\,d\mu(u)\in L_\iy(\BT^m).$$ \ep

Let us consider now the problem on approximation by the Bernstein method
of the class of functions with given modulus of continuity in $C(\BT)$
in the following form
$$A\om(f;\pi/n)\le||f(\cdot)-{S_n(\cdot)+S_n(\cdot+\pi/n)\over 2}||
\le B\om(f;\pi/n).$$

The upper bound was known even to S. N. Bernstein; for the lower
bound, see [T2] or [M6], p.241). The smallest constant $B$ in this inequality
under sufficiently general conditions on the modulus of continuity is
found by V. T. Gavrilyuk and A. I. Stepanets, 1973 (see also [M23]).
In the following statement the exact biggest constant $A$ in approximation
from below is found for the first time.
\proclaim{5.15} For every modulus of continuity the exact constant in
the above given inequality is $$A=\biggl(2+{4\over\pi}\Il_0^\pi
{\sin t\over t}\,dt\biggr)^{-1}.$$   \ep

A similar result in the lower bound is obtained also for the
Rogosinski sums ${1\over 2}(S_n(\cdot+\pi/2n)+S_n(\cdot-\pi/2n))$
when $\om(f;\pi/n)$ is replaced by $\om_2(f;\pi/2n)$ in the above
inequality.
\bigskip
\flushpar
{\bf Bibliographical remarks.}
\flushpar
Two-sided estimates of approximation {\bf 5.1} (by the Bernstein sums,
by the Rogosinski sums and the like) were firstly obtained by the author;
see [T2] or [M6].
\flushpar
For {\bf 5.3}, see [T3,15]; for {\bf 5.4}, see [No2];
for {\bf 5.5}, {\bf 5.6}, {\bf 5.8}, and {\bf 5.14}, see [T17]
(the case $r=1$ in {\bf 5.6}a) was obtained earlier in [Be3]);
for {\bf 5.7} and {\bf 5.10}, see [T23]; the result {\bf 5.9} was
obtained by O. I. Kuznetsova (see [KT]).  For $p\in(1,+\iy)$ the rate of
approximation by the J. Marcinkiewicz means was found by M. F. Timan and
V. G. Ponomarenko in 1975.
\flushpar
For {\bf 5.12}, see [Be19]; for {\bf 5.13}, see [T32], see also [KT];
for {\bf 5.15}, see [Kl].
\flushpar
General theorems of Voronovskaya type and converse to them are also
obtained by the multiplier method, see [T2], [Ku5,7]. See also [T18]
or [T27].

\head 6. Hardy spaces $H_p.$ \endhead
Let us consider the same problems as in \S\S 3 and 5 for the
Hardy spaces $H_p$ immediately in the multidimensional case
and the first with $p\in(0,1].$ Our investigation is based on
the multiplier method.

Let $D^m=\{z=(z_1,...,z_m): |z_j|<1, 1\le j\le m\}$ be the unit
polydisc in $\BC^m.$ Every function in $H_p(D^m),$ where $p>0,$
is expanded in $D^m$ in the absolutely convergent power series
$$f(z)=\Sl_{k\in\BZ_+^m}c_kz^k,$$ where $z^k=z_1^{k_1}...z_m^{k_m}$ and
$$||f||_{H_p}=||f(\cdot)||_{H_p}=\sup\limits_{0<r_j<1}\biggl(
\Il_{-\pi}^\pi\,du_1...\Il_{-\pi}^\pi|f(r_1e^{iu_1},...,r_me^{iu_m})|^p
\,du_m\biggr)^{1/p}<\iy.$$ For  $p\ge1$ this is a subspace of $L_p(\BT^m),$
namely every function from $H_p(D^m)$ is characterized in this case by
the fact that its limit function on $\BT^m,$ the skeleton of the boundary of
$D^m,$ belongs to $L_p$ and has the \Fs\ of the power type, that is with
the spectrum in $\BR^m_+$ (see, e.g., [M17]). Thus many results from
\S\S 2-3 and 5 are still valid in the same form as particular cases
(for instance, {\bf 5.1, 5.3, 5.6, 5.9, 5.10}, and the like). For
$p\in(1,+\iy)$ also the Riesz projector theorem as well as more general
Marcinkiewicz' multiplier theorem (see [M21]) help. But in $C$ as well as in
$L_1$ results may differ. Let us give an example (one-dimensional).

For every $f\in H_\iy(D)\cap C(\bar D)$
$$\max\limits_{z\in\bar D}|f(z)-\Sl_{k=0}^n(1-{k\over n+1})c_kz^k|
\asymp\om(f;1/n)_\iy,$$ and for $f\in C(\BT),$ in assumption that also
its trigonometrically conjugate $\tilde f\in C(\BT),$ we have only
$$||f-\s_n(f)||_\iy+||\tilde f-\s_n(\tilde f)||_\iy\asymp\om(f;1/n)_\iy
+\om(\tilde f;1/n)_\iy.$$
For the second relation, see [T4], while the first one immediately
follows from the second. Another situation is for $p\in(0,1).$ A function
from the quasinormed space $L_p$ is not expanded into the \Fs\ if it
is not from $L_1.$ Moreover, as it is well-known, there are no in $L_p$
for $p\in(0,1)$ nonzero linear continuous functionals. But in $H_p$
there are Tailor series for all $p>0$ and thus there are also
multipliers. It is possible to introduce multipliers in $L_p,$ though:
first on the dense set of polynomials and then to extend continuously.

A number sequence $\{\lm_k\}_{k\in\BZ_+^m}$ is called the multiplier
in $H_p(D^m),$ written $\{\lm_k\}\in M_p,$ if for every $f\in H_p(D^m)$
with the Tailor coefficients at zero $\{c_k\}$ we have
$$\Lm f)(z)=\Sl_{k\in\BZ_+^m}\lm_k c_kz^k\in H_p(D^m)$$ and there
exists a positive constant $\g$ such that for every $f\in H_p(D^m)$
$$||\Lm f||_{H_p}\le\g||f||_{H_p}$$ and $$||\{\lm_k\}||_{M_p}=\inf\g.$$
\proclaim{6.1} Let $0<p<q\le1<r\le\iy.$ Then $M_p\sbt M_q\sbt M_r$ and
this embedding is continuous. \ep

For $p\in(1,+\iy)$, as it is follows from the Riesz projector theorem,
a sequence $\{\lm_k\}_{k\in\BZ_+^m}$ is the multiplier in $H_p(D^m)$
if and only if its extension by zero on $\BZ^m$ is the multiplier of
\Fs\ in $L_p(\BT^m).$ For $p=\iy$ and $p=1$ such extensions do exist
(see {\bf 2.4} and [M7, II, 16.7.5], respectively) but it is difficult
to describe them.

For $p\in(0,1],$ let us start with the finite sequence.
\proclaim{6.2} For every $p\in(0,1]$ and every $N\in\BN^m$
$$||\{\lm_k\}_{0\le k\le N}||_{M_p}\le\g(m,p)\biggl(\Il_{\BT^m}
|\Sl_{-N\le k\le N}\lm_ke^{i(k,u)}|^p\,du\biggr)^{1/p},$$
where $\lm_k$ may be arbitrary for $k\in[-N,N]\sm[0,N]$ since $\lm_k=0$
for $k\not\in[0,N].$   \ep

Let $\vp:\BR^m_+\to\BC.$ We will write $\vp\in M_p$ if
$$||\vp||_{M_p}=\sup\limits_{\ve>0}||\{\vp(\ve_k)\}||_{M_p}<\iy.$$
\proclaim{6.3} {\rm a)} Let $\vp\in C(\BR^m)$ and $\supp\vp\sbt[-N,N].$
If for some $p\in(0,1]$ we have $\hat\vp\in L_p(\BR^m),$ then
$$||\vp||_{M_p}\le\g(m,p)(\prod\limits_{j=1}^m N_j)^{1/p-1}(\Il_{\BR^m}
|\hat\vp(x)|^p\,dx)^{1/p}.$$

{\rm b)} If $\vp\in M_p,$ then for every inner point from $\BR^m_+$
there exists a function $psi$ with compact support, which coincides with
$\vp$ in some neighborhood of this point, and $\hat\psi\in L_p(\BR^m).$\ep

As it follows from {\bf 1.8} now, each function with compact support
$\vp\in C^\a(\BR^m_+)$ for $\a>m/p-m/2$ belongs to $M_p.$ P. Oswald, 1982,
has proved by another method in essence that if $\vp\in C^\iy(\BR_+)$
equals to $1$ on $[0,1]$ and $\vp$ vanishes on $[2,+\iy),$ then
$\vp\in M_p$ for every $p\in(0,1)$ (see [S8]).
\proclaim{Example} Let $\vp(x)=(1-|x|^{2r})_+^\dl,$ where $r\in\BN.$
It belongs to $M_p$ with $p\in(0,]$ if and only if $\dl>m/p-(m+1)/2.$\ep
\proclaim{6.4} The spherical Bochner-Riesz means
$$\Sl_{k\in\BZ_+^m}(1-|k|^{2r}/n^{2r})^\dl_+c_ke_k,$$ where $r\in\BN$
and $\dl\ge0,$ are regularin $H_p(D^m)$ with $p\in(0,1]$ if and
only if $\dl>m/p-(m+1)/2.$  \ep

Let us give now sufficient conditions on $\vp$ near infinity.
\proclaim{6.5} Let $p\in(0,1]$ and $\vp\in C^r(\BR^m_+)$ for some
$r>m(1/p-1/2).$ If also $$|\vp(x)|\le A (1+|x|^\a)^{-1},\qquad
\Sl_{j=1}^m|{\p^r\vp\over\p x_j^r}(x)|\le B (1+|x|^\be)^{-1},$$
where $\a>0$ and $\be=r+\a$ or $\a=\be>m(1/p-1/2),$ then
$$||\vp||_{M_p}\le\g(m,p,r,\a,\be)(A+B).$$   \ep

Let us go on to two-sided estimates of approximation by the Bochner-Riesz
means of power series and to moduli of smoothness.

Note that for $m=1$ two-sided estimates of approximation by the
Abel-Poisson means in $H_p(D)$ with $p\in(0,1]$ were obtained earlier
by E. A. Storozhenko, 1982 (see [S8]), and those by the Riesz means
are due to P. Oswald, 1984 (upper bound) and to E. Belinskii, 1993
(lower bound). Another methods were used for these. When $m=1$ no
special moduli of smoothness are needed (see also {\bf 6.9} below).

Let us set $$(\Dl f)(z)=\Sl_{j=1}^m(z_j{\p\over\p z_j})^2f(z),$$
which is the Laplace operator on $\BT^m,$ and introduce a linearized
modulus of smoothness
$$\tilde\om^0_{2r}(f;h)_p=||\Il_{|x|\le1}dx\Il_{|y|\le1}dy...\Il_{|w|\le1}
\Sl_{\nu=0}^{2r}\binom{2r}{\nu}(-1)^\nu f((\cdot)e^{ih(\nu-r)(x+y+...+w)}
dw||_{H_p}.$$ Here $r\in\BN,$ $h>0,$ and the integral averaging is taken
over the Carthesian product of $q$ unit balls in $\BR^m$ with $q>2m(1/p-
1/2)/(m+1).$
\proclaim{6.6} For every $r\in\BN,$ for every $p\in(0,1],$ for every
$f\in H_p(D^m),$ and for every $\ve>0$

{\rm a)} $||f(\cdot)-\Sl_k(1-\ve^{2r}|k|^{2r})_+^\dl c_k(\cdot)^k||_{H_p}
\asymp\tilde\om_{2r}^0(f;\ve)_p$ for $\dl>m/p-(m+1)/2,$

{\rm b)} $K(\ve^{2r};f;H_p(D^m),\Dl^r)\asymp\tilde\om_{2r}^0(f;\ve)_p,$
\flushpar
where the two-sided estimates are with constants independent of $f$ and
$\ve.$   \ep

Let us set now $$D_{2r}(f;z)=\Sl_{j=1}^m(z_j{\p\over\p z_j})^{2r}f(z)$$
and for $h>0$
$$\tilde\om^+_{2r}(f;h)_p=||\Il_{||x||_0\le1}dx\Il_{||y||_0\le1}dy...
\Il_{||w||_0\le1}\Sl_{\nu=0}^{2r}\binom{2r}{\nu}(-1)^\nu f((\cdot)e^
{ih(\nu-r)(x+y+...+w)}dw||_{H_p}.$$ Here the integral averaging is taken
over Carthesian product of $q$ crosses in $\BR^m$ with $q>2m(1/p-
1/2)/(m+1),$ where the cross is the sum of all unit intervals along
the coordinate axes with center at the origin.
\proclaim{6.7} For every $r\in\BN,$ for every $p\in(0,1],$ for every
$f\in H_p(D^m),$ and for every $\ve>0$

{\rm a)} $||f(\cdot)-\Sl_k(1-\ve^{2r}\Sl_{j=1}^m k_j^{2r})_+^\dl c_k
(\cdot)^k||_{H_p}\asymp\tilde\om_{2r}^+(f;\ve)_p$ for $\dl>m/p-(m+1)/2,$

{\rm b)} $K(\ve^{2r};f;H_p(D^m),D_{2r})\asymp\tilde\om_{2r}^+(f;\ve)_p,$
\flushpar
where the two-sided estimates are with constants independent of $f$ and
$\ve.$   \ep

Let us give in the multi-dimensional case the Hardy-Littlewood type
inequalities; for $m=1$ a similar inequality  for every
$r\in\BN$ is due to E. A. Storozhenko, 1982 (see [S8]).
\proclaim{6.8} For every $r\in\BN,$ for every $f\in H_p(D^m),$ and
for every $\dl\in[1/2,1)$

{\rm a)} $||\Dl^rf(\dl(\cdot))||_{H_p}\le\g(r,m,p)(1-\dl)^{-2r}
\tilde\om_{2r}^0(f;1-\dl)_p,$

{\rm b)} $||D_{2r}f(\dl(\cdot))||_{H_p}\le\g(r,m,p)(1-\dl)^{-2r}
\tilde\om_{2r}^+(f;1-\dl)_p.$   \ep

The Hardy space is also studied on the ball $\{z:\Sl_{j=1}^m|z_j|^2<1\}$
in $\BC^m$ (see [M18]). Recently Vit. Volchkov [ViV] proved {\bf 6.2}
for this case, that is for Reinhardt domains,
and consequently all what follows from this. Let us return
to the one-dimensional case.
\proclaim{6.9} If $E_n(f)$ is best approximation by polynomials of
degree not higher than $n,$ then for every $N\in\BN$
$$\ln^{-1}N\Sl_{n=1}^N n^{-1}||f(\cdot)-\Sl_{k=0}^n(1-k^2/n^2)^{1/p-1}
c_k(\cdot)^k||_{H_p}\asymp\ln^{-1}N\Sl_{n=1}^N n^{-1}E_n(f)_{H_p},$$ where
the constants in the two-sided inequality depend only on $p\in(0,1].$ \ep

Let us consider different moduli of smoothness in the disc $D.$ For
moduli of smoothness in $C$ for analytic functions, see e.g., [M26].
We will study $H_p(D)$ for all $p\in(0,+\iy).$ The contour (boundary)
modulus for $f\in H_p(D)$ is defined as follows:
$$\om_r(f;h)_p=\sup\limits_{0<\dl\le h}||\Sl_{\nu=0}^r\binom r\nu(-1)^\nu
f((\cdot)e^{i\nu\dl})||_{H_p}.$$
One can replace $f$ by the limit function $f(e^{it})$ and integrate
in $L_p$ over the circle $\BT=\p D.$

Let us introduce also a linearized boundary modulus (averaged contour
modulus) for $q\in\BN$
$$\tilde\om_r(f;h)_p=||\Il_{[0,1]^q}\Sl_{\nu=0}^r\binom r\nu(-1)^\nu
f((\cdot)e^{ih\nu\Sl_1^q u_j})\,du||_{H_p}$$ as well as radial one
for $h\in(0,2/r)$
$$\om_r(f;\text{rad};h)_p=\biggl(\Il_{-\pi}^\pi|\Sl_{\nu=0}^r\binom r\nu
(-1)^\nu f(e^{it}(1-h\nu))|^p\,dt\biggr)^{1/p}.$$
This modulus annulate all polynomials of degree not higher than $r-1.$
\proclaim{6.10} For every $r\in\BN,$ for every $p\in(0,+\iy],$ for every
$f\in H_p(D),$ and for every $h\in(0,1/(r+1)]$ we have

{\rm a)} when $q=1$ and for every $p\ge1$ or for $q=1+[1/p-1/2]$ and $p\in
(0,1)$ $$\om_r(f;h)_p\asymp\tilde\om_r(f;h)_p;$$

{\rm b)} when $S_{r-1}(z)=\Sl_{k=0}^{r-1}f^{(k)}(0)z^k/k!$ for $r\ge2$
and $S_0\equiv 0$ $$\om_r(f;\text{rad};h)_p\asymp\om_r(f-S_{r-1};h)_p;$$
where the constants in the two-sided inequalities depend only on
$r$ and $p.$  \ep

Let us give also the Bernstein type inequality in $L_p$ for $p\in(0,1)$
when the case of fractional $r$ differs from that for $p\ge1.$
\proclaim{6.11} For every $r>0$ and for every $p\in(0,1)$
$$\sup\limits_{||T_n||_p\le1}||T^{(r)}_n||_p\asymp n^r;\quad n^{1/p-1};
\quad n^{1/p-1}\ln^{1/p}n$$ for $r\in\BZ_+$ or for $r>1/p-1,$ and $r\not
\in\BZ_+$ so that $r<1/p-1,$ and $r=1/p-1\not\in\BZ_+,$ respectively. \ep
\bigskip
\flushpar
{\bf Bibliographical remarks.}
\flushpar
For {\bf 6.1 - 6.6, 6.10}, see [T35]; {\bf 6.7 - 6.8} is due to
E. M. Klebanov (submitted for publication); for {\bf 6.9}, see [Be19];
for {\bf 6.11}, see [BL4].
\head 7. Positive definite functions and splines. \endhead
Let $E$ be a linear (vectorial) space over the field $\BR.$ A function
$f: E\to\BC$ is called positive definite on $E,$ written $f\in\Phi(E)_,$
if for any set of elements $\{x_k\}_1^n$ and for any set of numbers $\{\z_k\}_
1^n$ we have $$\Sl_{k,s=1}^n f(x_k-x_s)\z_k\bar\z_s\ge0.$$
If $E$ is the Hilbert space, then one can take, for example,
$f(x)=e^{i(x,y)}$ with $y\in E$ (see also {\bf 7.14} below).
\proclaim{7.1} For every $f\in\Phi(E)$ and for every $x,y\in E$
$$|f(x+y)-2f(x)+f(x-y)|\le 2{\Re}(f(0)-f(y)).$$  \ep
This simple necessary condition immediately yields the fact that
$e^{-|x|^\a}\in\Phi(\BR)$ only for $\a\in[0,2].$

For functions in $\Phi(\BR^m),$ which are the characteristic functions
if $f(0)=1,$ see [M15] for $m=1$ and [M1] for any $m.$
\proclaim{7.2} In order that $f\in\Phi(\BR)$ it is necessary and
sufficient that the following three conditions are satisfied:

{\rm a)} $f\in C(\BR)$ and bounded,

{\rm b)} the following improper integral converges $\Il_0^\iy(f(t)-
f(-t))/t\,dt,$

{\rm c)} $\Ll_{T\to+\iy}(2T)^{-1}\Il_{-T}^T f(t)\,dt\ge 0$ and there
exists $k_0\in\BN$ such that for every $k\ge k_0$ and for every $x\ne0$
$$(\sign x)^{k+1}\Il_{-\iy}^\iy(x+it)^{-k-1}f(t)\,dt\ge0.$$  \ep

This new criterion proved to be fruitful during the proof of {\bf 7.5}b)
(see below).

Let us go on to the multiple case and radial functions. We will write
$f\in\Phi_m$ if $f:[0,+\iy)\to\BR$ and $f(|x|)\in\Phi(\BR^m).$
It is of certain interest that $\Phi(\BR^m)\sbt\Phi(\BR^{m+1})$ but
$\Phi_{m+1}\sbt\Phi_m$ (see [M1, Ch.5, \S4]).
\proclaim{7.3} {\rm a)} Let $1\le n<m.$ In order that $f\in\Phi_m$
it is necessary and sufficient that there exists a function $g\in\Phi_n$
uniquely defined and such that
$$f(t)=\Il_0^1 g(ut)u^{n-1}(1-u^2)^{(m-n)/2-1}\,du.$$

{\rm b)} Let $m\ge3$ and odd. In order that $f\in\Phi_m$ it is necessary
and sufficient that $${d^k\over dt^k}\{t^{(m-2)/2}f_0(\sqrt{t})\}_{t=0}=0$$
for $0\le k\le(m-3)/2$ and the function $f_1,$ defined for $t\ge0$ by
the formula $$f_1(\sqrt{t})=\sqrt{t}{d^{(m-1)/2}\over dt^{(m-1)/2}}
\{t^{(m-3)/2} f_0(\sqrt{t})\},$$ belongs to $\Phi_1.$  \ep

For necessary conditions for membership in $\Phi_m,$ see [T21]. Let
us give a simple sufficient condition which is Polya's theorem when $m=1.$
\proclaim{7.4} Let $n=[(m+2)/2].$ If $f\in C[0,+\iy),$ $\Ll_{t\to+\iy}
f(t)\ge0,$ $f\in C^{n-1}(0,+\iy),$ $(-1)^{n-1}f^{(n-1)}$ is convex
downwards on $(0,+\iy)$ and $$\Ll_{t\to+0}t^n f^{(n)}(t)=
\Ll_{t\to+\iy}t^n f^{(n)}(t)=0$$ where, for example, the right derivative
may be taken, then $f\in\Phi_m.$ \ep

Shoenberg's $B$-splines are very well-known (see [M30]), namely
$B_0$ is the indicator function of the interval $(-1/2,1/2),$ and
for $n\ge1$ we have $B_n=B_{n-1}*B_0$ (convolution). Evidently,
$B_n$ is a spline of compact support and of degree $n$ which belongs
to $C^{n-1}(\BR)$ and sewn by $n+1$-th algebraic polynomial. It is
obvious also that for odd $n$ the Fourier transform $\hat B_n(x)\ge0$
for each $x\in\BR$ and consequently $B_n\in\Phi(\BR).$

Let us consider compact splines in $\Phi(\BR)$ sewn only by two
polynomials, the simplest in some sense. More precise they are of
kind $p_n(|x|)$ on $[-1,1],$ where $p_n$ are real algebraic 
polynomials of degree $n.$ An example for $n=1$ is given by $(1-|x|)_+.$
\proclaim{7.5} Let us give two examples of sequences of such splines
with $n=1,2,....$

{\rm a)} $\tilde e_{2n+1}=(-1)^n P_n*P_n,$ where $P_n$ is Legendre's
polynomial for the interval $[-1/2,1/2]$ extended to $\BR$ by zero.

{\rm b)} $e_n(\sqrt{t})=\sqrt{t}{d^{n-1}\over dt^{n-1}}\{t^{n-3/2}
(1-\sqrt{t})_+^n\}.$   \ep

A problem of splines of highest smoothness arises.
\proclaim{7.6} For each $n\ge2$ there exists a unique spline of
indicated type $A_{3n-2}$ of degree $3n-2$ in $C^{2n-2}(\BR)$
satisfying $A_{3n-2}(0)=1.$ For every $x\in\BR$ therewith $\hat A_
{3n-2}(x)>0.$ The graph of $A$-spline is, for any $n,$ of the bell form,
that is increases on $[0,1]$ with one point of inflection. \ep

There exist explicit formulas for $A$-splines. For example,
$A_4(x)=(1-|x|)_+^3(1+3|x|).$ 
\medskip
\flushpar
By Wiener's theorem linear combination
of shifts of $A$-spline are dense in $L(\BR).$ Studied is the question
of degree of convergence of $B$-splines, see [M30]. The \Ft\ $\hat B$
has zeros and shifts are supplemented by contractions.
\proclaim{7.7} Let $n\ge2$ and $A=A_{3n-2}.$ If a function $f$ is
such that, for some $r\in\BN$ and for some $\lm>1,$ both $f(x)$ and
$x^rf(x)\in C(\BR)\cap H_1^\lm(\BR),$ then for each $h\in(0,1)$
there exist $N\asymp(1/h)^{1+2n/r}$ and $\{c_k\}_N^N$ satisfying
the inequality
$$\sup\limits_{x\in\BR}|f(x)-\Sl_{k=-N}^N c_kA(x+kh)|\le Kh^{\lm-1},
\qquad Kh^{2n}\ln(1/h),\quad\text{and}\quad Kh^{2n},$$ for the cases
$\lm<2n+1,$ $\lm=2n+1,$ and $\lm>2n+1,$ respectively. A constant $K$
depends only on $\lm,$ $r,$ $n,$ and the norms of two functions
in the Nikolskii space $H_1^\lm(\BR).$  \ep

Note, that the change of degree of smoothness in the space $L,$
that is $\lm$ is replaced by $\lm-1,$ is based on the embedding
theorem from $L$ into $C.$

To apply splines to the finite element method it should be possible
to represent any polynomial of corresponding degree as a linear combination
of shifts of the spline.
\proclaim{7.8} Let $n\ge2$ and $A=A_{3n-2}.$

{\rm a)} For every $\ve\in(0,1),$ for different points $\{h_k\}_0^{3n-2}$
from $[0,\ve],$ and for any polynomial $p$ of degree $3n-2$ there exists
a unique representation of the form
$$p(x)=\Sl_{k=0}^{3n-2}c_k A(x+h_k)$$ with $x\in[0,1-\ve].$

{\rm b)} If $n\ge3$ and $b-a>1,$ then the fact $\sum c_k A(x+h_k)\in
C^{2n+1}[a,b]$ yields that this sum vanishes on $[a,b].$

{\rm c)} In the case $n=2$ for any $\ve\in(0,1)$ and for any $p$ of degree
$n\le4$ there exists a representation of the form
$$p(x)=\Sl_{k=-N}^N c_k A_4(x+kh)$$ where $x\in[0,2-\ve]$ and $h$ and
$N$ depend only on $\ve.$ This is not so for $\ve=0.$   \ep

In the convex set of splines of the type in question of fixed degree
and equal $1$ at the origin, there exist extreme points. For example,
$\lm e_n$ and $A_{3n-2}$ are the extreme points in the set of splines
of the same degree. But $\lm e_n$ is not an extreme point in the set
of splines of degree not higher than $n+1,$ while $\mu\tilde e_{2n+1}$
is the extreme point even for all the (characteristic) functions in $\Phi_1$
of support $[-1,1].$
\proclaim{7.9} For every odd $m\ge3$ and any $n>(m+1)/2$ there exists
non-trivial spline of type $p_n(|x|)$ for $|x|\le1$ and vanishing for
$|x|\ge1$ which belongs to $\Phi(\BR^m)\cap C^r(\BR^m),$ where
$r=2[(2n-m-1)/6],$ and this degree of smoothness is maximal. \ep

In connection with one problem in approximation theory (see {\bf 7.11}
below) arose a question of possible extension of a function defined,
outside the interval $(-a,a),$ to a function in $\Phi(\BR).$ If such
a continuation does exist, then it is desirable that its value at
the origin is minimal.
\proclaim{7.10} {\rm a)} If $f\in\Phi(\BR)$ and $|\Re f|
+\Re f\in L(\BR),$ then also $\Re f\in L(\BR).$

{\rm b)} If $f\in L[a,+\iy),$ locally absolutely continuous, and there
exists $\ve>0$ such that for $h\to+0$
$$\Il_a^\iy|f'(x)-f'(x+h)|\,dx=O(\ln^{-1-\ve}(1/h),$$
then $f$ is extendable on $\BR$ so that it will belong to $\Phi(\BR).$
For $\ve=0$ a similar statement is false.

{\rm c)} If on $[a,+\iy)$ we have $f(x)=cx^{-r},$ where $c\in\BC$ and
$r>0,$ then an extension into $\Phi(\BR)$ is possible only in the
following two cases: {\rm 1)} $r>1$ and $c\in\BC$ is arbitrary,
{\rm 2)} $r\in(0,1]$ and $|\operatorname{arg}c|\le r\pi/2.$  \ep
\proclaim{7.11} Let $r$ be an arbitrary positive number.

{\rm a)} There exists a function $\vp_r,$ supported on $[0,1]$ and
depending only on $r$ and $m,$ such that
$$||(-\Dl)^{r/2}f||_\iy\le M\quad\text{if and only if}\quad
||f-\Sl_k\vp_r(|k|/n)c_k e_k||_\iy\le Mn^{-r}(r+1)...(r+p)/p!$$
for every $n\in\BN,$ where $\Dl$ is the Laplace operator and $p$
is the least integer satisfying the inequalities $p\ge(m-1)/2$ and
$p\ge(r+1)/2.$

{\rm b)} There exists a function $\vp_r\in C[-1,1]$ and a constant
$\g(r)$ such that
$$||f^{(r)}||_\iy\le M\quad\text{if and only if}\quad
||f-\Sl_{k=-n}^n\vp_r(k/n)c_k e_k||_\iy\le\g(r) Mn^{-r}$$
for every $n\in\BN.$  \ep

Let $E$ be a linear normed space. When $f(||x||)\in\Phi(E)?$
This question is connected with the problem of isometric embedding
of $E$ into some $L_p$ space (see [S9]).
\proclaim{7.12} If there exist three linear independent elements
$a_1,$ $a_2,$ and $a_3$ from $E$ such that
$$||a_1+a_2y_1+a_3y_2||^{-1}\{{\p\over\p t}||ta_1+a_2y_1+a_3y_2||\}_
{t=1}\in L_1(\BR^2),$$ then $f(||x||)\in\Phi(E)$ only in the case
$f\equiv\operatorname{const}.$  \ep
\proclaim{7.13} The following spaces satisfy the condition {\bf 7.12}:

{\rm 1)} $L_p(\Si),$ where $\Si$ is a space with measure either
finite or infinite, when $\operatorname{dim}L_p(\Si)\ge3$ and $2<p\le\iy.$
In particular, $l_p^m$ with $m\ge3$ and $2<p\le\iy.$

{\rm 2)} The space $C$ of continuous functions on metric space
consisting of not less than three points.  \ep
Investigated are in the same way the spaces $E$ with $\operatorname{dim}=2.$

In the next statement answered is one question due to Shoenberg, 1938.
\proclaim{7.14} For $E=l_p^m$ with $m\ge2$ and $p\in(2,+\iy]$
$$e^{-||x||_p^\a}\in\Phi(\BR^m)\quad\text{if and only if}\quad
\cases m=2, &\a\in[0,1]\\ m\ge3, &\a=0.\endcases$$   \ep
\bigskip
\flushpar
{\bf Bibliographical remarks.}
\flushpar
For {\bf 7.1-7.3, 7.5}b), {\bf 7.11}a), see [T21];
for {\bf 7.2, 7.4, 7.5}b), see [T26]; for {\bf 7.6, 7.7, 7.9}, see
the abstarct by R. M. Trigub at the International Conference in 
Constructive Theory of Functions, Sophia, 1987,
and [ZT1]; for {\bf 7.8}, see [T33].
\flushpar
{\bf 7.10, 7.11}b) are due to V. P. Zastavny (see [ZT1]); for {\bf 7.12-7.14},
see [Z6-8]. Shoenberg's problem (see {\bf 7.14}) was independently
solved in other way by Koldobskii (Algebra i Analiz (St. Petersburg
Math. J.), 1991).

\head 8. Multiplicative Walsh system. \endhead
This system was introduced by J. Walsh in 1923 as an orthonormal system
of functions complete in $L_2[0,1]$ and such that each of the
functions in this system takes only the meanings $+1$ and $-1$
(excluding some dyadic rational points at which it vanishes). One of the 
applications is the coding theory (see [M9]).

Here we consider the same questions we have studied earlier for the
trigonometric system. Besides the natural similarity, one can
find some special peculiarities in this.

To each number $x\in\BR_+$ a dyadic expansion corresponds 
$$x=\Sl_{k=0}^\iy\th_{-k}(x)2^k+\Sl_{k=1}^\iy\th_k 2^{-k},$$
where $\th_k(x)=0$ or $1;$ the first sum is finite and for dyadic
rational numbers the second sum is chosen to be finite as well.

Consider the Walsh system $\{\psi_n\}_0^\iy$ in Paley's denumeration:
$$\psi_0\equiv 1,\quad\psi_{2^k}(x)=(-1)^{\th_{k+1}(x)},\quad
\psi_n(x)=\prod\limits_{k\ge0}(\psi_{2^k}(x))^{\th_{-k}(n)}=
\prod\limits_{k\ge0}(-1)^{\th_{k+1}(x)\th_{-k}(n)}.$$
Let us introduce the following group operation on the interval $[0,1]$
$$x\dot+y=\Sl_{k=1}^\iy2^{-k}(\th_k(x)+\th_k(y))(\operatorname{mod}2).$$
Then $\psi_n(x\dot+y)=\psi_n(x)\cdot\psi_n(y).$

The Fourier-Walsh series of a function $f:[0,1]\to\BR$ and $f\in L_1[0,1]$
(the Lebesgue measure is invariant with respect to the operation $\dot+$)
may be written as follows:
$$f\sim\Sl_{k=0}^\iy c_k(f)\psi_k,\qquad c_k(f)=\Il_0^1f(t)\psi_k(t)\,dt.$$
The following theorem is an analog of Baire's classical theorem on
superposition for the trigonometric system (for the modern approach to
this question, see [S11]).
\flushpar
\proclaim{8.1} For every $f\in C[0,1]$ there exists a function $\vp\in C[0,1],$
which is strictly increasing from $0$ to $1,$ so that for $F=f\circ\vp$
we have:

\flushpar a) $\Sl_{k=0}^nkc_k(F)=o(n)$ as $n\to\iy;$

\flushpar b) The Fourier-Walsh series of $F$ converges to it uniformly. \ep
\flushpar
\proclaim{8.2} In order the Bernstein-Rogozinski method
$$B_n(f;x,\a,\be,\nu)=\a S_n(f;x)+\be S_n(f;x\dot+{\nu\over n}),$$
where $S_n$ is the partial sum of the Fourier-Walsh series of $f,$
and $\a,\be\in\BR,$ $\nu\in\BR_+,$ $n\in\BN,$ to be regular in $C[0,1],$
it is necessary and sufficient that $\a=\be=1/2$ and $\nu=1.$\ep
\flushpar
A similar result is obtained also in the multidimensional case $C[0,1]^2.$

Let $Q_k=[0,k-1]^m,$ where $k\in\BN$ and $Q_0=\emptyset,$ 
is the cube in $\BR^m.$
\flushpar
\proclaim{8.3} If $\lim\lm_k=0$ as $k\to\iy$ and $\Sl_{k=0}^\iy
\max\limits_{s\ge k}|\lm_s-\lm_{s+1}|<\iy,$ then the series
$$\Sl_{k=0}^\iy\lm_k\Sl_{n\in Q_{k+1}\setminus Q_k}\psi_n(x),\qquad
\psi_n(x)=\psi_{n_1}(x_1)...\psi_{n_m}(x_m),$$
converges a.e. on $[0,1]^m$ to an integrable function, is the Fourier-Walsh
series of its sum and
$$\Il_{[0,1]^m}\biggl|\Sl_{k=0}^\iy\lm_k\Sl_{n\in Q_{k+1}\setminus Q_k}
\psi_n(x)  \biggr|\,dx\le \Sl_{k=0}^\iy\max\limits_{s\ge k}|\lm_s-\lm_{s+1}|.$$
\ep
This is an analog of the Sidon-Telyakovskii result for the one-dimensional
trigonometric series.

Obtained is also an analog of the estimate from above for the \Lc s from
{\bf 4.6} when $p=1.$
\flushpar
\proclaim{8.4} If an integer valued sequence $\{\nu_k\}_1^\iy$ is increasing
and convex and $\ln\nu_n=O(n^{1/2m}),$ then for every $f\in C[0,1]^m$
we have $$\Ll_{p\to\iy}{1\over p}\Sl_{k=1}^p|f(x)-S_{\nu_k}(f;x)|=0,$$
where $S_\nu$ are cubic partial sums of the Fourier-Walsh series.\ep

Let us return to one-dimensional series.
\flushpar
\proclaim{8.5} For every $\a>0$ and every $f\in L_q,$ where $q\in[1,+\iy],$
$$\g_1(\a)||f-\s_n(f)||_q\le||f-\s_n^\a(f)||_q\le\g_2(\a)||f-\s_n(f)||_q.$$
\ep
\proclaim{8.6} In $C[0,1]$ we have for every $n\in\BZ_+,$ $N\in(2^n,2^{n+1}],$
and $\a>0$
$$\g_1(\a)[\Om_n(f)+\om_{n+1}(f)]\le||f-\s_N^\a(f)||\le
\g_2(\a)[\Om_n(f)+\om_{n}(f)],$$ where
$$\Om_n(f)=\sup\limits_{k\ge n}\biggl\Vert{1\over 2^{k+1}}\Sl_{\nu=0}^k
2^{\nu-1}[f(\cdot)-f(\cdot\dot+{1\over 2^{n+1}})]\biggr\Vert,$$ and
$$\om_n(f)=\sup\limits_{0<t\le{1\over 2^n}}||f(\cdot\dot+t)-f(\cdot)||.$$\ep
For $\a=1$ somewhat weaker result was obtained earlier by M. F. Timan
and K. Tukhliev in 1982.

Now let us speak about the strong summability.
\flushpar
\proclaim{8.7} In order that for every $f\in L_\iy[0,1]$ the following estimate
$$\sup\limits_N{1\over N}\Sl_{k=1}^N|S_{\nu_k}(f;0)|\le\g||f||_\iy$$
holds, it is necessary, and in the case when $\{\nu_k\}$ is convex,
also sufficient that the following condition holds:
$$\sup\limits_N\Sl_{n=0}^n{1\over 2^{s+1}}\biggl({1\over N^2}\Sl_{k=1}^N
[(1-\mu)\nu+\mu(2^s-\nu)]^2\biggr)^{1/2}<\iy,$$
where $2^n\le\nu_N<2^{n+1}$ and $\nu_k=l2^{s+1}+\mu2^s+\nu,$ with
$l,$ $\mu,$ and $\nu$ integer defined by inequalities $\mu=\mu(k,s)\in
[1,2^s),$ $\nu=\nu(k,s)\in[0,1].$\ep
\flushpar
Unlike  the trigonometric system (see {\bf 3.6}a)), also the arithmetic nature
of $\{\nu_k\}_1^\iy$ is important here. For instance, the strong summability
takes place when $2^k\le\nu_k\le2^k+\sqrt{k}.$

A similar result is obtained also for sums like
$${1\over N}\Sl_{k=1}^N|S_{\nu_k}(f;0)|^{\th}$$ with $\th>1.$

Really speaking, Theorems {\bf 8.3} and {\bf 8.7} are obtained for
general orthonormal multiplicative systems. For the definition of the
Vilenkin multiplicative systems, see [M9], \S1.5.
\bigskip
\flushpar
{\bf Bibliographical remarks.}
\flushpar
For {\bf 8.1}, see [G6]; for {\bf 8.2}-{\bf 8.7}, see [G1-5].

\head 9. Pointwise approximation of functions by polynomials.
Approximation by polynomials with integral coefficients. \endhead
Let us consider approximation of functions by algebraic polynomials
on an interval of real axis with regard to position of a point.
Known are here the results due to S. M. Nikolskii, 1946, A. F. Timan,
1951 (direct theorems), and V. K. Dzyadyk, 1958 (inverse theorems);
see e.g., [M26,6].

In connection with results of type {\bf 5.14} arose a question in
the periodic case on characterization of the Lipschitz class of integral
order as regarded to degree of approximation. G. Alexits and G. Sunouchi,
1970, have made an attempt to answer this question. They succeeded
only partially, namely they had to recede from the ends of the interval.

Let us denote by $W_0^r,$ where $r\in\BN,$ the set of
functions on $[-1,1]$ which satisfy $f^{(r-1)}\in\Lip1$ and the following
boundary conditions
$$\Sl_{p=1}^{r-1}a_{p,s}f^{(p)}(1)/p!=0\quad\text{and}\quad
\Sl_{p=1}^{r-1}(-1)^p a_{p,s}f^{(p)}(-1)/p!=0.$$
Here $$a_{p,s}=\Sl_{q=1}^p\Sl_{k=p}^{r-1}\Sl_{\nu=1}^k(-1)^{p+q+k+\nu}
\binom pm\binom{q/2}k\binom{2k}{k-\nu}4^{-k}\nu^{2s}.$$
No conditions are for $r=1$ and for $r\ge2$ we have $s=[(r+1)/2],...,r-1.$
\proclaim{9.1} We have $f\in W_0^r$ if and only if
$$|f(x)-R_n(f;x)|=O(\sqrt{1-x^2}/n+o(1/n))^r$$ as $n\to\iy$ uniformly
with respect to $x\in[-1,1],$ where
$$\align R_n(f;x)&=a_0/2+\Sl_{k=1}^n(1-k^{2[(r+1)/2]}/n^{2[(r+1)/2]})a_k
C_k(x)\\&+{1+(-1)^{r+1}\over 2}\sqrt{1-x^2}\Sl_{k=1}^n(1-k/n)k^r/n^r
a_k S_{k-1}(x),\endalign$$
where $C_k(x)=\cos k\arccos x$ is the Chebyshev polynomial, $a_k$ are
the Fourier-Chebyshev coefficients, $S_{k-1}=C_k'/k$ is the Chebyshev
polynomial of second order.  \ep
For another approach to this problem, see [M6], p.268.

Let us go on to approximation by polynomials with integral coefficients.
Necessary and sufficient conditions of uniform approximation by such
polynomials are found by Fekete, 1954, Hewitt-Zukerman, 1959, and those
in complex domain by S. Ya. Alper, 1964 (see the survey [T5]). A
peculiar kind of intertwinning of analytic questions and arithmetic
ones appears here. These conditions are such that for functions defined
on sets of the complex plane we have the transfinite diameter to be less
than $1$ and the function, in addition, coincides with some polynomial
of integral coefficients. The problems of degree of approximation of
functions in various classes were studied as well. For instance,
theorems on the rate of convergence to functions of finite smoothness
as well as analytic on the interval $[0,1]$ are obtained by A. O. Gelfond,
1955 (see [S5]).
\proclaim{9.2} {\rm a)} Let $f$ be analytic inside and continuous in the
closed disk $K_\rho=\{z\in\BC:|z|\le\rho<1\}.$ If, in addition, the
numbers $f^{(\nu)}(0)/\nu!$ are integers for all $\nu=0,1,2,...$ (this
is also necessary), then for every $r\in\BN$ there exists a constant $\g$
independent of $n\in\BN$ such that
$$E_n^e(f;K_\rho)=\inf\limits_{Q_n}||f-Q_n||_{C(K_\rho)}\le\g\om_r(f;1/n).$$
Here $Q_n(z)=\Sl_{k=0}^n c_kz^k$ with $\Re c_k, \Im c_k\in\BZ$ for $0\le
k\le n,$ and $\om_r$ is the boundary modulus of smoothness (see before
{\bf 6.10}).

{\rm b)} Let $X(z)$ be a polynomial with integral coefficients and the
leading coefficient $1,$ and $K_\rho=\{z\in\BC:|X(z)|\le\rho<1\}.$
Suppose, further, that $f$ is analytic inside the lemniscate $|X(z)|=1.$
In order that $$\varlimsup\limits_{n\to\iy}(E_n^e(f;K_\rho))^{1/n}<1$$
it is necessary and sufficient that the polynomial $q_r$ defined by
conditions $q_r^{(s)}(z_\nu)=f^{(s)}(z_\nu)$ for $0\le s\le r,$ where
$\{z_\nu\}_1^k$ are all the different zeros of $X(z),$ has integral
coefficients for any $r.$  \ep

The following theorem is a general theorem on approximation of functions
of finite smoothness on arbitrary interval of $\BR.$
\proclaim{9.3} Let both $\{x_\nu\}_1^k,$ the set of all integral
algebraic numbers, and the set of their algebraically conjugate numbers
situate on $[a,b]$ with $b-a<4.$ If $f\in C^r[a,b]$ for $r\in\BZ_+$
and Hermitian polynomial $q,$ defined by conditions $q^{(s)}(x_\nu)
=f^{(s)}(x_\nu)$ for $1\le\nu\le k$ and $0\le s\le r,$ has integral
coefficients (this is also necessary), then a sequence of polynomials
$Q_n$ of degree not higher than $n$ with integral coefficients can be
picked out so that they satisfy on $[a,b]$ the following inequalities
for $0\le\nu\le r:$
$$|f^{(\nu)}(x)-Q_n^{(\nu)}(x)|=O[({\sqrt{(x-a)(b-x)}\over n}+1/n^2)^{r-\nu}
\om(f^{(r)};{\sqrt{(x-a)(b-x)}\over n}+1/n^2)].$$ 
Besides that, $$|Q_n^{(r+1)}(x)|=O[({\sqrt{(x-a)(b-x)}\over n}+1/n^2)^{-1}
\om(f^{(r)};{\sqrt{(x-a)(b-x)}\over n}+1/n^2)].$$
\ep
Although it is impossible to put $\om_k(f)$ with $k\ge2$ in this
theorem, the same theorem is true for functions satisfying Zygmund's
condition $\om_2(f^{(r)};h)=O(h).$ A similar theorem holds with $r=0$
for functions with given majorants of partial moduli of smoothness
on the parallelepiped with the edges parallel to the coordinate axes
of length less than $4,$ see [T5].
\proclaim{9.4} Let $f$ be analytic inside the square $[0,1]\times[0,1]$
and $f\in C^r$ on $[0,1]\times[0,1].$ Suppose, further, that Hermitian
interpolation polynomial $q(z),$ defined by conditions $q^{(s)}(z_k)
=f^{(s)}(z_k)$ for $0\le s\le r$ and $1\le k\le4,$ $z_1=0,$ $z_2=1,$
$z_3=i,$ $z_4=1+i,$ has integral coefficients. Then for each $n$ there
exists $Q_n$ such that for any $z\in\G,$ where $\G$ is the boundary
of the square, we have
$$|f^{(\nu)}(z)-Q_n^{(\nu)}(z)|\le\g\rho_{1+1/n}^{r-\nu}(z)\om
(f^{(r)};\rho_{1+1/n}(z)),$$ where $\rho_{1+\ve}(z)$ is the distance
from $z\in\G$ to the level curve $\G_{1+\ve}.$  \ep
This is the first result of Dzyadyk's type (see [M6], Ch.9) on
approximation by polynomials with integral coefficients on a set in $\BC$
with nonsmooth boundary.

In his lecture on the 1st All-Union Congress of Mathematicians,
S. N. Bernstein posed, among other questions, a question on best
approximation of arbitraru number $\lm$ by polynomials $Q_n$ on an
interval situated in $(0,1).$ A certain upper bound for the interval
$[\dl,1-\dl],$ where $\dl\in(0,1/2),$ was given not long after by
R. O. Kuzmin and L. V. Kantorovich (see [S5]). Let us give one precise
result.
\proclaim{9.5} For any $\lm\in(0,1)$ we have
$$E_n^e(\lm;\dl,1-\dl)=\inf\limits_{Q_n}\max\limits_{[\dl,1-\dl]}
|\lm-Q_n(x)|\le 2n\rho^n,$$ where $\rho=\max\{1/2,(1-2\dl)/(1+2\sqrt{\dl
(1-\dl)})\}.$ Here $\rho=1/2$ for $\dl\in[1/10,1/2),$ and we have
change-over for $\dl=1/10,$ and the factor $2n$ can be replaced by
a bounded one if $\dl\ne1/10;$ it is impossible to make $\rho$ smaller
if $\lm$ is not dyadic.  \ep

The case of dyadic $\lm$ is also investigated.
\proclaim{9.6} {\rm a)} If inside $[a,b]$ with $b-a<4$ there exists
at least one integer point, then for every $\lm\in(0,1)$ and for
every $p\in[1,+\iy)$ we have
$$E_n^e(\lm;a,b)_p=\inf\limits_{Q_n}(\Il_a^b|\lm-Q_n(x)|^p\,dx)^{1/p}
\asymp n^{-1/p}.$$

{\rm b)} For any $b\in(0,1]$ and for any $\lm\in(0,1)$ we have
$$E_n^e(\lm;0,b)_p\asymp n^{-2/p}.$$   \ep
Somewhat weaker result was obtained earlier by E. Aparisio ($p=2$)
and A. O. Gelfond (for any $p\ge1$) (see [S5]). A result similar to
{\bf 9.6} is proved also for the ball in $\BR^m.$

Strengthened is also one result due to A. O. Gelfond, 1966, on
approximation of functions on $[\a,\be]\sbt(0,1)$ by polynomials
with coefficients from
a given set (with restrictions to the growth of absolute values; see [T11]).
\proclaim{9.7} Let $E$ be a compactum in $\BC$ with connected complement,
the origin be a point of $\p E$ and not being a limit point for inner
points of $E$ if are those. If $\{w(k)\}_0^\iy$ is a positive sequence
satisfying $\lim(w(k))^{1/k}=\iy$ as $k\to\iy,$ then for every $f\in C(E)$
and analytic inside such that $f(0)=0$ and for every $\ve>0$ there exists
$\{c_k\}_0^n$ we have $$\max\limits_{z\in E}|f(z)-\Sl_{k=0}^n c_k z^k|<\ve
\quad\text{and}\quad |c_k|\le w(k)\quad\text{for all}\quad k\ge0.$$
If a sequence $\{(w(k))^{1/k}\}$ is bounded, then this statement is false.\ep
For sets without inner points this result was proved earlier by S. Ya.
Khavinson, 1969, and more general result than {\bf 9.7} was obtained later by
V. A. Martirosyan.

The problem of best approximation in $L_q(\BT)$ of periodic functions
of class $W_p^r(\BT)$ by trigonometric polynomials of degree not higher
than $n$ was solved for $q=\iy$ and $p=\iy,$ in the space $C(\BT),$
and for $q=1$ and $p\in[1,+\iy]$ (J. Favard, S. M. Nikolskii, L. V.
Taikov; see [M13]). In the case of approximation by polynomials on
interval with regard to position of a point no such results are known.
\proclaim{9.8} For every $r\in\BN,$ for every $f\in W_\iy^r[-1,1],$
and for every $n\ge r-1$ there exists a polynomial $p_n$ of degree
not higher than $n$ satisfying the inequality
$$|f(x)-p_n(x)|\le K_r(\sqrt{1-x^2}/(n+1))^r+\g(r)(\sqrt{1-x^2})^{r-1}/
(n+1)^{r+1},$$ where $K_r=(4/\pi)\Sl_{k=0}^\iy(-1)^{k(r+1)}(2k+1)^{-r-1}$
and necessarily $\g(r)\ge Ce^r.$  \ep
This result was known earlier only for $r=1$ due to V. N. Temlyakov, 1981.

The problems of approximation by polynomials in integral metric with
power weight were investigated by many mathematicians (see [S12]).
\proclaim{9.9} {\rm a)} For every $p\in[1,+\iy),$ for every $r\in\BN,$ for
every $f\in W_p^r[-1,1],$ and for every $n\ge r-1$ there exists a polynomial
$p_n$ of degree not higher than $n$ satisfying the inequality
$$\Il_{-1}^1|f(x)-p_n(x)|(1-x^2)^{r/2}\,dx\le||\vp_r||_{p'}(n+1)^{-r}
+o((n+1)^{-r}),$$ where $1/p+1/p'=1$ and $\vp_r(t)$ is the $r$-th
$2\pi$-periodic integral of $\sign\sin t.$

{\rm b)} For every $f\in W_\iy^r[-1,1],$
and for every $n\ge r-1$ there exists a polynomial $p_n$ of degree
not higher than $n$ satisfying the inequality
$$\Il_{-1}^1|f(x)-p_n(x)|(1-x^2)^{r/2}\,dx\le 4K_{r+1}(n+1)^{-r}
+\g(r)(n+1)^{-r-1}$$ and this result id asymptotically sharp
on the class. \ep

The proof of {\bf 9.8} differs from that of {\bf 9.9} mainly because
the operator $p_n=p_n(f)$ is nonlinear, while in {\bf 9.9} it is linear.

Direct theorems on approximation by polynomials are applied, for example,
to the investigation of the question of convergence of  Fourier-Jacobi
series (see [M20]).
\proclaim{9.10} {\rm a)} If there exist $\g\in(0,1],$ $\dl<2\g,$ and
$\lm<2\g$ such that for any $x_1, x_2\in[-1,1]$
$$|f(x_1)-f(x_2)|\le|x_1-x_2|^\g(\sqrt{1-x_1}+\sqrt{1-x_2})^{-\dl}
(\sqrt{1+x_1}+\sqrt{1+x_2})^{-\lm},$$ then there exists a sequence
of polynomials $p_n$ such that for each $x\in[-1,1]$  $$|f(x)-p_n(x)|
=O(n^{-\g}(\sqrt{1-x}+1/n)^{\g-\dl}(\sqrt{1+x}+1/n)^{\g-\lm}).$$

{\rm b)} Let $f(\pm1)=0$ and $f(x)(1-x)^{\dl/2}(1+x)^{\lm/2}\in\Lip\g$
for some $\g\in(0,1],$ $\dl<2\g,$ and $\lm<2\g,$ then its Fourier-Jacobi
series orthonormal in $L_2$ with the weight $(1-x)^\a(1-x)^\be$
converges uniformly on $[0,1]$ for $\a-\be<1+\lm-2\g.$

{\rm c)} If the Fourier-Jacobi series of $f$ converges at the endpoints
of $[-1,1]$ and the Dini-Lipschitz condition $\om(f;h)\ln(1/h)\to0$ as
$h\to0$ holds, then this series converges uniformly on $[-1,1].$  \ep
\flushpar
{\bf Bibliographical remarks.}
\flushpar
For {\bf 9.1}, see [T7]; for {\bf 9.2}, see [T15];
for {\bf 9.2, 9.5, 9.6}, see [T1]; see also [T5].
{\bf 9.4} is due to Vit. Volchkov (the paper is accepted in Mat. Zametki);
for {\bf 9.7}, see [T11]; for {\bf 9.8, 9.9} ($p=1$), see [T34];
for {\bf 9.10}, see [Bl1,2].
\flushpar
See also [S3], and for integral coefficients, see [M8], where some
of the results mentioned are given with proofs.
\head 10. Approximation by trigonometric polynomials with a given number
of harmonics. \endhead
Let $E$ be a Banach space, $W$ be a compactum in $E$ symmetric
with respect to the origin, $L$ be a subspace of $E.$ The (Kolmogorov)
width of $W$ of order $n$ is defined as follows:
$$d_n(W)_E=\inf\limits_{L:\dim L=n}\sup\limits_{f\in W}\inf\limits_{g\in L}
||f-g||.$$ The space of dimension $n$ on which the least lower bound
is achieved, if one exists, is called extremal (the best for $W$).
In the case $E=C(\BT)$ and $W=W^r,$ that is the cylinder with compact base,
extremal is the subspace of trigonometric polynomials with spectrum
on the interval and the width is calculated (V. M. Tikhomirov, 1959).

K. I. Babenko, 1960, posed and investigated the problem of the choice
of approximate polynomials for approximation in $L_p$ of function classes
on $\BT^m$ which are defined by a given differential operator.

R. S. Ismagilov (see [S7]) has found out that for usual classes of
periodic functions even for $m=1$ polynomials with fixed spectrum do not
form even asymptotically best subspace and introduced the following
trigonometric widths:
$$d_n^T(W)_q=\inf\limits_{\th_n}\sup\limits_{f\in W}\inf\limits_{T(\th_n)}
||f(\cdot)-T(\th_n;(\cdot))||_q,$$   where
$$T(\th_n;x)=\Sl_{k=1}^n c_ke^{i(p_k,x)},$$ and $\th_n=\{p_1,...,p_n\},$
that is minimization first with respect to $\{c_k\}_1^n$ and then
with respect to the spectrum $\th_n$ of $n$ harmonics.

The notion of trigonometric width makes a sense for separate functions
as well and is as follows:
$$e_n(f)_q=\inf\limits_{\th_n}\inf\limits_{T(\th_n)}||f-T(\th_n)||_q.$$
For example, for $f_0(x)=|x|$ with $x\in\BT$ we have $e_n(f_0)_\iy\asymp
n^{-3/2}$ (V. E. Majorov, 1986). This is essentially nonlinear problem
of adaptive approximation.

Let us start with the problem  on the \Lc.

In 1981, S. V. Konyagin and O.C.McGehee-L.Pigno-B.Smith independently
proved Littlewood's conjecture: for any collection $\th_n=\{p_k\}_1^n$ of
integers we have $$\Sl_{k=1}^ne^{ip_k(\cdot)}||_1=\Il_{-\pi}^\pi
|\Sl_{k=1}^ne^{ip_kx}|\,dx\ge\g\ln n.$$ For this, see [M12]. Using
S. A. Vinogradov's version of the proof due to the three authors (oral
communication, 1983) we obtain the following general fact.
\proclaim{10.1} For any $\{p_k\}_1^\iy$ such that $p_k\in\BN$ and
$p_{k+1}>p_k$ and for any $\{c_k\}_1^\iy$
$$\Il_{-\pi}^\pi|\Sl_{k=1}^ne^{ip_kx}|\,dx\ge\g\Sl_{s=1}^\iy(\Sl_
{2^{s-1}\le\nu<2^s}\nu^{-1}|c_\nu|^2)^{1/2}.$$  \ep

V. M. Tikhomirov suggested to find an order, as $n\to\iy,$ of the value
$$L_n(r,q)=\inf\limits_{\{p_k\}_1^n}||(\Sl_{k=1}^n e^{ip_k(\cdot)})^
{(r)}||_q,$$ where the $r$-th derivative in the Weyl sense for $r\ge0$
is written inside the norm sign. The afore-mentioned results immediately
give $L_n(0,1)\asymp\ln n.$
\proclaim{10.2} The following relations are true with $1/q+1/q'=1$
$$L_n(r,q)\asymp\cases n^r\ln n,& q=1, r\ge0\\ n^{r+1/q'},& 1<q\le2, r\ge0\\
&\text{or}\quad 2<q\le\iy, r>1/q\\ n^{(rq+1)/2},& 2<q<\iy, 0\le r<1/q\\
n(\ln n)^{-1/q'}, & 2<q<\iy, r=1/q. \endcases$$ \ep
The first two relations for  any $r>0$ and the third one for even $q$
are established by V. E. Majorov. The other cases are due to E. Belinskii,
1988. The third relation in the general case was independently proved
by S. V. Konyagin, 1988. There exist alsso a generalization to the
multiple case due to E. Belinskii and E. M. Galeev, 1991.

Inequalities of different metrics for polynomials with floating
spectrum look as follows.
\proclaim{10.3} For all $1\le p<q\le\iy$ we have
$$\sup\limits_{\th_n}{||T(\th_n)||_q\over||T(\th_n)||_p}\asymp\cases 
n^{(1/p-1/q)p/2}, & 2<p<\iy\\ n^{1/p-1/q},& 1\le p\le2.\endcases$$\ep
The upper bounds are due partially to R. J. Nessel and G. Wilmes, 1978,
and the rest to V. A. Rodin, 1983.

Let us go on to problems of approximation of basic classes.

We say that $f\in W_p^{r,\be}$ if
$$\Sl_{k\in\BZ}|k|^r e^{(i\be\pi\sign k)/2}c_k(f)e_k\sim f^{r,\be}\in
L_p(\BT)$$ and $||f^{r,\be}||_p\le 1.$ For definition of Nikolskii's classes
$H_p^r,$ see [M16].
\proclaim{10.4} For $1\le p\le2$ and $2<q<\iy$ the following relations hold:
$$e_n(H_p^r)_q=\sup\limits_{f\in H_p^r} e_n(f)_q\asymp\cases
n^{-(r-1/p+1/2)}, & r>1/p\\ n^{-(r-1/p+1/q)q/2}, & 1/p>r>1/p-1/q\\
n^{-1/2}\ln n, & r=1/p. \endcases$$   \ep
\proclaim{10.5} The following relations hold:
$$e_n(W_p^{r,\be})_q=\sup\limits_{f\in W_p^{r,\be}} e_n(f)_q\asymp\cases
n^{-(r-1/p+1/2)}, & r>1/p\\ n^{-(r-1/p+1/q)q/2}, & 1/p>r>1/p-1/q\\
n^{-1/2}\ln^{1-1/p}n, & r=1/p,\ p>1\\n^{-1/2}\ln n, & r=p=1.\endcases$$ \ep

For the critical value $r=1/p$ the difference has come to light in
approximate properties of $H_p^r$ and $W_p^{r,\be}$ (cf. {\bf 10.4}
and {\bf 10.5}). This is for the first time in the one-dimensional case.

The same problem is solved also for the scale of Besov's spaces 
(for definitions, see [M16]).
\proclaim{10.6} The following relations hold:
$$e_n(B_{ps}^r)_q=\sup\limits_{f\in B_{ps}^r} e_n(f)_q\asymp\cases
n^{-(r-1/p+1/2)}, & r>1/p\\ n^{-(r-1/p+1/q)q/2}, & 1/p>r>1/p-1/q\\
n^{-1/2}\ln^{1-1/s}n,& r=1/p,\ 1\le s<\iy\\n^{-1/2},& r=1/p,\ 0<s<1.
\endcases$$ \ep

Hence for $r\ne1/p$ the degree of approximation is the same for each
of the three classes and coincides with B. S. Kashin's results, 1977-1981,
on Kolmogorov's widths for $W_p^r.$

The proofs of the given theorems as well as the following ones are
based on the next statement.
\proclaim{10.7} {\rm a)} For any polynomial $T(\th_N)$ with spectrum
$\th_N,$ for every $n\le N,$ and for every $q\in(2,+\iy)$ there exist
$\th_n\sbt\th_N$ and $T(\th_n)$ such that
$$||T(\th_N)-T(\th_n)||_q\le\g(q)(N/n)^{1/2}||T(\th_N)||_2.$$

{\rm b)} For any polynomial $T_N$ with spectrum in $[0,N],$ for
every $n\in(1,N),$ and for every $p\in[2,+\iy)$ there exist $\th_n\sbt
[0,2N]$ and $T(\th_n)$ such that
$$||T_N-T(\th_n)||_\iy\le\g(p)((N/n)\ln(N/n))^{1/p}||T_N||_p.$$  \ep

It is obvious that $e_n(f)_2=(\Sl_{k=n+1}^\iy(c_k^*)^22\pi)^{1/2},$
where $\{c_k^*\}_1^\iy$ is decreasing rearrangement of $\{|c_k(f)|\},$
the absolute values of \Fc\ of $f.$
\proclaim{10.8} {\rm a)} Let $q\in(2,+\iy)$ and $p\in(0,1].$ Then
$$\{c_k(f)\}_{-\iy}^\iy\in l_p\quad\text{if and only if}\quad
\Sl_{k=1}^\iy(2^{k(1/p-1/2)}e_{2^k}(f)_q)^p<\iy.$$

{\rm b)} If for some $\ve>0$ we have $\sum(|c_k(f)|(\ln k)^{1+\ve})^p<\iy,$
then $$\align \Sl_{k=1}^\iy(2^{k(1/p-1/2)}e_{2^k}(f)_\iy)^p&<\iy\qquad
\text{for}\quad p\in(0,2/3),\\ \quad\\
\Sl_{k=1}^\iy(2^{2k(1/p-1)}e_{2^k}(f)_\iy)^p&<\iy\qquad
\text{for}\quad p\in(2/3,1),\\ \quad\\
\Sl_{k=1}^\iy(2^k e_{2^k}(f)_\iy)^{2/3}k^{-1}&<\iy\qquad
\text{for}\quad p=2/3. \endalign$$  \ep

Let us go on to the multiple case. Let $r=(r_1,...,r_m)\in\BR_+^m$ and
$0<r_1=r_2=...=r_\nu<r_{\nu+1}\le...\le r_m.$ We  denote by $H_p^r$ the
class of functions satisfying the conditions $\Il_{-\pi}^\pi f(x)\,dx_j=0$
for $1\le j\le m$ and $$||\Dl_h^l f||_p\le\prod\limits_{j=1}^m|h_j|^{r_j},$$
where $\Dl_h^1$ is the mixed difference with step $h_j$ in direction
$Ox_j,$ $1\le j\le m,$ and $\Dl_h^l=(\Dl_h^1)^l$ with $l\in\BN$ and $l>r_m.$
The study of problems of approximation of this class was started by
K. I. Babenko and continued by S. A. Telyakovskii, Ya. S. Bugrov and
others; see [M25].
\proclaim{10.9} Let $r_1>1/2$ and $p\in[2,+\iy].$ Then
$$\g_1(n^{-1}\ln^{\nu-1}n)^{r_1}\ln^{(\nu-1)/2}n\le e_n(H_p^r)_\iy\le
\g_2(n^{-1}\ln^{\nu-1}n)^{r_1}\ln^{\nu/2}n.$$   \ep

A similar result is obtained also for any dimension for
$e_n(W_p^{r,\be})_\iy.$

And at last let us give results in the final form for trigonometric
widths in the multiple case. These results partially rehabilitate the 
trigonometric system as the best in the sense considered.
\proclaim{10.10} The following relations valid:
$$d_n^T(H_p^r)_q\asymp\cases (n^{-1}\ln^{\nu-1}n)^{r_1-1/p+1/2}
\ln^{(\nu-1)/2}n,& r_1>1,\ 1\le p<2<q<p'\\
(n^{-q/2}\ln^{(\nu-1)(q-1)}n)^{r_1-1+1/q}\ln^{(\nu-1)/q}n,& p=1,\ 1-1/q<r_1<1, 
\\ &2<q<\iy\\ n^{-1/2}\ln^\nu n,& p=r_1=1, 2<q<\iy.\endcases$$ \ep
The case $1<p\le q\le2$ was investigated earlier by V. N. Temlyakov, 1982.
\proclaim{10.11} The following relations valid:
$$d_n^T(W_p^{r,\be})_q\asymp\cases (n^{-1}\ln^{\nu-1}n)^{r_1-1/p+1/2}
,& r_1>1,\ 1< p\le2\le q<p'\\ (n^{-1}\ln^{\nu-1}n)^{r_1-1/2}
\ln^{(\nu-1)/2}n,& p=1<r_1,\ 2<q<\iy \\
(n^{-q/2}\ln^{(\nu-1)(q-1)}n)^{r_1-1+1/q}\ln^{(\nu-1)/q}n,& p=1,\
1-1/q<r_1<1\\ n^{-1/2}\ln^\nu n,& p=r_1=1.\endcases$$ \ep

Some problems are still open. But let us note that this field is
ahead, in some respects, of results in Kolmogorov's widths (for
example, for $r_1=1$) having voluminous literature (see [M13] and [M27]).
\bigskip
\flushpar
{\bf Bibliographical remarks.}
\flushpar
For {\bf 10.2-10.3}, see [Be16]; for {\bf 10.4-10.6, 10.7}a), see [Be12,13];
for {\bf 10.7}b), see [Be10,20]; for {\bf 10.8}, see [Be18];
for {\bf 10.9}, see [Be15]; for {\bf 10.10, 10.11}, see [Be11].
Pay attention also to [Be20].
\flushpar
Starting in 1984, methods of random choice are used here (E. Belinskii,
Yu. Makovoz); see also [M11,12]. For interesting applications to the
problem of approximation by degenerate functions, see [M25].

\head 11. The Pompeiu problem. \endhead
Let $A$ be an open bounded set in $\BR^m$ for $m\ge 2,$
and $M(m)$ be the group of Euclidean motions of $\BR^m.$
The question arises (Pompeiu, 1929) if there exists a non-trivial locally
integrable function $f:\BR^m\to\BC$ such that 
$$\Il_{\s A}f(x)\,dx=0\qquad\text{for every}\quad \s\in M(m).\tag*$$
Are of certain interest connections of this problem and various of
its versions with mean-value theorems for differential equations (see
[S10]), with the classical Morera theorem in complex analysis (see 
{\bf 11.5} below), with uniqueness theorems (see {\bf 11.7}), etc.

The existence of non-trivial function with (*) condition immediately
yields an estimate of density of packing of an arbitrary compactum
$K\sbt\BR^m$ by sets $\s A\in K$ (B. D. Kotlyar, 1984).

Taken in integral (*) a ball as $A$ and $f(x)=e^{i(x,y)}$ we get
the Bessel (radial) function $J_{m/2}$ (see [M22, Ch.IV, \S4]), and hence
condition (*) is satisfied on the countable set of balls (and their
shifts) radia of which are defined by zeros of $J_{m/2}.$ There
already exist a number of examples of sets $A$ for which (*) implies
that $f$ is trivial, that is vanishes a.e. In particular, so is when
when the boundary of $A$ is non-analytic (Williams, 1976; see 
also [S16]). Various criteria are also obtained.
\flushpar
\proclaim{11.1} {\rm a)} If condition (*) is satisfied for some $A$ and some
non-trivial function of power growth, then the function $f(x)=e^{i(x,y)}$
satisfies condition (*) for some $y\in\BR^m.$

{\rm b)} Let for $m=2$ a set $A$ be star-like woth respect to the origin,
$\r=\r(\vp),$ with $0\le\vp\le2\pi,$ be the polar equation of $\p A,$
where $\r$ is piece-wise continuous. In order that (*) implies that
$f$ is trivial, it is necessary and sufficient that there exists $r>0$
such that $$\Il_0^{2\pi}e^{-ik\vp}\,d\vp\Il_0^{\r(\vp)}tJ_k(tr)\,dt=0$$
for each $k\in\BZ,$ where $J_k$ is the Bessel function of order $k.$\ep

A similar criterion is obtained also for every $m\ge3.$
\flushpar
\proclaim{11.2} Let $A$ be a set with connected complement. In order that
there exists a non-trivial function satisfying (*), it is necessary
and sufficient the indicator of $A$ to be a limit in $L_1(\BR^m)$ of
a sequence of linear combinations of indicators of balls with radia
proportional to positive zeros of $J_{m/2}.$ The coefficient of
proportionality depends only on $A.$\ep

The problem becomes much more difficult if a function is defined not
on the whole space $\BR^m,$ where the \Ft\ methods are applicable,
but on a bounded domain.
\flushpar
\proclaim{11.3} Let $B_r$ be the open ball in $\BR^m$ of radius $r$ 
with center at the origin and $f$ be locally integrable in $B_r.$

{\rm a)} If $A$ is the unit cube in $\BR^m$ and (*) holds for every $\s\in
M(m)$ such that $\s A\sbt B_r,$ then for $r\ge{1\over 2}\sqrt{m+3}$
we have $f=0$ and for $r<{1\over 2}\sqrt{m+3}$ there exist nontrivial
functions satisfying the given condition.

{\rm b)} If $A$ is the unit semiball in $\BR^m$ 
and (*) holds for every $\s\in
M(m)$ such that $\s A\sbt B_r,$ then for $r\ge{\sqrt{5}\over 2}$
we have $f=0$ and for $r<{\sqrt{5}\over 2}$ there exist nontrivial
functions satisfying the given condition.
 
{\rm c)} If $A$ is an ellipsoid in $\BR^m$ and $\lm$ and $\mu$ are its largest
and its smallest halfaxis, respectively; $R(x,y)=2x$ for $y\le\sqrt{2}x$
and $R(x,y)=y^2(y^2-x^2)^{-1/2}$ for $y>\sqrt{2}$; and (*) holds for every 
$\s\in M(m)$ such that $\s A\sbt B_r,$ then for $r\ge R(\lm,\mu)$
we have $f=0$ and for $r<R(\lm,\mu)$ there exist nontrivial
functions satisfying the given condition. \ep

Found is complete description of a class of functions with non-zero
integrals over all the balls of a given radius (in terms of coefficients
of the \Fs\ of the function $f(\r\eta)$  for $\r>0$ with respect to the
spherical harmonic expansion of $\eta$ on the unit sphere).
\flushpar
\proclaim{11.4} Let the integrals of $f$ over all the balls in $B_r$
of radia $r_1$ and $r_2$ vanish. If $r_1+r_2<r$ and $r_1/r_2$ does not
equal to ratio of two different zeros of $J_{m/2},$ then $f=0.$
If $r_1+r_2>r$ or $r_1/r_2$ is equal to ratio of two different zeros
of $J_{m/2},$ then there exists a nontrivial function satisfying
the condition indicated above.\ep
\flushpar
The case $r_1+r_2=r$ is also completely investigated.

The Morera theorem can be strengthened as follows.
\flushpar
\proclaim{11.5} Let $D$ be the unit disc in $\BC$ and $f\in\BC(D).$

{\rm a)} If the integrals over all the circles tangent to $\p D$ vanish,
then $f$ is analytic in $D.$

{\rm b)} If the integral of $f$ over the boundary of every square in $D$
with a fixed side $d$ vanish, the $f$ is analytic for $d\le2/\sqrt{5}$
and not necessarily for $d>2/\sqrt{5}.$

{\rm c)} If the integral of $f$ over the boundary of every equisided triangular
in $D$ with a fixed side $d$ vanish, the $f$ is analytic for $d\le2/\sqrt{3}$
and not necessarily for $d>2/\sqrt{3}.$\ep
 
Similar results are obtained also for the multidimensional case in $\BC^m.$
Such conditions are given for harmonic functions as well.

The following theorem strengthens a criterion of analyticity of
functions due to V. K. Dzyadyk, 1960.
\flushpar
\proclaim{11.6} Let real functions $u$ and $v$ be in $C^2(D),$ in the
sense of real analysis, and $r_1+r_2<1$ and $r_1/r_2$ is not equal to ratio
of two zeros of $J_1.$ In order the function $f=u+iv$ to be analytic
in $D,$ it is necessary and sufficient that areas of surfaces of the graphs
of functions $u,$ $v,$ and $\sqrt{u^2+v^2}$ placed over each disc from
$D$ of radius $r_1$ or $r_2$ to be equal.\ep

Let us give an example of the uniqueness theorem for lacunary trigonometric 
system.
\flushpar
\proclaim{11.7} Let $D$ be an arbitrary domain in $\BR^m$ and $\lm_k\in\BR^m$
for all $k\in\BN$ be such that $|\lm_{k+1}|-|\lm_k|\to\iy$ as $k\to\iy.$
Assume further that a sequence of polynomials $\Sl_{k=1}^n c_{k,n}e^{i
(\lm_k,x)}$ converges weakly in $L_1(D)$ as $n\to\iy$ to a function
$f$ which is vanishing in some ball. Then $f=0$ a.e. in $D.$\ep
\flushpar
\proclaim{11.8} Let $f$ be locally integrable in $\BR^m$ and $\{\nu_k\}_1^\iy$
be the sequence of all the positive zeros of the function $J_{m/2}.$
In order the integrals of $f$ over all the balls from $\BR^m$ with
radia $\nu_1,$ $\nu_2,...$ to be zero, it is necessary and sufficient
that $f$ coincides a.e. with a solution of the Helmholz equation
$\Dl f+f=0.$\ep
\proclaim{11.9} Let $G$ be a bounded Jordan domain in $\BC$ with
rectifiable boundary $\p G$ and $u=u(z)$ be a function harmonic in
$\BC\sm\bar G$ and continuous in $\BC\sm G.$ Let further $u$ satisfy the 
following conditions: $u=0$ on $\p G,$ ${\p u\over\p n}=1$ a.e. on $\p G,$ and
$u(z)=o(|z|^2)$ as $z\to\iy.$ Then $G$ is the disk of center $z_0$
and radius $R$ and $u(z)=(1/R)\ln|(z-z_0)/R|.$ The condition at infinity
is sharp. \ep
\flushpar
{\bf Bibliografical remarks.}
\flushpar
For {\bf 11.1}a), see [Za2]; for {\bf 11.1}b), see [ZT2];
for {\bf 11.2}, {\bf 11.8}, see [VaV8]; for {\bf 11.3}, see [VaV1,7];
for {\bf 11.4}, see [VaV9,10]; for {\bf 11.5}, see [VaV5,11];
for {\bf 11.6}, see [VaV5]; for {\bf 11.7}, see [VaV4]. See also
[VaV3,6,12-14].

\head 12. Entire functions of exponential type and 
the Fourier integral.\endhead
The \Ft\ of a function with compact support is the entire function of
exponential type (EFET), see e.g., [M14].

Well-known are Wiener's and Cartright's theorems on sufficient conditions
for membership to $L_1$ and $L_\iy,$ respectively, of restriction of EFET
on $\BR$ according to the behavior of the function at equidistant points
of $\BR$ (see the same references). Generalizations of these theorems are
obtained. Let us give here only the asymptotic formula.
\proclaim{12.1} Let $f$ be an entire function of exponential type
not higher $p\pi$ with $p\in\BN,$
$$\Ll_{y\to\pm\iy}|f(iy)|e^{-p|y|\pi}=0\quad\text{and}\quad
\Sl_{k\in\BZ}|f^{(\nu)}(k)|<\iy$$ for $0\le\nu\le p-1.$ Then for every
$n\in\BZ$ and for every $x\in[n,n+1]$
$$f(x)=(-1)^n\pi^{-p}f^{(p)}(n)\sin^p(x\pi)/p!+\th\Sl_{k\in\BZ}
(k-n+1/2)^{-2}\Sl_{\nu=0}^{p-1}|f^{(\nu)}(k)|$$ with $|\th|\le\g(p).$ \ep

If $f$ is of type not higher than $\s,$ $\lim x^{-1}f(x)=0$ as $|x|\to\iy,$
and $\sup\limits_k|f(k\pi/\s)|<\iy,$ then $f$ is not obliged to be bounded
on $\BR$ while the sum $f(x)+f(x+\pi/\s)$ is bounded. This phenomenon
was found by S. N. Bernstein, 1949, and he called it interference.
There are general theorems on interference in $L_\iy$ (see [M2,13,23]).
Let us consider the interference in $L.$

Let $L^1_\pi$ be a set of entire functions of exponential type not higher
than $\pi$ such that $\lim f(x)=0$ as $|x|\to\iy$ and $\Sl_k|f(k)|<\iy.$
Then $$f(z)=(2\pi)^{-1}\Il_{-\pi}^\pi\hf(u)e^{-iuz}\,du\quad\text{with}
\quad\hf\in A(\BT).$$
Let us call a linear operator $\Lm$ taking $L_\pi^1$ into itself to be
interference if $$\Il_{-\iy}^\iy|(\Lm f)(x)|\,dx\le
C\Sl_{k=-\iy}^\iy|f(k)|.$$
\proclaim{12.2} In order that a linear operator $\Lm$ continuous with
respect to the uniform convergence on any compactum is interference and
commuting with the unit shift operator, it is necessary and sufficient
that it is representable in the form
$$(\Lm f)(z)=\Il_{-\pi}^\pi\lm(u)\hf(u)e^{-iuz}\,du\quad\text{with}
\quad\hat\lm\in L(\BR).$$ Besides that,
$$\Il_{-\iy}^\iy|(\Lm f)(x)|\,dx\le\Il_{-\iy}^\iy|\hat\lm(u)|\,du
\Sl_{k=-\iy}^\iy|f(k)|$$ and the constant $||\hat\lm||_1$ is exact on the
class $L_\pi^1.$  \ep

The following theorem is the analog of B. Ya. Levin's theorem on the $L_\iy$
space (see [M14]).
\proclaim{12.3} In order that there exists an entire function of
exponential type not higher than $\pi$ in $L(\BR)$ and taking the values
$c_k$ at the points $k\in\BZ$ so that $\sum|c_k|<\iy,$ it is necessary
and sufficient that $$\Sl_{n\in\BZ}|\Sl_{k\in\BZ}(-1)^{k+n} c_{k+n}
k/(k^2+1)|<\iy.$$  \ep

Functions defined on $\BR$ are approximated, following S. N. Bernstein,
by EFET as their type $\s\to\iy.$ In particular, if a function is
periodical, then the approximating EFET in $L_\iy$ can be considered
periodical, that is trigonometric polynomial.

The following theorem is an analog of {\bf 5.1}.
\proclaim{12.4} Let $p\in[1,+\iy]$ and $f\in L_p(\BR).$ Then for any
$r\in\BN$ and for any $\s>0$ there exists an entire function $g_\s$ of
exponential type not higher than $\s$ such that
$$||f-g_\s||_p\asymp\om_r(f;1/\s)_p$$ with constants in the two-sided
inequality depending only on $r.$   \ep

The following statements are very well known for $p\ge1$ (see [M26,16]).
The first one is on equivalence of integral and difference norms,
while the second one is of Bernstein type.
\proclaim{12.5} {\rm a)} For any entire function of exponential type
not higher than $\s$ and for any $p\in(0,1)$    $$\g_1(p)||f||_p
\le((\pi/2\s)\Sl_k|f(k\pi/(2\s))|^p)^{1/p}\le\g_2(p)||f||_p.$$

{\rm b)} For any entire function of exponential type not higher than
$\s$ and for any $p\in(0,1)$    $$||f'||_p\le\g(p)\s||f||_p.$$   \ep

Statement a) for trigonometric polynomials was obtained by V. V. Peller
(see [S13]), and b) for trigonometric polynomials was proved in 1975
independently by several authors, among them that with the exact
constant $\g(p)=1$ is due to V. V. Arestov, 1981. Also {\bf 12.5}b)
is apparently true with $\g(p)=1.$ An inequality similar to {\bf 12.5}a)
is proved in $L_p(\BR^m)$ for any $m.$

Methods of summability of Fourier integrals are often considered on
place of the summability of the \Fs. The question is when in some sense
$$\Il_{-\iy}^\iy\vp(\ve t)\hf(t)e^{itx}\,dt\to f(x)$$ as $\ve\to+0,$
where $\vp$ is a fixed function. If a (multiplier) function $\vp$ is
of compact support, then the integral is EFET.

For all $p>0$ the multiplier can be defined first on $L_p\cap L_1(\BR)$
which is dense in $L_p(\BR)$ and then to extend to $L_p(\BR)$ by
continuity. The comparison principle may be formulated as well. Note
only that for $p\ge1$ there exists connection in both directions
between summability methods (multipliers) of the \Fs\ and the Fourier
integrals generated by the same function $\vp$ (K. de Leeuw, 1965;
see [M22]).

Let us consider now multipliers in the Hardy spaces $H_p$ on the halfplane
for $p\in(0,1].$

Let $f$ be analytic in the upper halfplane $\{z:\Im z>0\}.$ We say
that $f\in H_p$ if $$||f||_{H_p}^p=\sup\limits_{y>0}\Il_{-\iy}^\iy
|f(x+iy)|^p\,dx<\iy.$$
There exists a limit function, as $y\to+0,$ which is in $L_p(\BR).$
A bounded measurable function $\vp:[0,+\iy)\to\BC$ belongs to $M_p$
if it defines a continuous multiplier taking $H_p$ into $H_p.$
\proclaim{12.6} If $\vp\in C(\BR),$ $\vp(x)=0$ for $|x|\ge\s,$ and
$\hat\vp\in L_p(\BR),$ then $\vp\in M_p$ and for $p\in(0,1]$
$$||\vp||_{M_p}\le\g(p)\s^{1/p-1}||\hat\vp||_p.$$  \ep

Sufficient conditions similar to given in {\bf 6.5} are obtained too,
and also in the multiple case. Weaker conditions of exponential decrease
of $\vp\in C^\iy(\BR)$ and its derivatives were obtained earlier by
A. A. Solyanik, 1986.

Also an analog of the theorem {\bf 6.4} on approximation of the Fourier
integrals of functions in $H_p$ on the halfplane by the Bochner-Riesz means
takes place. The following statement is similar to {\bf 9.8}.
\proclaim{12.7} For every $r\in\BN,$ for every $f\in W_\iy^r[0,+\iy),$
and for every $\s>0$ there exists an entire function $g_\s$ of halfdegree
not higher than $\s$ satisfying for every $x\in[0,+\iy)$ the inequality $$|
f(x)-g_\s(x)|\le 2^r K_r(\sqrt{x}/\s)^r+\g(r)(\sqrt{x})^{r-1}/\s^{r+1}.$$\ep

A general direct theorem, without exact constants in the main term, on
approximation on the halfaxis by entire functions of given halfdegree
as well as the inverse to it was obtained by Yu. A. Brudnyi, 1959
(see [M26]).
\flushpar
{\bf Bibliographical remarks.}
\flushpar
For {\bf 12.1}, see [T15,18]; {\bf 12.2-12.3} are due to G. Z. Ber, see
[Ber] and {\bf 12.4} is due to V. V. Zhuk, Dokl.AN SSSR, 1967, in the case
of nonlinear operators and due to G. Z. Ber [Ber] in the case of linear
operators.
 {\bf 12.6} and further, {\bf 12.7} are due to A. V. Tovstolis,
see [To], [TT].

\bigskip
\flushpar
{\bf Concluding remarks.}
\flushpar
In conclusion, let us note that in these twelve sections some results
are even not mentioned, such as those on approximation by Riesz means
of the \Fs\ with respect to eigenfunctions of the Sturm-Liouville
operator (see [T15,18]), on the Hankel transform (see [T15,12]), on
$\ve$-entropy (see [Be15]), on best approximation by polynomials in
$L_2(\BR)$ and Riemann's hypothesis on zeta function (Val. Volchkov,
Ukr. Math. j.,1995), on equivalence of differential operators in
the complex and $p$-adic analysis (Vit. Volchkov), etc.

\head 13. Unsolved problems.  \endhead
\flushpar
{\bf 1.} A. Beurling, 1949, proved that if $f\in C(\BR),$ $\Ll_{|x|\to\iy}
f(x)=0,$ $g\in A^*(\BR)$ (for definition, see before {\bf 1.13}), and
for any $x, h\in\BR$ we have $$|f(x)-f(x+h)|\le|g(x)-g(x+h)|,$$ then
$f\in A(\BR)$ (a new proof is given in [T17] for the multiple case at once).
This yields that if $\vp\in\Lip1,$ $\vp(0)=0,$ $g\in A^*(\BR),$ then
$\vp \circ f\in A(\BR).$

a) Is the converse statement true: if $\vp( g)\in A(\BR)$ for any
$g\in A^*(\BR),$ then $\vp\in\Lip1?$

b) Does the following analog of Beurling's theorem hold for summability
of \Fs\ (see {\bf 3.4} and {\bf 3.2}): if $\lim\lm_{n,k}=1$ as $n\to\iy$
for every $k,$ then we have for every $f\in L_1(\BT)$ at each of its
Lebesgue points $$\sum\mu_{n,k}c_k(f)e^{ikx}\to f(x)$$ as $n\to\iy$ and
$|\lm_{n,k}-\lm_{n,k+s}|\le|\mu_{n,k}-\mu_{n,k+s}|,$ then for every
$f\in C(\BT)$ and for every $x\in\BT$ we have
$$\sum\lm_{n,k}c_k(f)e^{ikx}\to f(x)$$ as $n\to\iy?$
\bigskip
\flushpar
{\bf 2.} S. M. Nikolskii, 1946, has found the following formula for the norm
of a multiplier in $M_1(S)$ for $S=\BZ_0=\BZ\setminus\{0\}$ in the compact
case of $\Lm f=f*g\sim\sum\lm_k c_k(f)e_k:$
$$||\{\lm_k\}||_{M_1(\BZ_0)}=(1/2)\max\limits_x\Il_{-\pi}^\pi
|g(x+t)-g(t)|\,dt,$$ and by choice of the kernel $g$ shown that the strict
inequality $||\{\lm_k\}||_{M_1(\BZ_0)}<||\{\lm_k\}||_{M(\BZ_0)}.$
This yields immediately that in $L_1$ unlike in $C$ not every multiplier
is norm preserve extendible even in the case of the spectrum $\BZ_0,$
the "closest" to $\BZ$). A norm extension coefficient is connected in
this case with the exact Jackson inequality in $L_1(\BT)$ (see {\bf 2.4}
and {\bf 2.5}b)). In the same time
for any $S\sbt\BZ$ every multiplier in $M_1(S)$ is extendible to a
multiplier in $M_1$ (see [M7]). For $S=[-n,n]$ apparently
there exist multipliers norms of which after extension enlarge to
$C\ln n$ times. For which $p\ne2$ any multiplier in $M_p(S)$ is
extendible with preservation of norm or at least with the $\g(p)$ growth of
norm?
\bigskip
\flushpar
{\bf 3.} a) For which sequences $\{\nu_k\}_0^\iy$ (see {\bf 3.6} we have
$$\Ll_{n\to\iy}{1\over n+1}\Sl_{k=0}^n S_{\nu_k}(f;x)=f(x)$$ for every
$f\in L_1(\BT)$ at all its Lebesgue points? This problem posed by Z.
Zalcwasser in 1936 is still open. A similar question may be asked for
$\g$-points introduced by O. D. Gabisonia (Math. Notes, 1973).

b) E. Landau, 1929, proved that for every $f\in H_\iy=L_\iy(\BT,\BZ_+)$
$${1\over n+1}\Sl_{k=0}^n|S_k(f;x)|\le||f||_\iy.$$ Is is true that
for every $f\in L_\iy(\BT)$ the same inequality holds? I know this question
from V. V. Zhuk. My conjecture is that the answer is negative.

c) For which $\dl\ge0$ and $p>0$ we have for every $f\in L_\iy(\BT^2)$
$${1\over n+1}\Sl_{k=0}^n|S_k^\dl(f;x)|^p\le\g(p,\dl)||f||^p_\iy,$$
spherical Bochner-Riesz means are involved? For $p=1$ and $\dl=0$
(spherical partial sums) this inequality is apparently not true; cf. [Ku13].

d) All these questions have sense for functions in $H_p.$
\bigskip
\flushpar
{\bf 4.} a) To construct a modulus of continuity $\om^*$ such that
for every $f\in C(\BT)$
$$||{1\over n+1}\Sl_{k=0}^n|f(\cdot)-S_k(f;(\cdot))|||\asymp\om^*(f;\pi/n).$$

b) Which modulus of smoothness is defined by moxed derivative; see
{\bf 5.6, 5.9,} and {\bf 5.12}?

c) Using the multiplier comparison principle, to prove that for
Fourier-Chebyshev series on $[-1,1]$    $$|f(x)-\Sl_{k=0}^n(1-k^2/n^2)
c_k\cos k\arccos x|\le\g\om(f;\sqrt{1-x^2}/n+n^{-2}).$$
To build in $C[-1,1]$ a sequence $\{p_n(f)\}$ of linear polynomial
operators such that $$|f(x)-p_n(f;x)|\asymp\om_r(f;\sqrt{1-x^2}/n+n^{-2}).$$

d) To formulate a comparison principle for Faber series, in accordance with
position of a point.
\bigskip
\flushpar
{\bf 5.} a) To consider the problem of extension of multipliers with
spectrum, in the Hardy space $H_p$ for $p>0.$

b) For which $\a$ and $\be$ the function $\vp(x)=x^\a\sin x^{-\be}$ belongs
to $M_p$ in $H_p(D)$ near the origin and near the infinity? What are
sufficient conditions of membership $\vp\in M_p,$ in the Hardy space
for $p\in(0,1]$ near the origin and near the infinity (individually,
on the class)?

c) The function $\vp(x)=x_1(\Sl_{j=1}^m x_j)^{-1}\in M_1$ in the Hardy space
on the ball in $\BC^m$ (see [M18], 6.6.3). To point out sufficient
conditions of membership to $M_p$ which are satisfied by this function
for $p=1.$
\bigskip
\flushpar
{\bf 6.} Let $f$ be a function which is even, 
increasing, and convex upwards on $[0,\pi].$ To find an exact order of
increase of $e_n(f)_\iy.$
\bigskip
\flushpar
{\bf 7.} a) Let $f^{(n-1)}(x)>0$ on $[-1,1].$ What is the value
$$E_n(f)_1=\inl_{p_n}\Il_{-1}^1|f(x)-p_n(x)|\,dx$$ for $n\ge 2?$ For
$n=0$ and $n=1$ we have $E_0(f)_1\le(1/2)\Il_0^1\om(f;t)_\iy dt$ and
$E_1(f)_1\le\Il_{1/4}^{1/2}\om_2(f;t)_\iy dt,$ respectively.

b) To find an analog of M\"untz' criterion, the completeness of the system
of powers in $C[0,1],$ in $C[0,1]^m$ for $m\ge2,$ in $L_p$ with $p\in(0,1).$

c) To prove general direct theorems on approximation by polynomials $Q_n$
with integral coefficients on sets in $\BC$ (see \S9). What is the rate
of approximation of the constant by polynomials $Q_n$ on sets in $\BC?$
What is the value $$q(a,b)=\Ll_{n\to\iy}\inl_{Q_n\not\equiv0}
\max\limits_{[a,b]}|Q_n(x)|^{1/n}$$ where $b-a<4?$ For discussion of
these questions and their connection with distribution of prime
numbers, see [T5].
\vfill\eject
\centerline{\bf References}


\Refs\nofrills{Books:}
\widestnumber\key{M27}

\ref\key M1\by N. I. Akhiezer\book Classical moment problem\publ
Fizmatgiz\publaddr Moscow\yr 1961\lang Russian\endref

\ref\key M2\bysame\book Lectures in Approximation Theory, 2nd ed.\publ Nauka
\publaddr Moscow\yr 1965\lang Russian\transl\nofrills German transl. 
\book Vorlesungen \"uber Approximationstheorie\publ Akademie Verlag
\publaddr Berlin\yr 1967\endref

\ref \key M3 \by N. K. Bary \book Trigonometricheskie ryady
\publ Fizmatgiz\publaddr Moscow\yr 1961\transl\nofrills English transl.
\book A treatise on trigonometric series,
Vol. I, II  \publ Pergamon Press  \publaddr New York  \yr 1964  \endref

\ref\key M4\by J. Bergh, J. L\"ofstr\"m\book Interpolation Spaces.
An introduction.\publ Springer-Verlag\publaddr Berlin\yr 1976\endref

\ref\key M5  \by K. M. Davis, Y.-C. Chang\book Lectures on
Bochner-Riesz means \publ London Math. Soc. Lecture Note Series 114,
Cambridge Univ. Press \publaddr Cambridge \yr 1987  \endref

\ref\key M6\by V. K. Dzyadyk\book An Introduction to Theory of Uniform
Approximation of Functions by Polynomials\publ Nauka\publaddr Moscow
\yr 1977\lang Russian \endref

\ref\key M7\by R. E. Edwards\book Fourier series, 2nd ed. Vol.1-1979;
Vol.2-1982 \publ Springer-Verlag  \endref
\ref\key M8\by Le Baron O. Ferguson\book Approximation by Polynomials with
Integral Coefficients\publ Math. Surveys, 17, Amer. Math. Soc.\publaddr 
Providence \yr 1980 \endref

\ref\key M9\by B. I. Golubov, A. V. Efimov, V. A. Skvortsov\book
Walsh series and transforms. Theory and applications \publ Nauka
\publaddr Moscow\yr 1987\lang Russian\endref

\ref \key M10 \by J.-P. Kahane \book S\'eries de Fourier absolument
convergentes \publ Springer - Verlag \yr 1970   \endref

\ref\key M11\bysame\book Some random series of functions\publ Heath
\publaddr Lexington, Mass. \yr 1968\endref

\ref\key M12\by B. S. Kashin, A. A. Saakyan\book Orthogonal Series
\publ Nauka\publaddr Moscow\yr 1984\lang Russian\transl\nofrills
English transl. \publ Transl. Math. Monographs, 75; AMS
\publaddr Providence\yr 1989 \endref

\ref\key M13\by N. P. Korneichuk\book Exact Constants in Approximation
Theory\publ Nauka\publaddr Moscow\yr 1987\lang Russian\transl
\nofrills English transl. \publ Encycl. Math. 38; Cambridge Univ.
Press\publaddr Cambridge\yr 1991\endref

\ref\key M14\by B. Y. Levin\book Entire Functions (course of lectures)
\publ Moscow Univ. Press\publaddr Moscow\yr 1971\lang Russian\endref

\ref\key M15\by E. Lucacs\book Characteristic functions\publ Griffin
\publaddr London\yr 1970\endref

\ref \key M16\by S. M. Nikolskii \book Approximation of Functions
of Several Variables and Imbedding Theorems, 2nd ed.\publ Nauka\publaddr
Moscow\yr 1977\lang Russian\transl\nofrills English translation of the
1st ed. \publ Springer - Verlag
\publaddr Berlin-Heidelberg-New York  \yr 1975   \endref

\ref \key M17\by W. Rudin \book Function theory in polydiscs\publ
Benjamin\publaddr New York\yr 1969\endref

\ref\key M18\bysame\book Function Theory in the Unit Ball of $\Bbb C^n$
\publ Grund. der Math. Wiss. 241; Springer-Verlag\publaddr 
New York\yr 1980\endref

\ref\key M19\by I. A. Shevchuk\book Approximation by Polynomials and
Traces of Functions Continuous on Interval\publ Naukova Dumka
\publaddr Kiev\yr 1992\lang Russian\endref
 
\ref\key M20\by P. K. Suetin\book Classical orthogonal polynomials
\publ Nauka\publaddr Moscow\yr 1979\lang Russian\endref

\ref \key M21 \by  E. M. Stein\book Singular Integrals and Differentiability
Properties of Functions \publ Princeton Univ. Press
\publaddr Princeton, N. J.   \yr 1970     \endref

\ref \key M22\by E. M. Stein, G. Weiss \book Introduction to Fourier
Analysis on Euclidean Spaces  \publ Princeton Univ. Presss
\publaddr Princeton, N. J.  \yr 1971     \endref

\ref\key M23\by A. I. Stepanets\book Uniform Approximations by
Trigonometric Polynomials\publ Naukova Dumka\publaddr Kiev\yr 1981
\lang Russian\endref

\ref\key M24\bysame\book Classification and Approximation of Periodic
Functions\publ Naukova Dumka\publaddr Kiev\yr 1987\lang Russian\endref

\ref\key M25\by V. N. Temlyakov\book Approximation of Functions with
Bounded Mixed Derivative\publ Vol. 178, Trudy Mat. Inst. im. V. A.
Steklova\publaddr Moscow\yr 1986\lang Russian\transl\nofrills English
translation in\book Proc. of the Steklov Inst. of Math.\yr 1989\endref
 
\ref\key M26\by A. F. Timan\book Theory of Approximation of Functions
of a Real Variable\publ Fizmatgiz \publaddr Moscow\yr 1960
\lang Russian\transl\nofrills English transl. \publ Oxford\yr 1963\endref

\ref\key M27\by V. M. Tikhomirov\book Certain Questions of Approximation
Theory\publ Moscow Univ. Press\publaddr Moscow\yr 1976\lang Russian\endref
\ref \key M28\by N. Wiener \book The Fourier Integral and Certain of its
Applications \publ Camb. Univ. Press \publaddr Cambridge \yr 1933    \endref

\ref\key M29\by W. Trebels\book Multipliers for $(C,\a)$-Bounded Fourier 
Expansions in Banach Spaces and Approximation Theory\publ Lecture Notes
Math. 329; Springer\publaddr Berlin\yr 1973 \endref

\ref\key M30\by R. Varga\book Functional Analysis and Approximation
Theory in Numerical Analysis\publ SIAM\publaddr Philadelphia\yr 1971 \endref
 
\ref\key M31\by V. V. Zhuk\book Approximation of Periodic Functions
\publ Leningrad Univ. Press\publaddr Leningrad\yr 1982\lang Russian\endref
\ref \key Z32 \by A. Zygmund \book Trigonometric series, Vol.I, II
\publ Camb. Univ. Press \publaddr Cambridge     \yr 1959    \endref
\endRefs

\Refs\nofrills{Surveys:}
\widestnumber\key{S16}

\ref \key S1 \by S. A. Alimov, R. R. Ashurov, A. K. Pulatov \paper Multiple
Fourier series and integrals \inbook Itogy Nauki i Techniki \publaddr Moscow
\publ VINITI \vol 42 \yr 1989 \pages 7 -- 104 \lang Russian 
\transl\nofrills English translation in\ \ \by V. P. Khavin, N. K. Nikolskii
(Eds.)\book Commutative Harmonic Analysis IV\publ Encycl. Math. Sciences,
Vol.42; Springer-Verlag\publaddr New York\yr 1992\pages 1--95\endref
\ref\key S2\by K. I. Babenko\paper Some Questions of Approximation Theory
and Numerical Analysis\jour Uspekhi Mat. Nauk\yr 1985\vol 40:1\pages 3--27
\lang Russian\transl\nofrills English translation in\jour Russian Math.
Surveys\vol 40:1\yr 1985\pages  \endref
\ref\key S3\by Y. A. Brudnyi\paper A. F. Timan's results in polynomial
approximation of functions\inbook Materials of All-Union Conference in
the Theory of Approximation of Functions. 26-29 of June, 1990\publaddr
Dnepropetrovsk\yr 1991\pages 13--17\lang Russian\endref
\ref\key S4 \by M. I. Dyachenko\paper Some problems in the theory
of multiple trigonometric series\jour Uspekhi Mat. Nauk\vol 47:5
\yr 1992\pages 97--162\lang Russian\transl\nofrills Ehglish translation in
\jour Russian Math. Surveys\vol 47:5\yr 1992\pages 103--171\endref
\ref\key S5\by A. O. Gelfond\paper On uniform approximation by polynomials
with integer coefficients (1955)\inbook Selected papers\publ Nauka
\publaddr Moscow\yr 1973\pages 287--309\lang Russian\endref

\ref \key S6 \by B. I. Golubov \paper Multiple Fourier series and integrals
\inbook Itogy Nauki i Techniki, Matematicheskii Analiz,\publ VINITI \publaddr
Moscow \vol 19 \yr 1982 \pages 3 -- 54 \lang Russian \transl\nofrills English
translation in \jour J. of Soviet Math. \vol 24 \issue 6
\yr 1984 \pages 639 -- 673   \endref

\ref\key S7\by R. S. Ismagilov\paper Widths of sets in linear normed
spaces and approximation of functions by trigonometric polynomials
\jour Uspekhi Mat. Nauk\vol 29:3\yr 1974\pages 161--178
\lang Russian\transl\nofrills Ehglish translation in
\jour Russian Math. Surveys\vol 29:3\yr 1974\pages 103--171\endref

\ref\key S8\by V. G. Krotov, E. A. Storozhenko\paper Approximate and
differential-difference properties of functions from Hardy spaces
\inbook Theory of functions and approximations, Part 1\publ Proc.
of Saratov winter school, Jan.25--Febr.5,1982\publaddr Saratov Univ.
Press\yr 1983\pages 65--80\lang Russian\endref
\ref\key S9\by J. Misiewicz, C. Scheffer\paper Pseudo isotropic measures
\jour Nieuw Archief voor Wiskunde\vol 8/2\yr 1990\pages 111--152\endref
\ref\key S10\by I. Netuka, J. Vesely\paper Mean value property and
harmonic functions\paperinfo Preprint\yr 1994\endref
\ref\key S11\by A. M. Olevsky\paper Modifications of functions and
Fourier series \jour Uspekhi Mat. Nauk\vol 40:3\yr 1985\pages 157--193
\lang Russian\transl\nofrills Ehglish translation in
\jour Russian Math. Surveys\vol 40:3\yr 1985\pages 103--171\endref

\ref\key S12\by M. K. Potapov\paper Approximation by polynomials on
a finite interval of the real axis\inbook Proc. of the Intern. Conference
in Constructive Theory of Functions. 1981\publaddr Sofia\yr 1983
\pages 134--138\lang Russian\endref
\ref\key S13\by V. V. Peller\paper Description of Hankel operators of
the class $G_p$ for $p>0,$ investigation of speed of rational approximation
and other applications\jour Mat. Sbornik\vol 122\issue 4\yr 1983
\pages 481--510\lang Russian\transl\nofrills English translation in
\jour Soviet Mat. Sbornik\vol 122\issue 4\yr 1983\pages  \endref
\ref\key S14\by S. A. Telyakovskii\paper An estimate, of the norm of a
function via the Fourier coefficients of the function, convenient in
problems of Approximation Theory\inbook Trudy Mat. Inst. im. V. A. Steklova
\vol 109\yr 1971\pages 645--97\publ Nauka\publaddr Moscow\lang Russian
\transl\nofrills English translation in \book Proc. of the Steklov Inst.
of Math. \vol 109\pages   \endref
\ref\key S15\bysame\paper On approximation of functions of high smoothness
by Fourier sums\jour Ukrain. Mat. Z.\yr 1989\vol 41\issue 4\pages 510--517
\lang Russian\transl\nofrills English translation in Ukrainian Math. J.
\yr 1989\vol 41\issue 4\pages  \endref
\ref\key S16\by L. Zalcman\paper A bibliografic survey of the Pompeiu
problem\inbook Approximation by Solutions of Partial Differential
Equations\publ Kluwer Acad. Publ.\yr 1992\pages 185--194\endref

\endRefs

\Refs\nofrills{Papers}
\widestnumber\key{VaV14}
\ref\key T1\by R. M. Trigub\paper Aproximation of function by polynomials
with integral coefficients\jour Izv. Akad. Nauk SSSR, ser. matem.
\vol 26 \yr 1962 \issue 2 \pages 261--280 \lang Russian\endref
\ref\key T2\bysame\paper Constructive characteristics of some function
classes \jour Izv. Akad. Nauk SSSR, ser. matem.
\vol 29 \yr 1965 \issue 3 \pages 615--630 \lang Russian\endref
\ref \key T3\bysame \paper Linear summation methods and absolute
convergence of Fourier series \jour Izv. Akad. Nauk SSSR, ser. matem.
\vol 32 \yr 1968 \issue 1 \pages 24--49 \lang Russian \endref
\ref\key T4\bysame\paper Summability and absolute convergence of Fourier
series in whole\inbook Metric Questions in the Theory of Functions and
Mappings \publ Nauk. dumka\publaddr Kiev\yr 1971
\pages 173--266\lang Russian\endref
\ref\key T5\bysame\paper Approximation of functions with Diophantine
conditions by polynomials with integral coefficients
\inbook Metric Questions of the Theory of Functions and Mappings
\publ Nauk. dumka\publaddr Kiev\yr 1971\pages 267--333\lang Russian\endref
\ref\key T6\bysame\paper On linear methods of summability of Fourier series
and moduli of continuity of different orders\jour Sib. Mat. Zh.\yr 1971
\vol 12\issue 6\pages 1416--1421\lang Russian \endref
\ref\key T7\bysame\paper Characteristics of Lipschitz classes of integer
order on the interval by the rate of polynomial approximation\inbook
Theory of functions, functional analysis and their applications
\publaddr Kharkov\vol 18\yr 1973\pages 63--70\lang Russian\endref
\ref\key T8\bysame\paper Integrability of the Fourier transform of a
boundedly supported function\inbook
Theory of functions, functional analysis and their applications
\publaddr Kharkov\vol 23\yr 1975\pages 124--131\lang Russian\endref
\ref \key T9\bysame\paper Summability of Fourier series at Lebesgue
points and one Banach algebra \inbook Metric Questions in the Theory of
Functions and Mappings \publ Nauk. dumka \publaddr Kiev  \vol 6   \yr 1975
\pages 125--135  \lang Russian   \endref

\ref \key T10 \bysame\paper On integral norms for polynomials
\jour Mat. Sbornik \yr 1976 \vol 101(143) \issue 3 \pages 315--333
\lang Russian \transl Engl. transl. in
\jour Math. USSR Sbornik \yr 1976 \vol 30 \issue 3 \pages 279--295 \endref
\ref\key T11\bysame\paper On approximation of functions by polynomials with
special coefficients\jour Izv. VUZ, matem.\issue 1\yr 1977\pages 93--99
\lang Russian\endref
\ref \key T12\bysame\paper Linear methods for the summation of simple and
multiple Fourier series and their approximative properties \inbook Theory of
Approximation of Functions ( Proc. Internat. Conf., Kaluga, 1975 )\publ Nauka
\publaddr Moscow \yr 1977   \pages 383--390     \lang Russian   \endref

\ref \key T13\bysame\paper Integrability and asymptotic behavior
of Fourier transform of radial function \inbook Metric Questions
of the Theory of Functions and Mappings \publ Nauk. dumka
\publaddr Kiev \yr 1977 \pages 142--163 \lang Russian  \endref

\ref \key T14\bysame \paper Summability of multiple Fourier series
\inbook Investigations in the theory of functions of many real variables
\publaddr Yaroslavl \vol 2 \yr 1978 \pages 196--214 \lang Russian \endref

\ref\key T15\bysame\paper Comparison principle and some questions of
the theory of approximation of functions\inbook Theory of Functions and
Mappings\yr 1979\publ Nauk. dumka\publaddr Kiev\pages 149--173
\lang Russian \endref
\ref\key T16\bysame\paper Summability of multiple Fourier series. Growth of
Lebesgue constants \jour Analysis Math. \yr 1980 \vol 6
\issue 3 \pages 255 -- 267  \endref

\ref \key T17\bysame \paper Absolute convergence of Fourier integrals,
summability of Fourier series and polynomial approximation of functions on
the torus \jour Izv.Akad. Nauk SSSR, Ser. Mat.  \vol 44 \yr 1980
\pages 1378--1409  \issue 6 \lang Russian \transl Engl.transl. in
\jour Math. USSR Izv. \vol 17 \issue 3 \yr 1981 \pages 567--593  \endref

\ref \key T18\bysame\book Summability of Fourier series and certain questions
of approximation theory   \publ All-Union Institute of Scientific
and Technical Information, No. 5145--80 \yr 1980      \lang Russian   \endref

\ref \key T19\bysame \paper Absolute convergence of Fourier integrals and
approximation of functions by the linear means of their Fourier series
\inbook Constructive Function Theory'81   \publaddr Sofia
\yr  1983   \pages  178--180      \lang Russian   \endref

\ref\key T20\bysame\paper On the comparison principle for Fourier expansions
\inbook Proc. of the Seminar of the I. N. Vekua Inst. of Applied Math.
\vol 1\issue 2\yr 1985\pages 439--443\lang Russian \endref
\ref\key T21\bysame\paper Some properties of the \Ft\ of measure
and their applications\inbook Proc. Int. Conf. on the Theory of Approx.
of Functions\publ Nauka\publaddr Moscow\yr 1987\pages 439--443
\lang Russian \endref
\ref\key T22\bysame\paper On the comparison principle of the Fourier
expansions and existence subspaces in integral metric\inbook Trudy Mat. Inst.
im. V. A. Steklova \yr 1987 \vol 180  \pages 151--152 \lang Russian
\transl \nofrills Engl. transl. in \jour Proc. Steklov Math. Inst.
\yr 1989 \vol 180 \pages 176--177   \endref
\ref\key T23\bysame\paper A formula for the $K$-functional of a couple
of spaces of functions of several variables\inbook Investigations in
the theory of functions of many variables\publaddr Yaroslavl
\yr 1988\pages 122--127\lang Russian \endref

\ref \key T24\bysame\paper Multipliers of Fourier series and approximation
of functions by polynomials in spaces $C$ and $L$ \jour Dokl. Akad. Nauk SSSR
\vol 306 \yr 1989 \pages 292--296 \lang Russian \transl Engl. transl. in
\jour Soviet Math. Dokl. \vol 39 \issue 3 \yr 1989 \pages 494--498  \endref
\ref\key T25\bysame\paper Asymptotics of a sequence of the norms of
multipliers of the \Fs\ in the spaces $C$ and $L$\inbook Abstracts of the
All-Union School in Lutsk "Theory of Approximation of Functions"
\publaddr Kiev\publ Math. Inst.\yr 1989\lang Russian\endref
\ref\key T26\bysame\paper Criterion of characteristic function and the
Polya type test for radial functions of several variables\jour
Probability theory and its applications\vol 34\issue 4\yr 1989\pages
805--810\lang Russian \endref
\ref\key T27\bysame\book Summability of the \Fs\ and some questions in
approximation theory\publaddr Kiev\publ DS thesis\yr 1984\lang Russian\endref
\ref\key T28\bysame\paper Approximation of continuous periodic functions
with bounded derivative by polynomials\inbook Theory of Mappings and
Approximation of Functions\publ Nauk. dumka\publaddr Kiev\yr 1989
\pages 185--195\lang Russian\endref
\ref\key T29\bysame\paper Multipliers of the \Fs\ \inbook Abstracts
of the Intern. Conf. in Constr. Theory of Functions\publaddr Varna
\yr 1991\lang Russian \endref
\ref\key T30\bysame\paper Multipliers of the Fourier series\jour
Ukr. Mat. Zh.\yr 1991\issue 12\pages 1686--1693\lang Russian\endref
\ref\key T31\bysame\paper Investigations of A. F. Timan on the \Lc\
and linear methods of summability of the \Fs\ \inbook Materials
of the All-Union Conference on the Theory of Approximation of Functions
\publaddr Dnepropetrovsk \yr 1991 \pages 8--12 \lang Russian    \endref
\ref\key T32\bysame\paper Two-sided estimates of approximation by
polynomials and interpolation of spaces of smooth functions on the torus
\inbook Materials
of the All-Union Conference on the Theory of Approximation of Functions
\publaddr Dnepropetrovsk \yr 1991 \pages 70--71 \lang Russian    \endref
\ref\key T33\bysame\paper Positive definite functions and splines
\inbook Theory of functions and approximations. Proc. of the 5th Saratov
winter school Jan.25-Feb.4, 1990, Part 1\publaddr Saratov\yr 1992
\pages 68--75\lang Russian \endref
\ref\key T34\bysame\paper Direct theorems on approximation by algebraic
polynomials of smooth functions on interval\jour Matem. Zametki\vol 54
\issue 6\yr 1993\pages 113--121\lang Russian\transl Engl. transl. in
\jour Math. Notes\vol 54\yr 1993\issue 5-6\endref
\ref\key T35\bysame\paper Multipliers in the Hardy spaces $H_p(D^m)$
for $p\in(0,1]$ and approximate properties of summability methods
of power series\jour Dokl. AN  Rossii\yr 1994\vol 335\issue 6
\pages 697--699\lang Russian \endref
\ref\key T36\bysame\paper The Multiplicators of Power Series in Hardy Spaces
and Approximation Problems in the Polydisk\inbook Special Semester in
Approx. Theory\publ Tecnion\publaddr Haifa\yr 1994\endref

\quad

See also below [BT1,2], [BLT], [KT], [ZT1,2], [ZgT], [TT1].

\quad

\ref\key Be1\by Belinskii E.S.\paper Summability of multiple Fourier series
at Lebesgue points \inbook Theory of functions, functional analysis and their
applications \vol 23 \publaddr Kharkov \yr 1975 \lang Russian
\endref

\ref \key Be2 \bysame \paper On asymptotic behavior of integral norms
of trigonometric polynomials \inbook Metric Questions in the Theory of 
Functions
and Mappings \publaddr Kiev \publ Naukova dumka \yr 1975 \lang Russian
\endref

\ref \key Be3 \bysame \paper Approximation by Bochner-Riesz means and
spherical modulus of continuty \jour Dokl. AN
Ukraine, Ser. A \issue 7 \yr 1975 \lang Russian\endref
\ref \key Be4 \bysame \paper An application of the Fourier transform
to summability of Fourier series \jour Sib. Mat. Zh. \issue
3 \yr 1977 \lang Russian \transl\nofrills Engl. transl. in \jour Siberian
Math. J. \vol 18 \yr 1977 \issue 3 \pages 353-363\endref

\ref \key Be5 \bysame\paper Behavior of Lebesgue constants of some methods
of summation of multiple Fourier series \inbook Metric Questions in the
Theory of Functions and Mappings \publaddr Kiev \publ Naukova  Dumka \yr 1977 
\lang Russian\endref
\ref \key Be6 \bysame \book An application of the Fourier transform
to summability of Fourier series\publ PhD thesis\publaddr Donetsk
\yr 1977\lang Russian\endref

\ref \key Be7 \bysame \paper On some properties of hyperbolic partial
sums \inbook Theory of Functions and Mappings \publaddr Kiev
\publ Naukova dumka\yr 1979 \lang Russian\endref

\ref \key Be8 \bysame \paper On the growth of Lebesgue constants of partial
sums, generated by some unbounded sets \inbook Theory of Mappings and
Approximation of Functions \publaddr Kiev \publ Naukova dumka \yr 1983
\lang Russian\endref
\ref \key Be9\bysame\paper On the summability of Fourier series with the
method of lacunary arithmetic means \jour Analysis Math. \vol 10 \issue 4
\yr 1984\pages 275--282\endref

\ref \key Be10 \bysame \paper Approximation of periodic functions by
the "floating" system of exponents and trigonometric widths \inbook
Investigations in the Theory of functions of several real variables
\publaddr Yaroslavl \yr 1984 \lang Russian\endref
	    
\ref \key Be11 \bysame\paper Approximation of periodic functions of
several variables by a "floating" system of exponentials and trigonometric
widths\jour Dokl. AN SSSR \vol 284 \issue 6 \yr 1985\lang Russian
\transl\nofrills Engl. transl. in \jour Soviet Math. Dokl.
\vol 32 \yr 1985 \issue 2 \pages 571-574\endref
	    
\ref \key Be12 \bysame\paper Approximation by the "floating" system of
exponents on the classes of smooth periodic functions \jour Mat.
sbornik \vol 132 \issue 1 \yr 1987 \lang Russian \transl\nofrills Engl.
transl. in \jour Math. USSR Sbornik \vol 60 \yr 1988 \issue 1\endref

\ref \key Be13 \bysame \paper Approximation by "floating" system of
exponents on the classes of periodic functions with bounded mixed derivative
\inbook Investigations in the Theory of Functions of several real variables
\publaddr Yaroslavl \yr 1988 \lang Russian\endref
	    
\ref \key Be14 \bysame\paper Lebesgue constants of "step-hyperbolic"
partial sums \inbook Theory of Functions and Mappings \publaddr Kiev \publ
Naukova dumka \yr 1989 \lang Russian\endref

\ref \key Be15 \bysame\paper Approximation of functions of several
variables by trigonometric polynomials with given number of harmonics \jour
Analysis Math. \vol 15 \yr 1989 \pages 67-74\endref

\ref \key Be16 \bysame\paper Two extremal problems for trigonometric
polynomials with given number of harmonics \jour Mat. Zametki \vol 49
\issue 1 \yr 1991 \lang Russian \pages 12-19 \transl\nofrills Engl. transl.
in \jour Math. Notes Acad. Sci. USSR \vol 49 \yr 1991\endref

\ref \key Be17 \bysame\paper Asymptotic charachteristic of classes of
functions with conditions on mixed derivatives (mixed difference) \inbook
Investigations in the theory of functions of several real variables \publaddr
Yaroslavl \yr 1990 \pages 22-37 \lang Russian\endref

\ref \key Be18 \bysame \paper Classes $A^p$ and their constructive
characteristics \inbook Materials of the All Union Conference on the Theory
of Approximation of Functions \publaddr Dnepropetrovsk \yr 1991 \lang Russian
\endref

\ref \key Be19 \bysame\paper On strong summability of periodic functions
and embedding theorems \jour Dokl. Ross. Akad. Nauk \yr 1993 \vol 332
\issue 2\lang Russian \transl\nofrills Engl. transl. in\jour Russian Acad.
Sci. Dokl. Math.\vol 48\yr 1994\issue 2\pages 255-258\endref

\ref \key Be20 \bysame\paper Decomposition theorems and approximation
by "floating" system of exponentials \jour Trans. Amer. Math. Soc.
(accepted)\endref

\quad

\ref \key BL1\by E. S. Belinskii and E. R. Liflyand \paper Lebesgue constants
and integrability of Fourier transform of radial functions \jour Dokl. Acad.
Sci. of Ukraine, Ser. A \yr 1980 \issue 6 \pages 5--10  \lang Russian \endref

\ref \key BL2\bysame \paper On asymptotic behavior of Lebesgue
constants of radial summability methods \inbook Constructive Theory
of Functions and Theory of Mappings \publ Nauk. dumka
\publaddr Kiev  \yr 1981 \pages 49--70 \lang Russian   \endref
\ref\key BL3\bysame\paper Behavior of the Lebesgue Constants of Hyperbolic
Partial Sums\jour Mat. Zametki \vol 43
\issue 2 \yr 1988 \lang Russian \pages 192-196 \transl\nofrills Engl. transl.
in \jour Math. Notes Acad. Sci. USSR \vol 43 \yr 1988\pages 107--109\endref
\ref\key BL4\bysame\paper Approximation properties in $L_p, 0<p<1$
\jour Funct. et Appr.\vol XXII\yr 1994\pages    \endref

\quad

\ref\key BLT\by E. S. Belinskii, E. R. Liflyand, R. M. Trigub
\paper The Banach algebra $A^*$ and its properties\paperinfo accepted
in the J. Fourier Anal. Appl.\endref

\quad

\ref\key BT1\by E. S. Belinskii and R. M. Trigub\paper Some numerical
inequalities and their applications to the theory of summability of Fourier
series \inbook Constructive Theory of Functions and Theory of Mappings
\publ Nauk. dumka \publaddr Kiev \yr 1981 \pages 70--81 \lang Russian \endref

\ref \key BT2\bysame\paper Summability on Lebesgue set and one Banach
algebra \inbook Theory of functions and approximations, Part 2 \publaddr
Saratov \yr  1983  \pages 29--34 \lang Russian   \endref

\quad

\ref \key Ku1\by O. I. Kuznetsova \paper paper On some properties
of polynomial operators of triangular form in a space of continuous
periodic functions of two variables \jour Dokl. Akad. Nauk SSSR
\vol 223 \yr 1975 \issue 6 \pages 1304--1306 \lang Russian
\transl\nofrills Engl. transl. in \jour Soviet Math. Dokl.
\vol 16\yr 1975\issue 4\pages 1080--1083 \endref
\ref\key Ku2\bysame\paper The Asymptotic Behavior of the Lebesgue
Constants for a Sequence of Triangular Partial Sums of Double Fourier Series
\jour Sib. Mat. Zh.\vol XVIII\yr 1977\issue 3\pages 629--636\lang Russian
\transl Engl. transl. in\jour Siberian Math. J.\vol 18\yr 1977
\pages 449--454\endref
\ref\key Ku3\bysame \book Lebesgue constants and approximate properties of
linear means of multiple \Fs\ \publ PhD thesis\yr 1985\publaddr
Donetsk \lang Russian\endref
\ref\key Ku4\bysame\paper On one condition of integrability of multiple
trigonometric series\inbook Proc. of the Seminar of the I. N. Vekua Inst.
of Applied Math.\vol 1\issue 2\yr 1985\pages 87--90\lang Russian \endref
\ref\key Ku5\bysame\paper On asymptotics of approximation of
differentiable functions\inbook Proc. Int. Conf. on the Theory of Approx.
of Functions\publ Nauka\publaddr Moscow\yr 1987\pages 243--245\endref
\ref\key Ku6\bysame\paper Integrability and strong summability of multiple
trigonometric series\inbook Trudy Mat. Inst.
im. V. A. Steklova \yr 1987 \vol 180  \pages 143--144 \lang Russian
\transl \nofrills English translation in \jour Proc. Steklov Math. Inst.
\yr 1989 \vol 180 \pages 168--169   \endref
\ref\key Ku7\bysame\paper Asymptotic approximation of smooth functions
\inbook Theory of Mappings and
Approximation of Functions\publ Nauk. dumka\publaddr Kiev\yr 1989
\pages 75--81\lang Russian\endref
\ref\key Ku8\bysame \paper On the strong Karleman means of multiple
trigonometric \Fs\ \jour Ukr. Mat. Zh.\yr 1992\vol 44\issue 2
\pages 275--279\lang Russian \endref
\ref\key Ku9\bysame\paper Strong summability with gaps of multiple \Fs\
\inbook Proc. of the Seminar of the I. N. Vekua Inst.
of Applied Math.\vol 7\issue 2\yr 1992\pages \lang Russian \endref
\ref\key Ku10\bysame\paper On partial sums with respect to polyhedra of
the \Fs\ of bounded functions\jour Anal. Math.\yr 1993\vol 19\issue 4
\pages 267--272\endref
\ref\key Ku11\bysame\paper On the growth of partial sums by polyhedra of
\Fs\ of bounded functions\publ Abstracts of ICM94\publaddr Z\"urich
\yr 1994\endref
\ref\key Ku12\bysame\paper The asymptotic behavior of the Lebesgue function
of two variables\jour Ukr. Mat. Zh.\yr 1995\vol 47\issue 2\pages
220--226\lang Russian \endref
\ref\key Ku13\bysame\paper To the question of the strong summability
by circles\paperinfo  accepted in Ukr. Mat. Zh.\endref

\quad

\ref\key KT\by O. I. Kuznetsova, R. M. Trigub\paper Two-sided estimates
of approximation of functions by Riesz and Marzinkiewicz means\jour
Dokl. AN SSSR\yr 1980\vol 251\issue 1\pages 34--36\lang Russian
\transl Engl. transl. in\jour Soviet Math. Dokl. \endref

\quad
\ref\key L1\by E. R. Liflyand
\paper On integrability of the Fourier transform of a function of compact
support and summability of Fourier series of functions of two variables
\inbook Metric Questions in the Theory of Functions and Mappings
\publ Nauk. dumka \publaddr Kiev  \vol 6   \yr 1975   \pages 69--81
\lang Russian   \endref

\ref\key L2\bysame\paper Some questions of absolute convergence of multiple
Fourier integrals\inbook Theory of Functions and Mappings
\publ Nauk. dumka \publaddr Kiev\yr 1979   \pages 110--132
\lang Russian   \endref
\ref \key L3\bysame\paper On certain conditions of integrability of the
Fourier transform \jour Ukranian Math. J.  \yr 1980 \vol 32 \issue 1
\pages 110--118 \lang Russian \endref
\ref\key L4\bysame\book Absolute convergence and summability of multiple
Fourier series and Fourier integrals\publ PhD thesis\publaddr Donetsk
\yr 1982\lang Russian\endref
\ref\key L5\bysame Absolute convergence and summability with respect
to spheres of multiple Fourier integrals of weakly discontinuous
functions\inbook Theory of Mappings and Approximation of Functions
\publ Nauk. dumka \publaddr Kiev\yr 1983   \pages 77--88
\lang Russian   \endref
\ref\key L6\bysame  \paper Exact order of the Lebesgue constants
of hyperbolic partial sums of multiple Fourier series \jour Mat. Zametki
\yr 1986\vol 39\issue 5\pages 674 -- 683 \lang Russian \transl
\nofrills Engl. transl. in \jour Math. Notes Acad. Sci. USSR
\yr 1986\vol 39 \issue 5-6 \pages 369 -- 374  \endref
\ref \key L7\bysame  \paper Sharp estimates of the Lebesgue constants of
partial sums of multiple Fourier series \inbook Trudy Mat. Inst.
im. V. A. Steklova \yr 1987 \vol 180  \pages 151--152 \lang Russian
\transl \nofrills Engl. transl. in \jour Proc. Steklov Math. Inst.
\yr 1989 \vol 180 \pages 176--177   \endref

\ref \key L8\bysame  \paper Order of growth of Lebesgue constants
of hyperbolic means of multiple Fourier series \inbook Materials
of the All-Union Conference on the Theory of Approximation of Functions
\publaddr Dnepropetrovsk \yr 1991 \pages 59--61 \lang Russian    \endref

\ref \key L9\bysame  \paper On asymptotics of Fourier transform for
functions of certain classes \jour Anal. Math. \yr 1993
\vol 19 \issue 2 \pages 151--168     \endref 

\ref \key L10\bysame  \paper On the operator of division by a power
function in  the multi-dimensional case\jour Mat. Zametki
\yr 1994\vol 56\issue 2\pages 82--90 \lang Russian \transl
\nofrills Engl. transl. in \jour Math. Notes Acad. Sci. Russia
\yr 1994\vol 56\issue 1-2 \pages   \endref
\ref\key L11\bysame\paper Estimates of Lebesgue constants via
Fourier transforms. Many dimensions\paperinfo Preprint, Bar-Ilan Univ.
\yr 1995\endref

\quad

See also [BL1-4], [BLT].

\quad

\ref\key N1\by Yu. L. Nosenko\paper On comparison of linear methods of
summability of double Fourier series
\inbook Metric Questions in the Theory of Functions and Mappings
\publ Nauk. dumka \publaddr Kiev  \vol 6   \yr 1975   \pages 89--101
\lang Russian   \endref
\ref\key N2\bysame\paper Approximate properties of the Riesz means
of double Fourier series\jour Ukr. Mat. Zh.\yr 1979\issue 2
\pages 157--165\lang Russian\endref
\ref\key N3\bysame\paper Concerning the Sidon-type conditions for
integrability of double trigonometric series
\inbook Theory of Functions and Mappings
\publ Nauk. dumka \publaddr Kiev    \yr 1979   \pages 132--149
\lang Russian   \endref
\ref\key N4\bysame\paper The exact order of deviation of continuous
functions of two variables from their rectangular Riesz means
\inbook Constructive Theory of Functions and Theory of Mappings
\publ Nauk. dumka \publaddr Kiev\yr 1981   \pages 129--134
\lang Russian   \endref
\ref\key N5\bysame \paper Summability of double Fourier series by
Bernstein-Rogozinski type methods
\inbook Theory of functions and approximations. Proc. of the 1st Saratov
winter school Jan.24-Feb.5, 1982, Part 2\publaddr Saratov\yr 1983
\pages 115--118\lang Russian \endref
\ref\key N6\bysame\book Regularity and approximate properties of linear
methods of summability of double Fourier series
\publ PhD thesis\publaddr Donetsk\yr 1983\lang Russian\endref
\ref\key N7\bysame\paper Regularity and approximate properties of
Bernstein-Rogozinski type methods of double Fourier series\jour
Ukr. Mat. Zh.\yr 1985\vol 37\issue 5\pages 599--604\lang Russian\endref
\ref\key N8\bysame\paper Sufficient conditions of integrability of
double trigonometric series
\inbook Theory of functions and approximations. Proc. of the 2nd Saratov
winter school, 1984, Part 3\publaddr Saratov\yr 1986
\pages 55--59\lang Russian \endref
\ref\key N9\bysame\paper Sufficient conditions of integrability of
double cosine trigonometric series\inbook Trudy Mat. Inst.
im. V. A. Steklova \yr 1987 \vol 180  \pages 166--168 \lang Russian
\transl \nofrills Engl. transl. in \jour Proc. Steklov Math. Inst.
\yr 1989 \vol 180 \pages  \endref
\ref\key N10\bysame\paper Regularity of Bernstein-Rogozinski type means of
double Fourier series of continuous functions\inbook Theory of Mappings
and Approximation of Functions\publ Nauk. dumka \publaddr Kiev\yr 1989
\pages 133--142\lang Russian   \endref
\ref\key N11\bysame\paper Approximation of functions by Riesz means of
their Fourier series\inbook Harmonic Analysis and Progress in
Approximate Methods\publ Math. Inst. AN Ukraine\publaddr Kiev
\yr 1989\pages 83--84\lang Russian\endref
\ref\key N12 \bysame \paper Deviation of continuous functions from
some means of their Fourier series\inbook Materials
of the All-Union Conference on the Theory of Approximation of Functions
\publaddr Dnepropetrovsk \yr 1991 \pages 61--62 \lang Russian    \endref
\ref\key N13\bysame\paper On approximation of functions by $(C,\a)$ means
of their \Fs\ \inbook Abstracts of the 6th annual IWAA\publ Maine
\publaddr Orono, USA\yr 1992\endref
\ref\key N14\bysame \paper Approximation of functions by Riesz, $(C,\a)$
and typical means of \Fs\ of these functions\inbook Buletinul Stiintific
al Univ. din Baia mare, ser B\publ Matem. Inform X\publaddr Romania
\yr 1994 \pages 89--92\endref

\quad

See also [ZN]

\quad

\ref\key G1\by V. A. Glukhov\paper Comparison of $(C,\a)$ means and
Nikolskii type summability test for Fourier series with respect to
multiplicative systems\inbook Theory of functions and approximations.
Proc. of the 2nd Saratov winter school, 1984, Part 3\publaddr Saratov\yr 1986
\pages 78--84\lang Russian \endref
\ref\key G2\bysame\paper On summability of multiple Fourier series with
respect to multiplicative systems\jour Mat. Zametki
\yr 1986\vol 39\issue 5\pages 665--673 \lang Russian \transl
\nofrills Engl. transl. in \jour Math. Notes Acad. Sci. USSR
\yr 1986\vol 39 \issue 5-6 \pages 361--368   \endref
\ref\key G3\bysame \paper On summation of Fourier-Walsh series\jour
Ukr. Mat. Zh.\yr 1986\issue 3\pages 303--309\endref
\ref\key G4\bysame\book Some questions of summability of simple and
multiple Fourier series with respect to multiplicative systems
\publ PhD thesis\publaddr Donetsk\yr 1987\lang Russian\endref
\ref\key G5\bysame\paper Strong summability of Fourier series with
respect to periodic multiplicative systems \inbook Theory of Mappings
and Approximation of Functions\publ Nauk. dumka \publaddr Kiev\yr 1989
\pages 40--49\lang Russian   \endref

\ref\key G6\bysame\paper Uniform convergence of Fourier-Walsh series
of composition of functions\inbook Materials
of the All-Union Conference on the Theory of Approximation of Functions
\publaddr Dnepropetrovsk \yr 1991 \pages 36--37 \lang Russian    \endref

\quad

\ref\key Z1\by Zastavny V. P.\paper On the set of zeros of the \Ft\ of
measure and on summability of double \Fs\ by Bernstein-Rogozinski
type methods\jour Ukr. Mat. Zh.\yr 1984\vol 36\issue 5\pages
615--621\lang Russian\endref
\ref\key Z2\bysame\paper A theorem on zeros of the \Ft\ of the
indicator function and its applications\paperinfo Deposited at VINITI
No. 701-B, 1987\lang Russian\endref
\ref\key Z3\bysame The \Ft\ of measure and some questions of approximation
theory\publ PhD thesis\publaddr Donetsk\yr 1987\lang Russian\endref
\ref\key Z4\bysame\paper On the set of zeros of the \Ft\ of the indicator
function of a convex body and on summability at Lebesgue points
\inbook Theory of Mappings and Approximation of Functions\publ Nauk.
dumka \publaddr Kiev\yr 1989\pages 40--49\lang Russian   \endref
\ref\key Z5\bysame\paper On some properties of positive definite splines
\inbook Materials
of the All-Union Conference on the Theory of Approximation of Functions
\publaddr Dnepropetrovsk \yr 1991 \pages 47--48 \lang Russian    \endref
\ref\key Z6\bysame\paper Positive definite functions depending on the norm.
Solution of the Shoenberg problem\paperinfo Preprint of the Inst. Appl.
Math. and Mechanics, Donetsk, 1991\lang Russian\endref
\ref\key Z7\bysame\paper Positive definite functions depending on the norm
\jour Dokl. AN Russia\yr 1992\vol 352\issue 5\pages 901--903
\lang Russian\endref
\ref\key Z8\bysame\paper Positive definite functions depending on the norm
\jour Russian J. Math. Phisics\yr 1993\endref

\quad

\ref\key ZN\by V. P. Zastavny, Yu. L. Nosenko\paper Regularity of
Bernstein-Rogozinski type means\inbook Proc. Int. Conf. on the Theory of
Approx. of Functions\publ Nauka\publaddr Moscow\yr 1987\pages 183--184\endref

\quad

\ref\key ZT1\by V. P. Zastavny, R. M. Trigub\paper Positive definite splines
\paperinfo Deposited at Ukr. NIINTI No. 593-Uk., 1987\endref
\ref\key ZT2\bysame\paper On functions with zero integrals over sets
congruent to a given one\inbook Theory of functions and approximations.
Proc. of the 3d Saratov winter school, 1986, Part 2\publaddr Saratov\yr 1988
\pages 69--71\lang Russian \endref

\quad

\ref\key VaV1\by Val. V. Volchkov\paper On functions with zero integrals
over cubes\jour Ukr. Mat. Zh.\yr 1991\vol 43\issue 6\pages 859--863
\lang Russian\endref
\ref\key VaV2\bysame\book Pompeiu type problems on bounded domains
\publ PhD thesis\publaddr Donetsk\yr 1991\lang Russian\endref
\ref\key VaV3\bysame\paper A theorem on spherical means for some
differential equations\jour Dokl. AN Ukraine\issue 5\yr 1992
\pages 8-11\lang Russian\endref
\ref\key VaV4\bysame\paper A uniqueness theorem for multiple lacunary
trigonometric series\jour Mat. Zametki
\yr 1992\vol 51\issue 6\pages 27--31 \lang Russian \transl
\nofrills Engl. transl. in \jour Math. Notes Acad. Sci. USSR
\yr 1992\vol 51 \issue 5-6 \pages    \endref
\ref\key VaV5\bysame\paper Morera type theorems in the domains with
the weak cone condition\jour Izv. VUZ. Matem.\yr 1993\issue 10
\pages 15--20\lang Russian\endref
\ref\key VaV6\bysame\paper On one Zalcman's problem and its generalizations
\jour Mat. Zametki
\yr 1993\vol 53\issue 2\pages 30--36\lang Russian \transl
\nofrills Engl. transl. in \jour Math. Notes Acad. Sci. USSR
\yr 1993\vol 53 \issue 1-2 \pages    \endref
\ref\key VaV7\bysame\paper On the Pompeiu problem and certain of its
applications\jour Dokl. AN Ukraine\yr 1993\issue 11\pages 9--12
\lang Russian\endref
\ref\key VaV8\bysame\paper New mean-value theorems for solution of the
Helmholtz\jour Mat. Sbornik\yr 1993\vol 184\issue 7\pages 71--78
\lang Russian\transl Engl. transl. in\jour Russian Math. Sb.\yr 1993\endref
\ref\key VaV9\bysame\paper Two radii theorems on bounded domains of
Euclidean spaces\jour Dif. Urav.\yr 1994\vol 30\issue 10\pages
1719--1724\lang Russian\endref
\ref\key VaV10\bysame\paper Mean-value theorems for one class of polynomials
\jour Sib. Mat. Zh.\vol 35\issue 4\pages 737--745\lang Russian\transl
Engl. transl. in\jour Siberian Math. J.\yr 1994\vol 18\endref
\ref\key VaV11\bysame\paper Morera type theorems on the unit disk
\jour Analysis Math.\yr 1994\vol 20\pages 49--63\endref
\ref\key VaV12\bysame\paper New two radii theorems in the theory of
harmonic functions\jour Izv. AN Russia, ser. matem.\yr 1994\vol 58
\issue 1\pages 182--194\lang Russian\transl Engl. transl. in\jour
Russian Math. Izv.\yr 1994\endref
\ref\key VaV13\bysame\paper On the Pompeiu problem and certain of its
generalizations\jour Ukr. Mat. Zh.\yr 1994\vol 47\issue 10\pages
1444--1448\lang Russian\endref
\ref\key VaV14\bysame\paper Approximation of functions on bounded domains
in $\BR^n$ by linear combination of shifts\jour Dokl. AN Russia\yr 1994
\vol 334\issue 5\pages 560--561\endref

\quad

\ref\key Le1\by V. O. Leontyev\paper The second term in the Kolmogorov
asymptotic formula for approximation by partial sums of the \Fs\
\jour Dokl. AN Ukraine, ser. A\yr 1990\issue 5\pages 17--21
\lang Russian \endref
\ref\key Le2\bysame\paper Extension of a multiplier from the spectrum
\inbook Materials
of the All-Union Conference on the Theory of Approximation of Functions
\publaddr Dnepropetrovsk \yr 1991 \pages 58--59 \lang Russian    \endref
\ref\key Le3\bysame \book Asymptotic approximations of
differentiable functions
by the \Fs\ \publ PhD thesis\publaddr Donetsk\yr 1991\lang Russian\endref
\ref\key Le4\bysame \paper Asymptotic approximations of differentiable
functions by the \Fs\ \jour Dokl. AN Russia\vol 326\yr 1992\issue 1
\lang Russian\endref

\quad

\ref\key Ber\by G. Z. Ber\paper On the interference phenomenon in integral
metric and approximation by entire functions of exponential type\inbook
Theory of functions, functional analysis and their applications
\publaddr Kharkov\vol 34\yr 1980\pages 11--24\lang Russian\endref

\quad

\ref\key Bl1\by A. M. Belenkii\paper On expansion of functions in the
Fourier-Jacobi series\inbook Constructive Theory of Functions and Theory
of Mappings\publ Nauk. dumka \publaddr Kiev\yr 1981   \pages 35--48
\lang Russian   \endref
\ref\key Bl2\bysame\paper On uniform convergence of Fourier-Jacobi series
on the interval of orthogonality  \jour Mat. Zametki
\yr 1989\vol 46\issue 6\pages 18--26\lang Russian \transl
\nofrills Engl. transl. in \jour Math. Notes Acad. Sci. USSR
\yr 1989\vol 46 \issue 5-6 \pages    \endref

\quad

\ref\key ZgT\by N. A. Zagorodnii, R. M. Trigub\paper On one Salem's question
\inbook Theory of Functions and Mappings\publ Nauk. dumka \publaddr
Kiev\yr 1979   \pages 97--101 \lang Russian   \endref

\quad

\ref\key Kl\by E. M. Klebanov\paper Estimates of approximation by linear
means of Fourier series exact on class\paperinfo Deposited in GNTB
Ukraine No. 175-Uk, 1994\lang Russian \endref

\quad

\ref\key ViV\by Vit. V. Volchkov\paper Multipliers of power series in
Hardy spaces in the unit ball in $\BC^n$ \paperinfo Abstracts of
scientific conference of Donetsk Univ., April 1995\pages 186--187
\lang Russian \endref

\quad

\ref\key To\by A. V. Tovstolis\paper Multipliers in the Hardy spaces
$H_p$ in the upper halfspace for $p\in(0,1]$ and approximation by
the Bochner-Riesz means of Fourier integrals \paperinfo Abstracts of
scientific conference of Donetsk Univ., April 1995\pages 189--190
\lang Russian \endref

\quad

\ref\key TT\by A. V. Tovstolis, R. M. Trigub\paper Pointwise approximation
by polynomials on an interval and by entire functions on the exterior of
interval and halfaxis\inbook Theory of approximation and problems of
numerical mathematics\publ Abstracts of Int. Conf., May 26-28, 1993
\publaddr Dnepropetrovsk\yr 1993\endref

\endRefs
\enddocument